\documentclass[12pt]{article}

\usepackage{amssymb}
\usepackage{amsmath}
\usepackage{xcolor}
\catcode`\@=11

\newcount\printerresolution
\global\printerresolution=300

\newcount\headsize
\global\headsize=450

\font\dotfont = lcircle10 at 3pt

\newcount\epihead
\global\epihead=200

\newcount\defaultscale
\def\setdefaultscale#1{\global\defaultscale=#1}
\setdefaultscale{100}

\newcount\actualtextarrowspace
\newcount\actualtextarrowlength

\newcommand{\computetextparameters}%
{\global\actualtextarrowspace=\textarrowlength%
\global\advance\actualtextarrowspace by 3%
\global\actualtextarrowlength=\textarrowlength%
\global\multiply\actualtextarrowlength by 100}

\newcount\textarrowlength
\def\settextarrowlength#1{\global\textarrowlength=#1%
\computetextparameters}
\settextarrowlength{20}

\newcount\actualdisplayarrowspace
\newcount\actualdisplayarrowlength

\newcommand{\computedisplayparameters}%
{\global\actualdisplayarrowspace=\displayarrowlength%
\global\advance\actualdisplayarrowspace by 3%
\global\actualdisplayarrowlength=\displayarrowlength%
\global\multiply\actualdisplayarrowlength by 100}

\newcount\displayarrowlength
\def\setdisplayarrowlength#1{\global\displayarrowlength=#1%
\computedisplayparameters}
\setdisplayarrowlength{40}

\def\@ifnexttok#1#2#3{\let\@tempe #1\def\@tempa{#2}\def\@tempb{#3}%
\futurelet\@tempc\@ifntok}
\def\@ifntok{\ifx \@tempc \@tempe\let\@tempd\@tempa\else\let\@tempd\@tempb\fi%
\@tempd}

\def\DOT{{\dotfont q}}

\newcount\basicstep
\global\divide\basicstep by \printerresolution

\newcount\numberofsteps
\numberofsteps=\headsize
\global\divide\numberofsteps by \basicstep

\newcommand{\makehead}[3]{%
\begin{picture}(0,0)%
\multiput(0,0)(#1,#2){#3}{\DOT}%
\multiput(0,0)(-#2,#1){#3}{\DOT}%
\end{picture}}

\newcount\xstep
\newcount\ystep

\newsavebox{\northhead}
\savebox{\northhead}{%
\xstep=-\basicstep%
\multiply\xstep by 7071%
\divide\xstep by 10000%
\ystep=\xstep%
\makehead{\xstep}{\ystep}{\numberofsteps}}
\newcommand{\nhead}{\usebox{\northhead}}

\newsavebox{\easthead}
\savebox{\easthead}{%
\xstep=-\basicstep%
\multiply\xstep by 7071%
\divide\xstep by 10000%
\ystep=-\xstep%
\makehead{\xstep}{\ystep}{\numberofsteps}}
\newcommand{\ehead}{\usebox{\easthead}}

\newsavebox{\southhead}
\savebox{\southhead}{%
\xstep=\basicstep%
\multiply\xstep by 7071%
\divide\xstep by 10000%
\ystep=\xstep%
\makehead{\xstep}{\ystep}{\numberofsteps}}
\newcommand{\shead}{\usebox{\southhead}}

\newsavebox{\westhead}
\savebox{\westhead}{%
\xstep=\basicstep%
\multiply\xstep by 7071%
\divide\xstep by 10000%
\ystep=-\xstep%
\makehead{\xstep}{\ystep}{\numberofsteps}}
\newcommand{\whead}{\usebox{\westhead}}

\newsavebox{\northwesthead}
\savebox{\northwesthead}{%
\makehead{0}{-\basicstep}{\numberofsteps}}
\newcommand{\nwhead}{\usebox{\northwesthead}}

\newsavebox{\northeasthead}
\savebox{\northeasthead}{%
\makehead{-\basicstep}{0}{\numberofsteps}}
\newcommand{\nehead}{\usebox{\northeasthead}}

\newsavebox{\southwesthead}
\savebox{\southwesthead}{%
\makehead{\basicstep}{0}{\numberofsteps}}
\newcommand{\swhead}{\usebox{\southwesthead}}

\newsavebox{\southeasthead}
\savebox{\southeasthead}{%
\makehead{0}{\basicstep}{\numberofsteps}}
\newcommand{\sehead}{\usebox{\southeasthead}}

\newsavebox{\eastnortheasthead}
\savebox{\eastnortheasthead}{%
\xstep=-\basicstep%
\multiply\xstep by 9486%
\divide\xstep by 10000%
\ystep=\xstep%
\divide\ystep by -3%
\makehead{\xstep}{\ystep}{\numberofsteps}}
\newcommand{\enehead}{\usebox{\eastnortheasthead}}

\newsavebox{\northnortheasthead}
\savebox{\northnortheasthead}{%
\xstep=-\basicstep%
\multiply\xstep by 9486%
\divide\xstep by 10000%
\ystep=\xstep%
\divide\ystep by 3%
\makehead{\xstep}{\ystep}{\numberofsteps}}
\newcommand{\nnehead}{\usebox{\northnortheasthead}}

\newsavebox{\southsouthwesthead}
\savebox{\southsouthwesthead}{%
\xstep=\basicstep%
\multiply\xstep by 9486%
\divide\xstep by 10000%
\ystep=\xstep%
\divide\ystep by 3%
\makehead{\xstep}{\ystep}{\numberofsteps}}
\newcommand{\sswhead}{\usebox{\southsouthwesthead}}

\newsavebox{\westsouthwesthead}
\savebox{\westsouthwesthead}{%
\xstep=\basicstep%
\multiply\xstep by 9486%
\divide\xstep by 10000%
\ystep=\xstep%
\divide\ystep by -3%
\makehead{\xstep}{\ystep}{\numberofsteps}}
\newcommand{\wswhead}{\usebox{\westsouthwesthead}}

\newsavebox{\westnorthwesthead}
\savebox{\westnorthwesthead}{%
\xstep=\basicstep%
\multiply\xstep by 3162%
\divide\xstep by 10000%
\ystep=\xstep%
\multiply\ystep by -3%
\makehead{\xstep}{\ystep}{\numberofsteps}}
\newcommand{\wnwhead}{\usebox{\westnorthwesthead}}

\newsavebox{\eastsoutheasthead}
\savebox{\eastsoutheasthead}{%
\xstep=-\basicstep%
\multiply\xstep by 3162%
\divide\xstep by 10000%
\ystep=\xstep%
\multiply\ystep by -3%
\makehead{\xstep}{\ystep}{\numberofsteps}}
\newcommand{\esehead}{\usebox{\eastsoutheasthead}}

\newsavebox{\northnorthwesthead}
\savebox{\northnorthwesthead}{%
\xstep=-\basicstep%
\multiply\xstep by 3162%
\divide\xstep by 10000%
\ystep=\xstep%
\multiply\ystep by 3%
\makehead{\xstep}{\ystep}{\numberofsteps}}
\newcommand{\nnwhead}{\usebox{\northnorthwesthead}}

\newsavebox{\southsoutheasthead}
\savebox{\southsoutheasthead}{%
\xstep=\basicstep%
\multiply\xstep by 3162%
\divide\xstep by 10000%
\ystep=\xstep%
\multiply\ystep by 3%
\makehead{\xstep}{\ystep}{\numberofsteps}}
\newcommand{\ssehead}{\usebox{\southsoutheasthead}}

\newsavebox{\easteastnortheasthead}
\savebox{\easteastnortheasthead}{%
\xstep=-\basicstep%
\multiply\xstep by 8944%
\divide\xstep by 10000%
\ystep=\xstep%
\divide\ystep by -2%
\makehead{\xstep}{\ystep}{\numberofsteps}}
\newcommand{\eenehead}{\usebox{\easteastnortheasthead}}

\newsavebox{\northnorthnortheasthead}
\savebox{\northnorthnortheasthead}{%
\xstep=-\basicstep%
\multiply\xstep by 8944%
\divide\xstep by 10000%
\ystep=\xstep%
\divide\ystep by 2%
\makehead{\xstep}{\ystep}{\numberofsteps}}
\newcommand{\nnnehead}{\usebox{\northnorthnortheasthead}}

\newsavebox{\southsouthsouthwesthead}
\savebox{\southsouthsouthwesthead}{%
\xstep=\basicstep%
\multiply\xstep by 8944%
\divide\xstep by 10000%
\ystep=\xstep%
\divide\ystep by 2%
\makehead{\xstep}{\ystep}{\numberofsteps}}
\newcommand{\ssswhead}{\usebox{\southsouthsouthwesthead}}

\newsavebox{\westwestsouthwesthead}
\savebox{\westwestsouthwesthead}{%
\xstep=\basicstep%
\multiply\xstep by 8944%
\divide\xstep by 10000%
\ystep=\xstep%
\divide\ystep by -2%
\makehead{\xstep}{\ystep}{\numberofsteps}}
\newcommand{\wwswhead}{\usebox{\westwestsouthwesthead}}

\newsavebox{\westwestnorthwesthead}
\savebox{\westwestnorthwesthead}{%
\xstep=\basicstep%
\multiply\xstep by 4472%
\divide\xstep by 10000%
\ystep=\xstep%
\multiply\ystep by -2%
\makehead{\xstep}{\ystep}{\numberofsteps}}
\newcommand{\wwnwhead}{\usebox{\westwestnorthwesthead}}

\newsavebox{\easteastsoutheasthead}
\savebox{\easteastsoutheasthead}{%
\xstep=-\basicstep%
\multiply\xstep by 4472%
\divide\xstep by 10000%
\ystep=\xstep%
\multiply\ystep by -2%
\makehead{\xstep}{\ystep}{\numberofsteps}}
\newcommand{\eesehead}{\usebox{\easteastsoutheasthead}}

\newsavebox{\northnorthnorthwesthead}
\savebox{\northnorthnorthwesthead}{%
\xstep=-\basicstep%
\multiply\xstep by 4472%
\divide\xstep by 10000%
\ystep=\xstep%
\multiply\ystep by 2%
\makehead{\xstep}{\ystep}{\numberofsteps}}
\newcommand{\nnnwhead}{\usebox{\northnorthnorthwesthead}}

\newsavebox{\southsouthsoutheasthead}
\savebox{\southsouthsoutheasthead}{%
\xstep=\basicstep%
\multiply\xstep by 4472%
\divide\xstep by 10000%
\ystep=\xstep%
\multiply\ystep by 2%
\makehead{\xstep}{\ystep}{\numberofsteps}}
\newcommand{\sssehead}{\usebox{\southsouthsoutheasthead}}

\newsavebox{\northeasteastnortheasthead}
\savebox{\northeasteastnortheasthead}{%
\xstep=-\basicstep%
\multiply\xstep by 9806%
\divide\xstep by 10000%
\ystep=\xstep%
\divide\ystep by -5%
\makehead{\xstep}{\ystep}{\numberofsteps}}
\newcommand{\neenehead}{\usebox{\northeasteastnortheasthead}}

\newsavebox{\northeastnorthnortheasthead}
\savebox{\northeastnorthnortheasthead}{%
\xstep=-\basicstep%
\multiply\xstep by 9806%
\divide\xstep by 10000%
\ystep=\xstep%
\divide\ystep by 5%
\makehead{\xstep}{\ystep}{\numberofsteps}}
\newcommand{\nennehead}{\usebox{\northeastnorthnortheasthead}}

\newsavebox{\southwestsouthsouthwesthead}
\savebox{\southwestsouthsouthwesthead}{%
\xstep=\basicstep%
\multiply\xstep by 9806%
\divide\xstep by 10000%
\ystep=\xstep%
\divide\ystep by 5%
\makehead{\xstep}{\ystep}{\numberofsteps}}
\newcommand{\swsswhead}{\usebox{\southwestsouthsouthwesthead}}

\newsavebox{\southwestwestsouthwesthead}
\savebox{\southwestwestsouthwesthead}{%
\xstep=\basicstep%
\multiply\xstep by 9806%
\divide\xstep by 10000%
\ystep=\xstep%
\divide\ystep by -5%
\makehead{\xstep}{\ystep}{\numberofsteps}}
\newcommand{\swwswhead}{\usebox{\southwestwestsouthwesthead}}

\newsavebox{\northwestwestnorthwesthead}
\savebox{\northwestwestnorthwesthead}{%
\xstep=\basicstep%
\multiply\xstep by 1961%
\divide\xstep by 10000%
\ystep=\xstep%
\multiply\ystep by -5%
\makehead{\xstep}{\ystep}{\numberofsteps}}
\newcommand{\nwwnwhead}{\usebox{\northwestwestnorthwesthead}}

\newsavebox{\southeasteastsoutheasthead}
\savebox{\southeasteastsoutheasthead}{%
\xstep=-\basicstep%
\multiply\xstep by 1961%
\divide\xstep by 10000%
\ystep=\xstep%
\multiply\ystep by -5%
\makehead{\xstep}{\ystep}{\numberofsteps}}
\newcommand{\seesehead}{\usebox{\southeasteastsoutheasthead}}

\newsavebox{\northwestnorthnorthwesthead}
\savebox{\northwestnorthnorthwesthead}{%
\xstep=-\basicstep%
\multiply\xstep by 1961%
\divide\xstep by 10000%
\ystep=\xstep%
\multiply\ystep by 5%
\makehead{\xstep}{\ystep}{\numberofsteps}}
\newcommand{\nwnnwhead}{\usebox{\northwestnorthnorthwesthead}}

\newsavebox{\southeastsouthsoutheasthead}
\savebox{\southeastsouthsoutheasthead}{%
\xstep=\basicstep%
\multiply\xstep by 1961%
\divide\xstep by 10000%
\ystep=\xstep%
\multiply\ystep by 5%
\makehead{\xstep}{\ystep}{\numberofsteps}}
\newcommand{\sessehead}{\usebox{\southeastsouthsoutheasthead}}

\newsavebox{\isomorphismmark}
\newcommand{\isomark}[1]{\savebox{\isomorphismmark}{#1}}
\isomark{$\cong$}
\newcommand{\isosign}{\usebox{\isomorphismmark}}

\newif\ifuserdist
\newsavebox{\distributormark}
\newcommand{\distmark}[1]{\ifx#1\distcircle\userdistfalse\else%
\userdisttrue\savebox{\distributormark}{#1}\fi}
\newsavebox{\distributorcircle}
\savebox{\distributorcircle}{\begin{picture}(0,0)%
\put(0,0){\circle{4}}\end{picture}}
\newcommand{\distsign}{\ifuserdist\usebox{\distributormark}\else%
\usebox{\distributorcircle}\fi}
\userdistfalse

\newcount\monolength
\newcount\epilength
\newcount\bimolength
\newcount\secondmonolength
\newcount\secondepilength

\newcount\monotail%
\global\monotail=\headsize%
\global\multiply\monotail by 7070%
\global\divide\monotail by 10000%

\newcount\truemonotail%
\newcount\trueepihead%

\def\truetail{\truemonotail=\monotail%
\multiply\truemonotail by 100%
\divide\truemonotail by \SCALE}

\def\truehead{\trueepihead=\epihead%
\multiply\trueepihead by 100%
\divide\trueepihead by \SCALE}

\newcount\Monotail%
\global\Monotail=\headsize%
\global\divide\Monotail by 2%

\newcount\Epihead%
\global\Epihead=\epihead%
\global\multiply\Epihead by 7070%
\global\divide\Epihead by 10000%

\newcount\Truemonotail%
\newcount\Trueepihead%

\def\Truetail{\Truemonotail=\Monotail%
\multiply\Truemonotail by 100%
\divide\Truemonotail by \SCALE}%

\def\Truehead{\Trueepihead=\Epihead%
\multiply\Trueepihead by 100%
\divide\Trueepihead by \SCALE}

\newcount\MonoTail%
\global\MonoTail=\headsize%
\global\multiply\MonoTail by 6324%
\global\divide\MonoTail by 10000%

\newcount\EpiHead%
\global\EpiHead=\epihead%
\global\multiply\EpiHead by 8944%
\global\divide\EpiHead by 10000%

\newcount\TrueMonoTail%
\newcount\TrueEpiHead%

\def\TrueTail{\TrueMonoTail=\MonoTail%
\multiply\TrueMonoTail by 100%
\divide\TrueMonoTail by \SCALE}%

\def\TrueHead{\TrueEpiHead=\EpiHead%
\multiply\TrueEpiHead by 100%
\divide\TrueEpiHead by \SCALE}

\newcount\monotaiL%
\global\monotaiL=\headsize%
\global\multiply\monotaiL by 3162%
\global\divide\monotaiL by 10000%

\newcount\epiheaD%
\global\epiheaD=\epihead%
\global\multiply\epiheaD by 4472%
\global\divide\epiheaD by 10000%

\newcount\truemonotaiL%
\newcount\trueepiheaD%

\def\truetaiL{\truemonotaiL=\monotaiL%
\multiply\truemonotaiL by 100%
\divide\truemonotaiL by \SCALE}%

\def\trueheaD{\trueepiheaD=\epiheaD%
\multiply\trueepiheaD by 100%
\divide\trueepiheaD by \SCALE}

\newcount\SCALE%

\newcount\NUMBER%

\newcount\LINE%

\newcount\COLUMN%

\newcount\WIDTH%

\newcount\SOURCE%

\newcount\ARROW%

\newcount\TARGET%

\newcount\ARROWLENGTH%

\newcount\NUMBEROFDOTS%

\newcounter{x}%

\newcounter{y}%

\newcounter{z}%

\newcounter{horizontal}%

\newcounter{vertical}%

\newskip\itemlength%

\newskip\firstitem%

\newskip\seconditem%

\newcommand{\printarrow}{}%

\newcount\X%

\newcount\Y%

\newcount\Z%

\newcommand{\truex}[1]{%
\NUMBER=#1%
\multiply\NUMBER by 100%
\divide\NUMBER by \SCALE%
\setcounter{x}{\NUMBER}}%

\newcommand{\truey}[1]{%
\NUMBER=#1%
\multiply\NUMBER by 100%
\divide\NUMBER by \SCALE%
\setcounter{y}{\NUMBER}}%

\newcommand{\truez}[1]{%
\NUMBER=#1%
\multiply\NUMBER by 100%
\divide\NUMBER by \SCALE%
\setcounter{z}{\NUMBER}}%

\newcommand{\changecounters}[1]{%
\SOURCE=\ARROW%
\ARROW=\TARGET%
\settowidth{\itemlength}{#1}%
\ifdim \itemlength > 2800\unitlength%
\addtolength{\itemlength}{-2800\unitlength}%
\TARGET=\itemlength%
\divide\TARGET by 1310%
\multiply\TARGET by 100%
\divide\TARGET by \SCALE%
\else%
\TARGET=0%
\fi%
\ARROWLENGTH=5000%
\advance\ARROWLENGTH by -\SOURCE%
\advance\ARROWLENGTH by -\TARGET%
\divide\ARROWLENGTH by 100%
\advance\SOURCE by -\TARGET}%

\newcommand{\initialize}[1]{%
\LINE=0%
\COLUMN=0%
\WIDTH=0%
\ARROW=0%
\TARGET=0%
\changecounters{#1}%
\renewcommand{\printarrow}{#1}%
\begin{center}%
\vspace{2pt}%
\begin{picture}(0,0)}%

\newcommand{\DIAGV}[2]{%
\SCALE=#1%
\setlength{\unitlength}{655sp}%
\multiply\unitlength by \SCALE%
\divide\unitlength by 100%
\initialize{\mbox{$#2$}}}%

\newcommand{\n}[1]{%
\changecounters{\mbox{$#1$}}%
\put(\COLUMN,\LINE){\makebox(0,0){\printarrow}}%
\thinlines%
\renewcommand{\printarrow}{\mbox{$#1$}}%
\advance\COLUMN by 4000}%

\newcommand{\nn}[1]{%
\put(\COLUMN,\LINE){\makebox(0,0){\printarrow}}%
\thinlines%
\ifnum \WIDTH < \COLUMN%
\WIDTH=\COLUMN%
\else%
\fi%
\advance\LINE by -4000%
\COLUMN=0%
\ARROW=0%
\TARGET=0%
\changecounters{\mbox{$#1$}}%
\renewcommand{\printarrow}{\mbox{$#1$}}}%

\newcommand{\conclude}{%
\put(\COLUMN,\LINE){\makebox(0,0){\printarrow}}%
\thinlines%
\ifnum \WIDTH < \COLUMN%
\WIDTH=\COLUMN%
\else%
\fi%
\setcounter{horizontal}{\WIDTH}%
\setcounter{vertical}{-\LINE}%
\end{picture}}%

\newcommand{\diag}{%
\conclude%
\raisebox{0pt}[0pt][\value{vertical}\unitlength]{}%
\hspace*{\value{horizontal}\unitlength}%
\vspace{12pt}%
\end{center}%
\setlength{\unitlength}{1pt}%
}%

\newcommand{\diagv}[3]{%
\conclude%
\NUMBER=#1%
\rule{0pt}{\NUMBER pt}%
\hspace*{-#2pt}%
\raisebox{0pt}[0pt][\value{vertical}\unitlength]{}%
\hspace*{\value{horizontal}\unitlength}%6
\NUMBER=#3%
\advance\NUMBER by 12%
\vspace*{\NUMBER pt}%
\end{center}%
\setlength{\unitlength}{1pt}%
}%

\def\movename(#1,#2)#3{%
\hspace{#1pt}%
\raisebox{#2pt}[5pt][2pt]{\raisebox{#2pt}{$#3$}}%
\hspace{-#1pt}}%

\def\movearrow(#1,#2)#3{%
\makebox[0pt]{%
\hspace{#1pt}\hspace{#1pt}%
\raisebox{#2pt}[0pt][0pt]{\raisebox{#2pt}{$#3$}}}}%

\def\movevertex(#1,#2)#3{%
\mbox{\hspace{#1pt}%
\raisebox{#2pt}{\raisebox{#2pt}{$#3$}}%
\hspace{-#1pt}}}%

\newcommand{\CrossLength}[2]{%
\settowidth{\firstitem}{#1}%
\settowidth{\seconditem}{#2}%
\ifdim\firstitem < \seconditem%
\itemlength=\seconditem%
\else%
\itemlength=\firstitem%
\fi%
\divide\itemlength by 2%
\hspace{\itemlength}}%

\def\basicDIAG#1|{\DIAGV{\defaultscale}{#1}\@ifnexttok|{\finishline}{\basicn}}
\def\basicDIAGV[#1]#2|{\DIAGV{#1}{#2}\@ifnexttok|{\finishline}{\basicn}}
\def\basicn#1|{\n{#1}\@ifnexttok|{\finishline}{\basicn}}
\def\basicnn#1|{\nn{#1}\@ifnexttok|{\finishline}{\basicn}}
\def\finishline#1{\@ifnextchar\end{\diag}%
{\@ifnextchar\spacing{\relax}{\basicnn}}}
\def\spacing(#1,#2,#3){\diagv{#1}{#2}{#3}}

\newif\ifcaption%

\newenvironment{diagram}{%
\iffloatdiag\relax\else
\global\def\diagramcaption##1{%
\global\captiontrue%
\global\def\@diagcaption{##1}}%
\global\def\@diagcaption{}\fi%
\@ifnextchar[{\basicDIAGV}{\basicDIAG}}%
{\iffloatdiag\relax\else%
\ifcaption
\begin{center}\mbox{}\@diagcaption\end{center}%
\else\relax\fi\fi\global\captionfalse}

\global\captionfalse

\gdef\@diaglabel{Diagram}
\gdef\diagramlabel#1{\gdef\@diaglabel{#1}}

\newcounter{Diagram}
\def\theDiagram{\@arabic\c@Diagram}
\def\fps@Diagram{tpbh}
\def\ftype@Diagram{1}
\def\ext@Diagram{lof}
\def\fnum@Diagram{\@diaglabel\ \theDiagram}
\def\Diagram{\@float{Diagram}}
\let\endDiagram\end@float
\@namedef{Diagram*}{\@dblfloat{Diagram}}
\@namedef{endDiagram*}{\end@dblfloat}

\newif\iffloatdiag

{\end{diagram}\caption{\@diagcaption}\end{Diagram}%
\global\floatdiagfalse}

{\end{diagram}\caption{\@diagcaption}\end{Diagram}%
\global\floatdiagfalse%
\typeout{*****}%
\typeout{WARNING: Floating diagram forced to be on the page}%
\typeout{*****}
}

{\end{diagram}\caption{\@diagcaption}\end{Diagram}%
\global\floatdiagfalse%
}

\global\floatdiagfalse

\def\dlabel#1{\immediate\write\@auxout{\string\newlabel{#1}%
{{\theChapter.\@arabic\c@Diagram}{\thepage}}}}

\newcommand{\TUP}[1]{\raisebox{0pt}[0pt][3pt]{}#1}

\newcommand{\TDOWN}[1]{\raisebox{0pt}[6pt][0pt]{}#1}

\newcommand{\tlowername}[2]%
{$\stackrel{\makebox[1pt]{#1}}%
{\begin{picture}(0,0)%
\put(0,0){\makebox(0,6)[t]{\makebox[1pt]{$\scriptstyle#2$}}}%
\end{picture}}$}%

\newcommand{\tcase}[1]{%
\setlength{\unitlength}{0.01pt}%
\makebox[\actualtextarrowspace pt]%
{\raisebox{2.5pt}{#1{\actualtextarrowlength}}}%
\setlength{\unitlength}{1pt}}%

\newcommand{\Tcase}[2]{%
\setlength{\unitlength}{0.01pt}%
\makebox[\actualtextarrowspace pt]%
{\raisebox{2.5pt}{$\stackrel{\scriptstyle #2}{#1{\actualtextarrowlength}}$}}%
\setlength{\unitlength}{1pt}}%

\newcommand{\tbicase}[1]{%
\setlength{\unitlength}{0.01pt}%
\makebox[\actualtextarrowspace pt]%
{\raisebox{1pt}{#1{\actualtextarrowlength}}}%
\setlength{\unitlength}{1pt}}%

\newcommand{\Tbicase}[3]{%
\setlength{\unitlength}{0.01pt}%
\makebox[\actualtextarrowspace pt]%
{\raisebox{-1pt}%
{$\stackrel{\scriptstyle #2}%
{\mbox{\tlowername{#1{\actualtextarrowlength}}%
{\scriptstyle #3}}}$}}%
\setlength{\unitlength}{1pt}}%

\newcommand{\ttricase}[1]{%
\setlength{\unitlength}{0.01pt}%
\makebox[\actualtextarrowspace pt]%
{\raisebox{4pt}{#1{\actualtextarrowlength}%
{\rule{0pt}{6pt}}{\rule{0pt}{6pt}}{\rule{0pt}{6pt}}}}%
\setlength{\unitlength}{1pt}}%

\newcommand{\Ttricase}[4]{%
\setlength{\unitlength}{0.01pt}%
\makebox[\actualtextarrowspace pt]%
{\raisebox{4pt}%
{#1{\actualtextarrowlength}{\rule{0pt}{6pt}\scriptstyle #2}%
{\rule{0pt}{6pt}\scriptstyle #3}{\rule{0pt}{6pt}\scriptstyle #4}}}%
\setlength{\unitlength}{1pt}}%

\newcommand{\tmulticase}[1]{%
\setlength{\unitlength}{0.01pt}%
\makebox[\actualtextarrowspace pt]%
{\raisebox{1pt}{#1{\actualtextarrowlength}}}%
\setlength{\unitlength}{1pt}}%

\newcommand{\DUP}[1]{\raisebox{0pt}[0pt][4pt]{}#1}

\newcommand{\DDOWN}[1]{\raisebox{0pt}[9pt][0pt]{}#1}

\newcommand{\dlowername}[2]%
{$\stackrel{\makebox[1pt]{#1}}%
{\begin{picture}(0,0)%
\put(0,0){\makebox(0,6)[t]{\makebox[1pt]{$\textstyle#2$}}}%
\end{picture}}$}%

\newcommand{\dcase}[1]{%
\setlength{\unitlength}{0.01pt}%
\makebox[\actualdisplayarrowspace pt]%
{\raisebox{2.5pt}{#1{\actualdisplayarrowlength}}}%
\setlength{\unitlength}{1pt}}%

\newcommand{\Dcase}[2]{%
\setlength{\unitlength}{0.01pt}%
\makebox[\actualdisplayarrowspace pt]%
{\raisebox{2.5pt}{$\stackrel{\textstyle #2}{#1{\actualdisplayarrowlength}}$}}%
\setlength{\unitlength}{1pt}}%

\newcommand{\dbicase}[1]{%
\setlength{\unitlength}{0.01pt}%
\makebox[\actualdisplayarrowspace pt]%
{\raisebox{1pt}{#1{\actualdisplayarrowlength}}}%
\setlength{\unitlength}{1pt}}%

\newcommand{\Dbicase}[3]{%
\setlength{\unitlength}{0.01pt}%
\makebox[\actualdisplayarrowspace pt]%
{\raisebox{-1pt}%
{$\stackrel{\textstyle #2}%
{\mbox{\tlowername{#1{\actualdisplayarrowlength}}%
{\textstyle #3}}}$}}%
\setlength{\unitlength}{1pt}}%

\newcommand{\dtricase}[1]{%
\setlength{\unitlength}{0.01pt}%
\makebox[\actualdisplayarrowspace pt]%
{\raisebox{4pt}{#1{\actualdisplayarrowlength}%
{\rule{0pt}{6pt}}{\rule{0pt}{6pt}}{\rule{0pt}{6pt}}}}%
\setlength{\unitlength}{1pt}}%

\newcommand{\Dtricase}[4]{%
\setlength{\unitlength}{0.01pt}%
\makebox[\actualdisplayarrowspace pt]%
{\raisebox{4pt}%
{#1{\actualdisplayarrowlength}{\rule{0pt}{6pt}\textstyle #2}%
{\rule{0pt}{6pt}\textstyle #3}{\rule{0pt}{6pt}\textstyle #4}}}%
\setlength{\unitlength}{1pt}}%

\newcommand{\dmulticase}[1]{%
\setlength{\unitlength}{0.01pt}%
\makebox[\actualdisplayarrowspace pt]%
{\raisebox{1pt}{#1{\actualdisplayarrowlength}}}%
\setlength{\unitlength}{1pt}}%

\newcommand{\AR}[1]%
{\begin{picture}(#1,0)%
\put(0,0){\line(1,0){#1}}%
\put(#1,0){\ehead}%
\end{picture}}%

\newcommand{\DIST}[1]%
{\begin{picture}(#1,0)%
\put(0,0){\line(1,0){#1}}%
\put(#1,0){\ehead}%
\NUMBER=#1%
\divide\NUMBER by 2%
\put(\NUMBER,0){\distsign}%
\end{picture}}%

\newcommand{\DOTAR}[1]%
{\NUMBEROFDOTS=#1%
\divide\NUMBEROFDOTS by 300%
\advance\NUMBEROFDOTS by 1%
\begin{picture}(#1,0)%
\multiput(0,0)(300,0){\NUMBEROFDOTS}{\circle*{100}}%
\put(#1,0){\ehead}%
\end{picture}}%

\newcommand{\MONO}[1]%
{\monolength=#1%
\advance\monolength by -\monotail%
\begin{picture}(#1,0)%
\put(\monotail,0){\line(1,0){\monolength}}%
\put(#1,0){\ehead}%
\put(\monotail,0){\ehead}%
\end{picture}}%

\newcommand{\EPI}[1]%
{\epilength=#1%
\advance\epilength by -\epihead%
\begin{picture}(#1,0)(-#1,0)%
\put(-#1,0){\line(1,0){\epilength}}%
\put(-\epihead,0){\ehead}%
\put(0,0){\ehead}%
\end{picture}}%

\newcommand{\BIMO}[1]%
{\monolength=#1%
\advance\monolength by -\monotail%
\epilength=\monolength%
\advance\epilength by -\epihead%
\begin{picture}(#1,0)(-#1,0)%
\put(-\monolength,0){\line(1,0){\epilength}}%
\put(-\monolength,0){\ehead}%
\put(-\epihead,0){\ehead}%
\put(0,0){\ehead}%
\end{picture}}%

\newcommand{\BIAR}[1]%
{\begin{picture}(#1,700)%
\put(0,0){\line(1,0){#1}}%
\put(#1,0){\ehead}%
\put(0,700){\line(1,0){#1}}%
\put(#1,700){\ehead}%
\end{picture}}%

\newcommand{\BIDIST}[1]%
{\begin{picture}(#1,700)%
\put(0,0){\line(1,0){#1}}%
\put(#1,0){\ehead}%
\put(0,700){\line(1,0){#1}}%
\put(#1,700){\ehead}%
\NUMBER=#1%
\divide\NUMBER by 2%
\put(\NUMBER,0){\distsign}%
\put(\NUMBER,700){\distsign}%
\end{picture}}%

\newcommand{\EQL}[1]%
{\begin{picture}(#1,0)%
\put(0,100){\line(1,0){#1}}%
\put(0,-100){\line(1,0){#1}}%
\end{picture}}%

\newcommand{\LLINE}[1]%
{\begin{picture}(#1,0)%
\put(0,0){\line(1,0){#1}}%
\end{picture}}%

\newcommand{\ADJAR}[1]%
{\begin{picture}(#1,700)%
\put(0,0){\line(1,0){#1}}%
\put(#1,0){\ehead}%
\put(#1,700){\line(-1,0){#1}}%
\put(0,700){\whead}
\end{picture}}%

\newcommand{\ADJDIST}[1]%
{\begin{picture}(#1,700)%
\put(0,0){\line(1,0){#1}}%
\put(#1,0){\ehead}%
\put(#1,700){\line(-1,0){#1}}%
\put(0,700){\whead}
\NUMBER=#1%
\divide\NUMBER by 2%
\put(\NUMBER,0){\distsign}%
\put(\NUMBER,700){\distsign}%
\end{picture}}%

\newcommand{\TRIAR}[4]%
{%
\begin{tabular}{c}%
\mbox{$#2$}\\ \AR{#1}\\%
\mbox{$#3$}\\ \AR{#1}\\%
\mbox{$#4$}\\ \AR{#1}%
\end{tabular}%
}%

\newcommand{\TRIADJAR}[4]%
{%
\begin{tabular}{c}%
\mbox{$#2$}\\ \BKAR{#1}\\%
\mbox{$#3$}\\ \AR{#1}\\%
\mbox{$#4$}\\ \BKAR{#1}%
\end{tabular}%
}%

\newcommand{\QUADRIAR}[1]%
{%
\begin{tabular}{c}%
\AR{#1}\\%
\rule{0pt}{3pt}\\ \AR{#1}\\%
\rule{0pt}{3pt}\\ \AR{#1}\\%
\rule{0pt}{3pt}\\ \AR{#1}%
\end{tabular}%
}%

\newcommand{\QUADRIADJAR}[1]%
{%
\begin{tabular}{c}%
\BKAR{#1}\\%
\rule{0pt}{3pt}\\ \AR{#1}\\%
\rule{0pt}{3pt}\\ \BKAR{#1}\\%
\rule{0pt}{3pt}\\ \AR{#1}%
\end{tabular}%
}%

\newcommand{\QUINTIAR}[1]%
{%
\begin{tabular}{c}%
\AR{#1}\\%
\rule{0pt}{2pt}\\ \AR{#1}\\%
\rule{0pt}{2pt}\\ \AR{#1}\\%
\rule{0pt}{2pt}\\ \AR{#1}\\%
\rule{0pt}{2pt}\\ \AR{#1}%
\end{tabular%
}%

\newcommand{\QUINTIADJAR}[1]%
{%
\begin{tabular}{c}%
\AR{#1}\\%
\rule{0pt}{2pt}\\ \BKAR{#1}\\%
\rule{0pt}{2pt}\\ \AR{#1}\\%
\rule{0pt}{2pt}\\ \BKAR{#1}\\%
\rule{0pt}{2pt}\\ \AR{#1}%
\end{tabular}}%
}%

\newcommand{\ar}{\ifinner\tcase{\AR}\else\dcase{\AR}\fi}%

\newcommand{\Ar}[1]{\ifinner\Tcase{\AR}{#1}\else\Dcase{\AR}{#1}\fi}%

\newcommand{\dist}{\ifinner\tcase{\DIST}\else\dcase{\DIST}\fi}%

\newcommand{\Dist}[1]{\ifinner\Tcase{\DIST}{\TUP{#1}}%
\else\Dcase{\DIST}{\TUP{#1}}\fi}%

\newcommand{\dotar}{\ifinner\tcase{\DOTAR}\else\dcase{\DOTAR}\fi}%

\newcommand{\Dotar}[1]{\ifinner\Tcase{\DOTAR}{#1}%
\else\Dcase{\DOTAR}{#1}\fi}%

\newcommand{\mono}{\ifinner\tcase{\MONO}\else\dcase{\MONO}\fi}%

\newcommand{\Mono}[1]{\ifinner\Tcase{\MONO}{#1}\else\Dcase{\MONO}{#1}\fi}%

\newcommand{\epi}{\ifinner\tcase{\EPI}\else\dcase{\EPI}\fi}%

\newcommand{\Epi}[1]{\ifinner\Tcase{\EPI}{#1}\else\Dcase{\EPI}{#1}\fi}%

\newcommand{\bimo}{\ifinner\tcase{\BIMO}\else\dcase{\BIMO}\fi}%

\newcommand{\Bimo}[1]{\ifinner\Tcase{\BIMO}{#1}%
\else\Dcase{\BIMO}{#1}\fi}%

\newcommand{\iso}{\ifinner\Tcase{\AR}{\isosign}\else\Dcase{\AR}{\isosign}\fi}%

\newcommand{\Iso}[1]{\ifinner\Tcase{\AR}{\isosign{#1}}%
\else\Dcase{\AR}{\isosign{#1}}\fi}%

\newcommand{\biar}{\ifinner\tbicase{\BIAR}\else\dbicase{\BIAR}\fi}%

\newcommand{\Biar}[2]{\ifinner\Tbicase{\BIAR}{#1}{#2}%
\else\Dbicase{\BIAR}{#1}{#2}\fi}%

\newcommand{\bidist}{\ifinner\tbicase{\BIDIST}\else\dbicase{\BIDIST}\fi}%

\newcommand{\Bidist}[2]{\ifinner\Tbicase{\BIDIST}{\TUP{#1}}{\TDOWN{#2}}%
\else\Dbicase{\BIDIST}{\DUP{#1}}{\DDOWN{#2}}\fi}%

\newcommand{\eql}{\ifinner\tcase{\EQL}\else\dcase{\EQL}\fi}%

\newcommand{\Eql}[1]{\ifinner\Tcase{\EQL}{\TUP{#1}}%
\else\Dcase{\EQL}{\DUP{#1}}\fi}%

\newcommand{\lline}{\ifinner\tcase{\LLINE}\else\dcase{\LLINE}\fi}%

\newcommand{\Line}[1]{\ifinner\Tcase{\LLINE}{\TUP{#1}}%
\else\Dcase{\LLINE}{\DUP{#1}}\fi}%

\newcommand{\adjar}{\ifinner\tbicase{\ADJAR}\else\dbicase{\ADJAR}\fi}%

\newcommand{\Adjar}[2]{\ifinner\Tbicase{\ADJAR}{#1}{#2}%
\else\Dbicase{\ADJAR}{#1}{#2}\fi}%

\newcommand{\adjdist}{\ifinner\tbicase{\ADJDIST}\else\dbicase{\ADJDIST}\fi}%

\newcommand{\Adjdist}[2]{\ifinner\Tbicase{\ADJDIST}{\TUP{#1}}{\TDOWN{#2}}%
\else\Dbicase{\ADJDIST}{\DUP{#1}}{\DDOWN{#2}}\fi}%

\newcommand{\triar}{\ifinner\ttricase{\TRIAR}\else\dtricase{\TRIAR}\fi}%

\newcommand{\Triar}[3]{\ifinner\Ttricase{\TRIAR}{#1}{#2}{#3}%
\else\Dtricase{\TRIAR}{#1}{#2}{#3}\fi}%

\newcommand{\triadjar}{\ifinner\ttricase{\TRIADJAR}%
\else\dtricase{\TRIADJAR}\fi}%

\newcommand{\Triadjar}[3]{\ifinner\Ttricase{\TRIADJAR}{#1}{#2}{#3}%
\else\Dtricase{\TRIADJAR}{#1}{#2}{#3}\fi}%

\newcommand{\quadriar}{\ifinner\tquadricase{\QUADRIAR}%
\else\dmulticase{\QUADRIAR}\fi}%

\newcommand{\quadriadjar}{\ifinner\tmulticase{\QUADRIADJAR}%
\else\dmulticase{\QUADRIADJAR}\fi}%

\newcommand{\quintiar}{\ifinner\tmulticase{\QUINTIAR}%
\else\dmulticase{\QUINTIAR}\fi}%

\newcommand{\quintiadjar}{\ifinner\tmulticase{\QUINTIADJAR}%
\else\dmulticase{\QUINTIADJAR}\fi}%

\newcommand{\BKAR}[1]%
{\begin{picture}(#1,0)%
\put(#1,0){\line(-1,0){#1}}%
\put(0,0){\whead}%
\end{picture}}%

\newcommand{\BKDIST}[1]%
{\begin{picture}(#1,0)%
\put(#1,0){\line(-1,0){#1}}%
\put(0,0){\whead}%
\NUMBER=#1%
\divide\NUMBER by 2%
\put(\NUMBER,0){\distsign}%
\end{picture}}%

\newcommand{\BKDOTAR}[1]%
{\NUMBEROFDOTS=#1%
\divide\NUMBEROFDOTS by 300%
\advance\NUMBEROFDOTS by 1%
\begin{picture}(#1,0)%
\multiput(#1,0)(-300,0){\NUMBEROFDOTS}{\circle*{100}}%
\put(0,0){\whead}%
\end{picture}}%

\newcommand{\BKMONO}[1]%
{\monolength=#1%
\advance\monolength by -\monotail%
\begin{picture}(#1,0)(-#1,0)%
\put(-\monotail,0){\line(-1,0){\monolength}}%
\put(-\monotail,0){\whead}%
\put(-#1,0){\whead}%
\end{picture}}%

\newcommand{\BKEPI}[1]%
{\epilength=#1%
\advance\epilength by -\epihead%
\begin{picture}(#1,0)%
\put(#1,0){\line(-1,0){\epilength}}%
\put(\epihead,0){\whead}%
\put(0,0){\whead}%
\end{picture}}%

\newcommand{\BKBIMO}[1]%
{\monolength=#1%
\advance\monolength by -\monotail%
\epilength=\monolength%
\advance\epilength by -\epihead%
\begin{picture}(#1,0)%
\put(\monolength,0){\line(-1,0){\epilength}}%
\put(\monolength,0){\whead}%
\put(\epihead,0){\whead}%
\put(0,0){\whead}%
\end{picture}}%

\newcommand{\BKBIAR}[1]%
{\begin{picture}(#1,700)%
\put(#1,0){\line(-1,0){#1}}%
\put(0,0){\whead}%
\put(#1,700){\line(-1,0){#1}}%
\put(0,700){\whead}%
\end{picture}}%

\newcommand{\BKBIDIST}[1]%
{\begin{picture}(#1,700)%
\put(#1,0){\line(-1,0){#1}}%
\put(0,0){\whead}%
\put(#1,700){\line(-1,0){#1}}%
\put(0,700){\whead}%
\NUMBER=#1%
\divide\NUMBER by 2%
\put(\NUMBER,0){\distsign}%
\put(\NUMBER,700){\distsign}%
\end{picture}}%

\newcommand{\BKADJAR}[1]%
{\begin{picture}(#1,700)%
\put(0,700){\line(1,0){#1}}%
\put(#1,700){\ehead}%
\put(#1,0){\line(-1,0){#1}}%
\put(0,0){\whead}%
\end{picture}}%

\newcommand{\BKADJDIST}[1]%
{\begin{picture}(#1,700)%
\put(0,700){\line(1,0){#1}}%
\put(#1,700){\ehead}%
\put(#1,0){\line(-1,0){#1}}%
\put(0,0){\whead}%
\NUMBER=#1%
\divide\NUMBER by 2%
\put(\NUMBER,0){\distsign}%
\put(\NUMBER,700){\distsign}%
\end{picture}}%

\newcommand{\BKTRIAR}[4]%
{%
\begin{tabular}{c}%
\mbox{$#2$}\\ \BKAR{#1}\\%
\mbox{$#3$}\\ \BKAR{#1}\\%
\mbox{$#4$}\\ \BKAR{#1}%
\end{tabular}%
}%

\newcommand{\BKTRIADJAR}[4]%
{%
\begin{tabular}{c}%
\mbox{$#2$}\\ \AR{#1}\\%
\mbox{$#3$}\\ \BKAR{#1}\\%
\mbox{$#4$}\\ \AR{#1}%
\end{tabular}%
}%

\newcommand{\BKQUADRIAR}[1]%
{%
\begin{tabular}{c}%
\BKAR{#1}\\%
\rule{0pt}{3pt}\\ \BKAR{#1}\\%
\rule{0pt}{3pt}\\ \BKAR{#1}\\%
\rule{0pt}{3pt}\\ \BKAR{#1}%
\end{tabular}%
}%

\newcommand{\BKQUADRIADJAR}[1]%
{%
\begin{tabular}{c}%
\AR{#1}\\%
\rule{0pt}{3pt}\\ \BKAR{#1}\\%
\rule{0pt}{3pt}\\ \AR{#1}\\%
\rule{0pt}{3pt}\\ \BKAR{#1}%
\end{tabular}%
}%

\newcommand{\BKQUINTIAR}[1]%
{%
\begin{tabular}{c}%
\BKAR{#1}\\%
\rule{0pt}{2pt}\\ \BKAR{#1}\\%
\rule{0pt}{2pt}\\ \BKAR{#1}\\%
\rule{0pt}{2pt}\\ \BKAR{#1}\\%
\rule{0pt}{2pt}\\ \BKAR{#1}%
\end{tabular}%
}%

\newcommand{\BKQUINTIADJAR}[1]%
{%
\begin{tabular}{c}%
\BKAR{#1}\\%
\rule{0pt}{2pt}\\ \AR{#1}\\%
\rule{0pt}{2pt}\\ \BKAR{#1}\\%
\rule{0pt}{2pt}\\ \AR{#1}\\%
\rule{0pt}{2pt}\\ \BKAR{#1}%
\end{tabular}%
}%

\newcommand{\bkar}{\ifinner\tcase{\BKAR}\else\dcase{\BKAR}\fi}%

\newcommand{\Bkar}[1]{\ifinner\Tcase{\BKAR}{#1}\else\Dcase{\BKAR}{#1}\fi}%

\newcommand{\bkdist}{\ifinner\tcase{\BKDIST}\else\dcase{\BKDIST}\fi}%

\newcommand{\Bkdist}[1]{\ifinner\Tcase{\BKDIST}{\TUP{#1}}%
\else\Dcase{\BKDIST}{\TUP{#1}}\fi}%

\newcommand{\bkdotar}{\ifinner\tcase{\BKDOTAR}\else\dcase{\BKDOTAR}\fi}%

\newcommand{\Bkdotar}[1]{\ifinner\Tcase{\BKDOTAR}{#1}%
\else\Dcase{\BKDOTAR}{#1}\fi}%

\newcommand{\bkmono}{\ifinner\tcase{\BKMONO}\else\dcase{\BKMONO}\fi}%

\newcommand{\Bkmono}[1]{\ifinner\Tcase{\BKMONO}{#1}%
\else\Dcase{\BKMONO}{#1}\fi}%

\newcommand{\bkepi}{\ifinner\tcase{\BKEPI}\else\dcase{\BKEPI}\fi}%

\newcommand{\Bkepi}[1]{\ifinner\Tcase{\BKEPI}{#1}%
\else\Dcase{\BKEPI}{#1}\fi}%

\newcommand{\bkbimo}{\ifinner\tcase{\BKBIMO}\else\dcase{\BKBIMO}\fi}%

\newcommand{\Bkbimo}[1]{\ifinner\Tcase{\BKBIMO}{\hspace{9pt}#1}%
\else\Dcase{\BKBIMO}{\hspace{9pt}#1}\fi}%

\newcommand{\bkiso}{\ifinner\Tcase{\BKAR}{\isosign}%
\else\Dcase{\BKAR}{\isosign}\fi}%

\newcommand{\Bkiso}[1]{\ifinner\Tcase{\BKAR}{\isosign{#1}}%
\else\Dcase{\BKAR}{\isosign{#1}}\fi}%

\newcommand{\bkbiar}{\ifinner\tbicase{\BKBIAR}\else\dbicase{\BKBIAR}\fi}%

\newcommand{\Bkbiar}[2]{\ifinner\Tbicase{\BKBIAR}{#1}{#2}%
\else\Dbicase{\BKBIAR}{#1}{#2}\fi}%

\newcommand{\bkbidist}{\ifinner\tbicase{\BKBIDIST}%
\else\dbicase{\BKBIDIST}\fi}%

\newcommand{\Bkbidist}[2]{\ifinner\Tbicase{\BKBIDIST}{\TUP{#1}}{\TDOWN{#2}}%
\else\Tbicase{\BKBIDIST}{\DUP{#1}}{\DDOWN{#2}}\fi}%

\newcommand{\bkadjar}{\ifinner\tbicase{\BKADJAR}%
\else\dbicase{\BKADJAR}\fi}%

\newcommand{\Bkadjar}[2]{\ifinner\Tbicase{\BKADJAR}{#1}{#2}%
\else\Dbicase{\BKADJAR}{#1}{#2}\fi}%

\newcommand{\bkadjdist}{\ifinner\tbicase{\BKADJDIST}%
\else\dbicase{\BKADJDIST}\fi}%

\newcommand{\Bkadjdist}[2]{\ifinner\Tbicase{\BKADJDIST}{\TUP{#1}}{\TDOWN{#2}}%
\else\Dbicase{\BKADJDIST}{\TUP{#1}}{\TDOWN{#2}}\fi}%

\newcommand{\bktriar}{\ifinner\ttricase{\BKTRIAR}\else\dtricase{\BKTRIAR}\fi}%

\newcommand{\Bktriar}[3]{\ifinner\Ttricase{\BKTRIAR}{#1}{#2}{#3}%
\else\Dtricase{\BKTRIAR}{#1}{#2}{#3}\fi}%

\newcommand{\bktriadjar}{\ifinner\ttricase{\BKTRIADJAR}%
\else\dtricase{\BKTRIADJAR}\fi}%

\newcommand{\Bktriadjar}[3]{\ifinner\Ttricase{\BKTRIADJAR}{#1}{#2}{#3}%
\else\Dtricase{\BKTRIADJAR}{#1}{#2}{#3}\fi}%

\newcommand{\bkquadriar}{\ifinner\tquadricase{\BKQUADRIAR}%
\else\dmulticase{\BKQUADRIAR}\fi}%

\newcommand{\bkquadriadjar}{\ifinner\tmulticase{\BKQUADRIADJAR}%
\else\dmulticase{\BKQUADRIADJAR}\fi}%

\newcommand{\bkquintiar}{\ifinner\tmulticase{\BKQUINTIAR}%
\else\dmulticase{\BKQUINTIAR}\fi}%

\newcommand{\bkquintiadjar}{\ifinner\tmulticase{\BKQUINTIADJAR}%
\else\dmulticase{\BKQUINTIADJAR}\fi}%

\def\maketopbox#1{\vbox to 0\unitlength%
{\vfill\hbox to 0\unitlength{\hss#1\hss}\vss}}

\def\makebottombox#1{\vbox to 0\unitlength%
{\vss\hbox to 0\unitlength{\hss#1\hss}\vfill}}

\newcommand{\hcase}[2]%
{\makebox[0pt]%
{\raisebox{0pt}[0pt][0pt]%
{\begin{picture}(0,0)%
\put(0,0){\makebox(0,0){#1{#2}}}%
\end{picture}%
}}}%

\newcommand{\Hcase}[3]%
{\makebox[0pt]%
{\raisebox{0pt}[0pt][0pt]%
{\begin{picture}(0,0)%
\put(0,0){\makebox(0,0){#1{#3}}}%
\truex{200}%
\put(0,\value{x}){\makebottombox{$#2$}}%
\end{picture}%
}}}%

\newcommand{\hcasE}[3]%
{\makebox[0pt]%
{\raisebox{0pt}[0pt][0pt]%
{\begin{picture}(0,0)%
\put(0,0){\makebox(0,0){#1{#3}}}%
\truex{200}%
\put(0,-\value{x}){\maketopbox{$#2$}}%
\end{picture}%
}}}%

\newcommand{\Hisocase}[4]%
{\makebox[0pt]
{\raisebox{0pt}[0pt][0pt]%
{\begin{picture}(0,0)%
\put(0,0){\makebox(0,0){#1{#4}}}%
\truex{200}%
\put(0,\value{x}){\makebottombox{$#2$}}%
\put(0,-\value{x}){\maketopbox{$#3$}}%
\end{picture}%
}}}%

\newcommand{\hbicase}[2]%
{\makebox[0pt]%
{\raisebox{0pt}[0pt][0pt]%
{\begin{picture}(0,0)%
\put(0,0){\makebox(0,0){#1{#2}}}%
\end{picture}%
}}}%

\newcommand{\Hbicase}[4]%
{\makebox[0pt]
{\raisebox{0pt}[0pt][0pt]%
{\begin{picture}(0,0)%
\put(0,0){\makebox(0,0){#1{#4}}}%
\truex{700}%
\put(0,\value{x}){\makebottombox{$#2$}}%
\put(0,-\value{x}){\maketopbox{$#3$}}%
\end{picture}%
}}}%

\newcommand{\htricase}[2]%
{\makebox[0pt]%
{\raisebox{-12pt}[0pt][0pt]{#1{#2}{}{}{}{}}}}%

\newcommand{\Htricase}[6]%
{\makebox[0pt]%
{\raisebox{-12pt}[0pt][0pt]{#1{#6}{#2}{#3}{#4}{#5}}}}%

\newcommand{\hmulticase}[2]%
{\makebox[0pt]{\raisebox{-12pt}[0pt][0pt]{#1{#2}}}}%

\newcommand{\EAR}[1]%
{\begin{picture}(#1,0)%
\put(0,0){\line(1,0){#1}}%
\put(#1,0){\ehead}%
\end{picture}}%

\newcommand{\EDIST}[1]%
{\begin{picture}(#1,0)%
\put(0,0){\line(1,0){#1}}%
\put(#1,0){\ehead}%
\truex{400}
\NUMBER=#1%
\divide\NUMBER by 2%
\put(\NUMBER,0){\distsign}
\end{picture}}%

\newcommand{\EDOTAR}[1]%
{\truex{100}\truey{300}%
\NUMBEROFDOTS=#1%
\divide\NUMBEROFDOTS by \value{y}%
\advance\NUMBEROFDOTS by 1%
\begin{picture}(#1,0)%
\multiput(0,0)(\value{y},0){\NUMBEROFDOTS}%
{\circle*{\value{x}}}%
\put(#1,0){\ehead}%
\end{picture}}%

\newcommand{\EMONO}[1]%
{\truetail
\monolength=#1%
\advance\monolength by -\truemonotail%
\begin{picture}(#1,0)%
\put(\truemonotail,0){\line(1,0){\monolength}}%
\put(#1,0){\ehead}%
\put(\truemonotail,0){\ehead}%
\end{picture}}%

\newcommand{\EEPI}[1]%
{\truehead%
\epilength=#1%
\advance\epilength by -\trueepihead%
\begin{picture}(#1,0)(-#1,0)%
\put(-#1,0){\line(1,0){\epilength}}%
\put(-\trueepihead,0){\ehead}%
\put(0,0){\ehead}%
\end{picture}}%

\newcommand{\EBIMO}[1]%
{\truehead\truetail%
\monolength=#1%
\advance\monolength by -\truemonotail%
\epilength=\monolength%
\advance\epilength by -\trueepihead%
\begin{picture}(#1,0)(-#1,0)%
\put(-\monolength,0){\line(1,0){\epilength}}%
\put(-\monolength,0){\ehead}%
\put(-\trueepihead,0){\ehead}%
\put(0,0){\ehead}%
\end{picture}}%

\newcommand{\EBIAR}[1]%
{\truex{700}%
\begin{picture}(#1,\value{x})%
\put(0,0){\line(1,0){#1}}%
\put(#1,0){\ehead}%
\put(0,\value{x}){\line(1,0){#1}}%
\put(#1,\value{x}){\ehead}%
\end{picture}}%

\newcommand{\EBIDIST}[1]%
{\truex{700}%
\begin{picture}(#1,\value{x})%
\put(0,0){\line(1,0){#1}}%
\put(#1,0){\ehead}%
\put(0,\value{x}){\line(1,0){#1}}%
\put(#1,\value{x}){\ehead}%
\NUMBER=#1%
\divide\NUMBER by 2%
\put(\NUMBER,0){\distsign}
\put(\NUMBER,\value{x}){\distsign}%
\end{picture}}%

\newcommand{\EEQL}[1]%
{\begin{picture}(#1,0)%
\truex{200}%
\put(0,\value{x}){\line(1,0){#1}}%
\put(0,0){\line(1,0){#1}}%
\end{picture}}%

\newcommand{\ELINE}[1]%
{\begin{picture}(#1,0)%
\truex{100}%
\put(0,\value{x}){\line(1,0){#1}}%
\end{picture}}%

\newcommand{\EADJAR}[1]%
{\truex{700}%
\begin{picture}(#1,\value{x})%
\put(0,0){\line(1,0){#1}}%
\put(#1,0){\ehead}%
\put(#1,\value{x}){\line(-1,0){#1}}%
\put(0,\value{x}){\whead}%
\end{picture}}%

\newcommand{\EADJDIST}[1]%
{\truex{700}%
\begin{picture}(#1,\value{x})%
\put(0,0){\line(1,0){#1}}%
\put(#1,0){\ehead}%
\put(#1,\value{x}){\line(-1,0){#1}}%
\put(0,\value{x}){\whead}%
\NUMBER=#1%
\divide\NUMBER by 2%
\put(\NUMBER,0){\distsign}
\put(\NUMBER,\value{x}){\distsign}%
\end{picture}}%

\newcommand{\ETRIAR}[5]%
{\truex{1200}%
\X=\value{x}%
\multiply\X by 2%
\begin{picture}(#1,\X)%
\put(0,0){\line(1,0){#1}}%
\put(#1,0){\ehead}%
\put(0,\value{x}){\line(1,0){#1}}%
\put(#1,\value{x}){\ehead}%
\put(0,\X){\line(1,0){#1}}%
\put(#1,\X){\ehead}%
\X=#1%
\divide\X by 2%
\truex{-600}%
\put(\X,\value{x}){\makebox(0,0){$#5$}}%
\truex{600}%
\put(\X,\value{x}){\makebox(0,0){$#4$}}%
\truex{1800}%
\put(\X,\value{x}){\makebox(0,0){$#3$}}%
\truex{3000}%
\put(\X,\value{x}){\makebox(0,0){$#2$}}%
\end{picture}}%

\newcommand{\ETRIADJAR}[5]%
{\truex{1200}%
\X=\value{x}%
\multiply\X by 2%
\begin{picture}(#1,\X)%
\put(#1,0){\line(-1,0){#1}}%
\put(0,0){\whead}%
\put(0,\value{x}){\line(1,0){#1}}%
\put(#1,\value{x}){\ehead}%
\put(#1,\X){\line(-1,0){#1}}%
\put(0,\X){\whead}%
\X=#1%
\divide\X by 2%
\truex{-600}%
\put(\X,\value{x}){\makebox(0,0){$#5$}}%
\truex{600}%
\put(\X,\value{x}){\makebox(0,0){$#4$}}%
\truex{1800}%
\put(\X,\value{x}){\makebox(0,0){$#3$}}%
\truex{3000}%
\put(\X,\value{x}){\makebox(0,0){$#2$}}%
\end{picture}}%

\newcommand{\EQUADRIAR}[1]%
{\truex{800}%
\truey{1600}%
\X=\value{x}%
\multiply\X by 3%
\begin{picture}(#1,\X)%
\put(0,0){\line(1,0){#1}}%
\put(#1,0){\ehead}%
\put(0,\value{x}){\line(1,0){#1}}%
\put(#1,\value{x}){\ehead}%
\put(0,\value{y}){\line(1,0){#1}}%
\put(#1,\value{y}){\ehead}%
\put(0,\X){\line(1,0){#1}}%
\put(#1,\X){\ehead}%
\end{picture}}%

\newcommand{\EQUADRIADJAR}[1]%
{\truex{800}%
\truey{1600}%
\X=\value{x}%
\multiply\X by 3%
\begin{picture}(#1,\X)%
\put(0,0){\line(1,0){#1}}%
\put(#1,0){\ehead}%
\put(#1,\value{x}){\line(-1,0){#1}}%
\put(0,\value{x}){\whead}%
\put(0,\value{y}){\line(1,0){#1}}%
\put(#1,\value{y}){\ehead}%
\put(#1,\X){\line(-1,0){#1}}%
\put(0,\X){\whead}%
\end{picture}}%

\newcommand{\EQUINTIAR}[1]%
{\truex{600}%
\truey{1200}%
\truez{1800}%
\X=\value{x}%
\multiply\X by 4%
\begin{picture}(#1,\X)%
\put(0,0){\line(1,0){#1}}%
\put(#1,0){\ehead}%
\put(0,\value{x}){\line(1,0){#1}}%
\put(#1,\value{x}){\ehead}%
\put(0,\value{y}){\line(1,0){#1}}%
\put(#1,\value{y}){\ehead}%
\put(0,\value{z}){\line(1,0){#1}}%
\put(#1,\value{z}){\ehead}%
\put(0,\X){\line(1,0){#1}}%
\put(#1,\X){\ehead}%
\end{picture}}%

\newcommand{\EQUINTIADJAR}[1]%
{\truex{600}%
\truey{1200}%
\truez{1800}%
\X=\value{x}%
\multiply\X by 4%
\begin{picture}(#1,\X)%
\put(0,0){\line(1,0){#1}}%
\put(#1,0){\ehead}%
\put(#1,\value{x}){\line(-1,0){#1}}%
\put(0,\value{x}){\whead}%
\put(0,\value{y}){\line(1,0){#1}}%
\put(#1,\value{y}){\ehead}%
\put(#1,\value{z}){\line(-1,0){#1}}%
\put(0,\value{z}){\whead}%
\put(0,\X){\line(1,0){#1}}%
\put(#1,\X){\ehead}%
\end{picture}}%

\def\basicear[#1]{%
\Z=#1%
\multiply \Z by 100%
\hcase{\EAR}{\Z}}%

\newcommand{\ear}{\@ifnextchar[{\basicear}%
{\hspace{\SOURCE\unitlength}\basicear[\ARROWLENGTH]}}%

\def\basicEar[#1]#2{%
\Z=#1%
\multiply \Z by 100%
\Hcase{\EAR}{#2}{\Z}}%

\newcommand{\Ear}{\@ifnextchar[{\basicEar}%
{\hspace{\SOURCE\unitlength}\basicEar[\ARROWLENGTH]}}%

\def\basiceaR[#1]#2{%
\Z=#1%
\multiply \Z by 100%
\hcasE{\EAR}{#2}{\Z}}%

\newcommand{\eaR}{\@ifnextchar[{\basiceaR}%
{\hspace{\SOURCE\unitlength}\basiceaR[\ARROWLENGTH]}}%

\def\basicedist[#1]{%
\Z=#1%
\multiply \Z by 100%
\hcase{\EDIST}{\Z}}%

\newcommand{\edist}{\@ifnextchar[{\basicedist}%
{\hspace{\SOURCE\unitlength}\basicedist[\ARROWLENGTH]}}%

\def\basicEdist[#1]#2{%
\Z=#1%
\multiply \Z by 100%
\Hcase{\EDIST}{\DUP{#2}}{\Z}}%

\newcommand{\Edist}{\@ifnextchar[{\basicEdist}%
{\hspace{\SOURCE\unitlength}\basicEdist[\ARROWLENGTH]}}%

\def\basicedisT[#1]#2{%
\Z=#1%
\multiply \Z by 100%
\hcasE{\EDIST}{\DDOWN{#2}}{\Z}}%

\newcommand{\edisT}{\@ifnextchar[{\basicedisT}%
{\hspace{\SOURCE\unitlength}\basicedisT[\ARROWLENGTH]}}%

\def\basicedotar[#1]{%
\Z=#1%
\multiply \Z by 100%
\hcase{\EDOTAR}{\Z}}%

\newcommand{\edotar}{\@ifnextchar[{\basicedotar}%
{\hspace{\SOURCE\unitlength}\basicedotar[\ARROWLENGTH]}}%

\def\basicEdotar[#1]#2{%
\Z=#1%
\multiply \Z by 100%
\Hcase{\EDOTAR}{#2}{\Z}}%

\newcommand{\Edotar}{\@ifnextchar[{\basicEdotar}%
{\hspace{\SOURCE\unitlength}\basicEdotar[\ARROWLENGTH]}}%

\def\basicedotaR[#1]#2{%
\Z=#1%
\multiply \Z by 100%
\hcasE{\EDOTAR}{#2}{\Z}}%

\newcommand{\edotaR}{\@ifnextchar[{\basicedotaR}%
{\hspace{\SOURCE\unitlength}\basicedotaR[\ARROWLENGTH]}}%

\def\basicemono[#1]{%
\Z=#1%
\multiply \Z by 100%
\hcase{\EMONO}{\Z}}%

\newcommand{\emono}{\@ifnextchar[{\basicemono}%
{\hspace{\SOURCE\unitlength}\basicemono[\ARROWLENGTH]}}%

\def\basicEmono[#1]#2{%
\Z=#1%
\multiply \Z by 100%
\Hcase{\EMONO}{#2}{\Z}}%

\newcommand{\Emono}{\@ifnextchar[{\basicEmono}%
{\hspace{\SOURCE\unitlength}\basicEmono[\ARROWLENGTH]}}%

\def\basicemonO[#1]#2{%
\Z=#1%
\multiply \Z by 100%
\hcasE{\EMONO}{#2}{\Z}}%

\newcommand{\emonO}{\@ifnextchar[{\basicemonO}%
{\hspace{\SOURCE\unitlength}\basicemonO[\ARROWLENGTH]}}%

\def\basiceepi[#1]{%
\Z=#1%
\multiply \Z by 100%
\hcase{\EEPI}{\Z}}%

\newcommand{\eepi}{\@ifnextchar[{\basiceepi}%
{\hspace{\SOURCE\unitlength}\basiceepi[\ARROWLENGTH]}}%

\def\basicEepi[#1]#2{%
\Z=#1%
\multiply \Z by 100%
\Hcase{\EEPI}{#2}{\Z}}%

\newcommand{\Eepi}{\@ifnextchar[{\basicEepi}%
{\hspace{\SOURCE\unitlength}\basicEepi[\ARROWLENGTH]}}%

\def\basiceepI[#1]#2{%
\Z=#1%
\multiply \Z by 100%
\hcasE{\EEPI}{#2}{\Z}}%

\newcommand{\eepI}{\@ifnextchar[{\basiceepI}%
{\hspace{\SOURCE\unitlength}\basiceepI[\ARROWLENGTH]}}%

\def\basicebimo[#1]{%
\Z=#1%
\multiply \Z by 100%
\hcase{\EBIMO}{\Z}}%

\newcommand{\ebimo}{\@ifnextchar[{\basicebimo}%
{\hspace{\SOURCE\unitlength}\basicebimo[\ARROWLENGTH]}}%

\def\basicEbimo[#1]#2{%
\Z=#1%
\multiply \Z by 100%
\Hcase{\EBIMO}{#2}{\Z}}%

\newcommand{\Ebimo}{\@ifnextchar[{\basicEbimo}%
{\hspace{\SOURCE\unitlength}\basicEbimo[\ARROWLENGTH]}}%

\def\basicebimO[#1]#2{%
\Z=#1%
\multiply \Z by 100%
\hcasE{\EBIMO}{#2}{\Z}}%

\newcommand{\ebimO}{\@ifnextchar[{\basicebimO}%
{\hspace{\SOURCE\unitlength}\basicebimO[\ARROWLENGTH]}}%

\def\basiceiso[#1]{%
\Z=#1%
\multiply \Z by 100%
\Hisocase{\EAR}{\isosign}{}{\Z}}%

\newcommand{\eiso}{\@ifnextchar[{\basiceiso}%
{\hspace{\SOURCE\unitlength}\basiceiso[\ARROWLENGTH]}}%

\def\basicEiso[#1]#2{%
\Z=#1%
\multiply \Z by 100%
\Hisocase{\EAR}{#2}{\isosign}{\Z}}%

\newcommand{\Eiso}{\@ifnextchar[{\basicEiso}%
{\hspace{\SOURCE\unitlength}\basicEiso[\ARROWLENGTH]}}%

\def\basiceisO[#1]#2{%
\Z=#1%
\multiply \Z by 100%
\Hisocase{\EAR}{\isosign}{#2}{\Z}}%

\newcommand{\eisO}{\@ifnextchar[{\basiceisO}%
{\hspace{\SOURCE\unitlength}\basiceisO[\ARROWLENGTH]}}%

\def\basiceeql[#1]{%
\Z=#1%
\multiply \Z by 100%
\hcase{\EEQL}{\Z}}%

\newcommand{\eeql}{\@ifnextchar[{\basiceeql}%
{\hspace{\SOURCE\unitlength}\basiceeql[\ARROWLENGTH]}}%

\def\basicEeql[#1]#2{%
\Z=#1%
\multiply \Z by 100%
\Hcase{\EEQL}{\DUP{#2}}{\Z}}%

\newcommand{\Eeql}{\@ifnextchar[{\basicEeql}%
{\hspace{\SOURCE\unitlength}\basicEeql[\ARROWLENGTH]}}%

\def\basiceeqL[#1]#2{%
\Z=#1%
\multiply \Z by 100%
\hcasE{\EEQL}{#2}{\Z}}%

\newcommand{\eeqL}{\@ifnextchar[{\basiceeqL}%
{\hspace{\SOURCE\unitlength}\basiceeqL[\ARROWLENGTH]}}%

\def\basiceline[#1]{%
\Z=#1%
\multiply \Z by 100%
\hcase{\ELINE}{\Z}}%

\newcommand{\eline}{\@ifnextchar[{\basiceline}%
{\hspace{\SOURCE\unitlength}\basiceline[\ARROWLENGTH]}}%

\def\basicEline[#1]#2{%
\Z=#1%
\multiply \Z by 100%
\Hcase{\ELINE}{\DUP{#2}}{\Z}}%

\newcommand{\Eline}{\@ifnextchar[{\basicEline}%
{\hspace{\SOURCE\unitlength}\basicEline[\ARROWLENGTH]}}%

\def\basicelinE[#1]#2{%
\Z=#1%
\multiply \Z by 100%
\hcasE{\ELINE}{#2}{\Z}}%

\newcommand{\elinE}{\@ifnextchar[{\basicelinE}%
{\hspace{\SOURCE\unitlength}\basicelinE[\ARROWLENGTH]}}%

\def\basicebiar[#1]{%
\Z=#1%
\multiply \Z by 100%
\hbicase{\EBIAR}{\Z}}%

\newcommand{\ebiar}{\@ifnextchar[{\basicebiar}%
{\hspace{\SOURCE\unitlength}\basicebiar[\ARROWLENGTH]}}%

\def\basicEbiar[#1]#2#3{%
\Z=#1%
\multiply \Z by 100%
\Hbicase{\EBIAR}{#2}{#3}{\Z}}%

\newcommand{\Ebiar}{\@ifnextchar[{\basicEbiar}%
{\hspace{\SOURCE\unitlength}\basicEbiar[\ARROWLENGTH]}}%

\def\basicebidist[#1]{%
\Z=#1%
\multiply \Z by 100%
\hbicase{\EBIDIST}{\Z}}%

\newcommand{\ebidist}{\@ifnextchar[{\basicebidist}%
{\hspace{\SOURCE\unitlength}\basicebidist[\ARROWLENGTH]}}%

\def\basicEbidist[#1]#2#3{%
\Z=#1%
\multiply \Z by 100%
\Hbicase{\EBIDIST}{\DUP{#2}}{\DDOWN{#3}}{\Z}}%

\newcommand{\Ebidist}{\@ifnextchar[{\basicEbidist}%
{\hspace{\SOURCE\unitlength}\basicEbidist[\ARROWLENGTH]}}%

\def\basiceadjar[#1]{%
\Z=#1%
\multiply \Z by 100%
\hbicase{\EADJAR}{\Z}}%

\newcommand{\eadjar}{\@ifnextchar[{\basiceadjar}%
{\hspace{\SOURCE\unitlength}\basiceadjar[\ARROWLENGTH]}}%

\def\basicEadjar[#1]#2#3{%
\Z=#1%
\multiply \Z by 100%
\Hbicase{\EADJAR}{#2}{#3}{\Z}}%

\newcommand{\Eadjar}{\@ifnextchar[{\basicEadjar}%
{\hspace{\SOURCE\unitlength}\basicEadjar[\ARROWLENGTH]}}%

\def\basiceadjdist[#1]{%
\Z=#1%
\multiply \Z by 100%
\hbicase{\EADJDIST}{\Z}}%

\newcommand{\eadjdist}{\@ifnextchar[{\basiceadjdist}%
{\hspace{\SOURCE\unitlength}\basiceadjdist[\ARROWLENGTH]}}%

\def\basicEadjdist[#1]#2#3{%
\Z=#1%
\multiply \Z by 100%
\Hbicase{\EADJDIST}{\DUP{#2}}{\DDOWN{#3}}{\Z}}%

\newcommand{\Eadjdist}{\@ifnextchar[{\basicEadjdist}%
{\hspace{\SOURCE\unitlength}\basicEadjdist[\ARROWLENGTH]}}%

\def\basicetriar[#1]{%
\Z=#1%
\multiply \Z by 100%
\htricase{\ETRIAR}{\Z}}%

\newcommand{\etriar}{\@ifnextchar[{\basicetriar}%
{\hspace{\SOURCE\unitlength}\basicetriar[\ARROWLENGTH]}}%

\def\basicEtriar[#1]#2#3#4{%
\Z=#1%
\multiply \Z by 100%
\Htricase{\ETRIAR}{#2}{#3}{#4}{}{\Z}}%

\newcommand{\Etriar}{\@ifnextchar[{\basicEtriar}%
{\hspace{\SOURCE\unitlength}\basicEtriar[\ARROWLENGTH]}}%

\def\basicetriaR[#1]#2#3#4{%
\Z=#1%
\multiply \Z by 100%
\Htricase{\ETRIAR}{}{#2}{#3}{#4}{\Z}}%

\newcommand{\etriaR}{\@ifnextchar[{\basicetriaR}%
{\hspace{\SOURCE\unitlength}\basicetriaR[\ARROWLENGTH]}}%

\def\basicetriadjar[#1]{%
\Z=#1%
\multiply \Z by 100%
\htricase{\ETRIADJAR}{\Z}}%

\newcommand{\etriadjar}{\@ifnextchar[{\basicetriadjar}%
{\hspace{\SOURCE\unitlength}\basicetriadjar[\ARROWLENGTH]}}%

\def\basicEtriadjar[#1]#2#3#4{%
\Z=#1%
\multiply \Z by 100%
\Htricase{\ETRIADJAR}{#2}{#3}{#4}{}{\Z}}%

\newcommand{\Etriadjar}{\@ifnextchar[{\basicEtriadjar}%
{\hspace{\SOURCE\unitlength}\basicEtriadjar[\ARROWLENGTH]}}%

\def\basicetriadjaR[#1]#2#3#4{%
\Z=#1%
\multiply \Z by 100%
\Htricase{\ETRIADJAR}{}{#2}{#3}{#4}{\Z}}%

\newcommand{\etriadjaR}{\@ifnextchar[{\basicetriadjaR}%
{\hspace{\SOURCE\unitlength}\basicetriadjaR[\ARROWLENGTH]}}%

\def\basicequadriar[#1]{%
\Z=#1%
\multiply \Z by 100%
\hmulticase{\EQUADRIAR}{\Z}}%

\newcommand{\equadriar}{\@ifnextchar[{\basicequadriar}%
{\hspace{\SOURCE\unitlength}\basicequadriar[\ARROWLENGTH]}}%

\def\basicequadriadjar[#1]{%
\Z=#1%
\multiply \Z by 100%
\hmulticase{\EQUADRIADJAR}{\Z}}%

\newcommand{\equadriadjar}{\@ifnextchar[{\basicequadriadjar}%
{\hspace{\SOURCE\unitlength}\basicequadriadjar[\ARROWLENGTH]}}%

\def\basicequintiar[#1]{%
\Z=#1%
\multiply \Z by 100%
\hmulticase{\EQUINTIAR}{\Z}}%

\newcommand{\equintiar}{\@ifnextchar[{\basicequintiar}%
{\hspace{\SOURCE\unitlength}\basicequintiar[\ARROWLENGTH]}}%

\def\basicequintiadjar[#1]{%
\Z=#1%
\multiply \Z by 100%
\hmulticase{\EQUINTIADJAR}{\Z}}%

\newcommand{\equintiadjar}{\@ifnextchar[{\basicequintiadjar}%
{\hspace{\SOURCE\unitlength}\basicequintiadjar[\ARROWLENGTH]}}%

\newcommand{\WAR}[1]%
{\begin{picture}(#1,0)%
\put(#1,0){\line(-1,0){#1}}%
\put(0,0){\whead}%
\end{picture}}%

\newcommand{\WDIST}[1]%
{\begin{picture}(#1,0)%
\put(#1,0){\line(-1,0){#1}}%
\put(0,0){\whead}%
\NUMBER=#1%
\divide\NUMBER by 2%
\put(\NUMBER,0){\distsign}%
\end{picture}}%

\newcommand{\WDOTAR}[1]%
{\truex{100}\truey{300}%
\NUMBEROFDOTS=#1%
\divide\NUMBEROFDOTS by \value{y}%
\advance\NUMBEROFDOTS by 1%
\begin{picture}(#1,0)%
\multiput(#1,0)(-\value{y},0){\NUMBEROFDOTS}%
{\circle*{\value{x}}}%
\put(0,0){\whead}%
\end{picture}}%

\newcommand{\WMONO}[1]%
{\truetail%
\monolength=#1%
\advance\monolength by -\truemonotail%
\begin{picture}(#1,0)(-#1,0)%
\put(-\truemonotail,0){\line(-1,0){\monolength}}%
\put(-\truemonotail,0){\whead}%
\put(-#1,0){\whead}%
\end{picture}}%

\newcommand{\WEPI}[1]%
{\truehead%
\epilength=#1%
\advance\epilength by -\trueepihead%
\begin{picture}(#1,0)%
\put(#1,0){\line(-1,0){\epilength}}%
\put(\trueepihead,0){\whead}%
\put(0,0){\whead}%
\end{picture}}%

\newcommand{\WBIMO}[1]%
{\truehead\truetail%
\monolength=#1
\advance\monolength by -\truemonotail%
\epilength=\monolength%
\advance\epilength by -\trueepihead%
\begin{picture}(#1,0)%
\put(\monolength,0){\line(-1,0){\epilength}}%
\put(\monolength,0){\whead}%
\put(\trueepihead,0){\whead}%
\put(0,0){\whead}%
\end{picture}}%

\newcommand{\WBIAR}[1]%
{\truex{700}%
\begin{picture}(#1,\value{x})%
\put(#1,0){\line(-1,0){#1}}%
\put(0,0){\whead}%
\put(#1,\value{x}){\line(-1,0){#1}}%
\put(0,\value{x}){\whead}%
\end{picture}}%

\newcommand{\WBIDIST}[1]%
{\truex{700}%
\begin{picture}(#1,\value{x})%
\put(#1,0){\line(-1,0){#1}}%
\put(0,0){\whead}%
\put(#1,\value{x}){\line(-1,0){#1}}%
\put(0,\value{x}){\whead}%
\NUMBER=#1%
\divide\NUMBER by 2%
\put(\NUMBER,0){\distsign}%
\put(\NUMBER,\value{x}){\distsign}%
\end{picture}}%

\newcommand{\WADJAR}[1]%
{\truex{700}%
\begin{picture}(#1,\value{x})%
\put(0,\value{x}){\line(1,0){#1}}%
\put(#1,\value{x}){\ehead}%
\put(#1,0){\line(-1,0){#1}}%
\put(0,0){\whead}%
\end{picture}}%

\newcommand{\WADJDIST}[1]%
{\truex{700}%
\begin{picture}(#1,\value{x})%
\put(0,\value{x}){\line(1,0){#1}}%
\put(#1,\value{x}){\ehead}%
\put(#1,0){\line(-1,0){#1}}%
\put(0,0){\whead}%
\NUMBER=#1%
\divide\NUMBER by 2%
\put(\NUMBER,0){\distsign}%
\put(\NUMBER,\value{x}){\distsign}%
\end{picture}}%

\newcommand{\WTRIAR}[5]%
{\truex{1200}%
\X=\value{x}%
\multiply\X by 2%
\begin{picture}(#1,\X)%
\put(#1,0){\line(-1,0){#1}}%
\put(0,0){\whead}%
\put(#1,\value{x}){\line(-1,0){#1}}%
\put(0,\value{x}){\whead}%
\put(#1,\X){\line(-1,0){#1}}%
\put(0,\X){\whead}%
\X=#1%
\divide\X by 2%
\truex{-600}%
\put(\X,\value{x}){\makebox(0,0){$#5$}}%
\truex{600}%
\put(\X,\value{x}){\makebox(0,0){$#4$}}%
\truex{1800}%
\put(\X,\value{x}){\makebox(0,0){$#3$}}%
\truex{3000}%
\put(\X,\value{x}){\makebox(0,0){$#2$}}%
\end{picture}}%

\newcommand{\WTRIADJAR}[5]%
{\truex{1200}%
\X=\value{x}%
\multiply\X by 2%
\begin{picture}(#1,\X)%
\put(0,0){\line(1,0){#1}}%
\put(#1,0){\ehead}%
\put(#1,\value{x}){\line(-1,0){#1}}%
\put(0,\value{x}){\whead}%
\put(0,\X){\line(1,0){#1}}%
\put(#1,\X){\ehead}%
\X=#1%
\divide\X by 2%
\truex{-600}%
\put(\X,\value{x}){\makebox(0,0){$#5$}}%
\truex{600}%
\put(\X,\value{x}){\makebox(0,0){$#4$}}%
\truex{1800}%
\put(\X,\value{x}){\makebox(0,0){$#3$}}%
\truex{3000}%
\put(\X,\value{x}){\makebox(0,0){$#2$}}%
\end{picture}}%

\newcommand{\WQUADRIAR}[1]%
{\truex{800}%
\truey{1600}%
\X=\value{x}%
\multiply\X by 3%
\begin{picture}(#1,\X)%
\put(#1,0){\line(-1,0){#1}}%
\put(0,0){\whead}%
\put(#1,\value{x}){\line(-1,0){#1}}%
\put(0,\value{x}){\whead}%
\put(#1,\value{y}){\line(-1,0){#1}}%
\put(0,\value{y}){\whead}%
\put(#1,\X){\line(-1,0){#1}}%
\put(0,\X){\whead}%
\end{picture}}%

\newcommand{\WQUADRIADJAR}[1]%
{\truex{800}%
\truey{1600}%
\X=\value{x}%
\multiply\X by 3%
\begin{picture}(#1,\X)%
\put(#1,0){\line(-1,0){#1}}%
\put(0,0){\whead}%
\put(0,\value{x}){\line(1,0){#1}}%
\put(#1,\value{x}){\ehead}%
\put(#1,\value{y}){\line(-1,0){#1}}%
\put(0,\value{y}){\whead}%
\put(0,\X){\line(1,0){#1}}%
\put(#1,\X){\ehead}%
\end{picture}}%

\newcommand{\WQUINTIAR}[1]%
{\truex{600}%
\truey{1200}%
\truez{1800}%
\X=\value{x}%
\multiply\X by 4%
\begin{picture}(#1,\X)%
\put(#1,0){\line(-1,0){#1}}%
\put(0,0){\whead}%
\put(#1,\value{x}){\line(-1,0){#1}}%
\put(0,\value{x}){\whead}%
\put(#1,\value{y}){\line(-1,0){#1}}%
\put(0,\value{y}){\whead}%
\put(#1,\value{z}){\line(-1,0){#1}}%
\put(0,\value{z}){\whead}%
\put(#1,\X){\line(-1,0){#1}}%
\put(0,\X){\whead}%
\end{picture}}%

\newcommand{\WQUINTIADJAR}[1]%
{\truex{600}%
\truey{1200}%
\truez{1800}%
\X=\value{x}%
\multiply\X by 4%
\begin{picture}(#1,\X)%
\put(#1,0){\line(-1,0){#1}}%
\put(0,0){\whead}%
\put(0,\value{x}){\line(1,0){#1}}%
\put(#1,\value{x}){\ehead}%
\put(#1,\value{y}){\line(-1,0){#1}}%
\put(0,\value{y}){\whead}%
\put(0,\value{z}){\line(1,0){#1}}%
\put(#1,\value{z}){\ehead}%
\put(#1,\X){\line(-1,0){#1}}%
\put(0,\X){\whead}%
\end{picture}}%

\def\basicwar[#1]{%
\Z=#1%
\multiply \Z by 100%
\hcase{\WAR}{\Z}}%

\newcommand{\war}{\@ifnextchar[{\basicwar}%
{\hspace{\SOURCE\unitlength}\basicwar[\ARROWLENGTH]}}%

\def\basicWar[#1]#2{%
\Z=#1%
\multiply \Z by 100%
\Hcase{\WAR}{#2}{\Z}}%

\newcommand{\War}{\@ifnextchar[{\basicWar}%
{\hspace{\SOURCE\unitlength}\basicWar[\ARROWLENGTH]}}%

\def\basicwaR[#1]#2{%
\Z=#1%
\multiply \Z by 100%
\hcasE{\WAR}{#2}{\Z}}%

\newcommand{\waR}{\@ifnextchar[{\basicwaR}%
{\hspace{\SOURCE\unitlength}\basicwaR[\ARROWLENGTH]}}%

\def\basicwdist[#1]{%
\Z=#1%
\multiply \Z by 100%
\hcase{\WDIST}{\Z}}%

\newcommand{\wdist}{\@ifnextchar[{\basicwdist}%
{\hspace{\SOURCE\unitlength}\basicwdist[\ARROWLENGTH]}}%

\def\basicWdist[#1]#2{%
\Z=#1%
\multiply \Z by 100%
\Hcase{\WDIST}{\DUP{#2}}{\Z}}%

\newcommand{\Wdist}{\@ifnextchar[{\basicWdist}%
{\hspace{\SOURCE\unitlength}\basicWdist[\ARROWLENGTH]}}%

\def\basicwdisT[#1]#2{%
\Z=#1%
\multiply \Z by 100%
\hcasE{\WDIST}{\DDOWN{#2}}{\Z}}%

\newcommand{\wdisT}{\@ifnextchar[{\basicwdisT}%
{\hspace{\SOURCE\unitlength}\basicwdisT[\ARROWLENGTH]}}%

\def\basicwdotar[#1]{%
\Z=#1%
\multiply \Z by 100%
\hcase{\WDOTAR}{\Z}}%

\newcommand{\wdotar}{\@ifnextchar[{\basicwdotar}%
{\hspace{\SOURCE\unitlength}\basicwdotar[\ARROWLENGTH]}}%

\def\basicWdotar[#1]#2{%
\Z=#1%
\multiply \Z by 100%
\Hcase{\WDOTAR}{#2}{\Z}}%

\newcommand{\Wdotar}{\@ifnextchar[{\basicWdotar}%
{\hspace{\SOURCE\unitlength}\basicWdotar[\ARROWLENGTH]}}%

\def\basicwdotaR[#1]#2{%
\Z=#1%
\multiply \Z by 100%
\hcasE{\WDOTAR}{#2}{\Z}}%

\newcommand{\wdotaR}{\@ifnextchar[{\basicwdotaR}%
{\hspace{\SOURCE\unitlength}\basicwdotaR[\ARROWLENGTH]}}%

\def\basicwmono[#1]{%
\Z=#1%
\multiply \Z by 100%
\hcase{\WMONO}{\Z}}%

\newcommand{\wmono}{\@ifnextchar[{\basicwmono}%
{\hspace{\SOURCE\unitlength}\basicwmono[\ARROWLENGTH]}}%

\def\basicWmono[#1]#2{%
\Z=#1%
\multiply \Z by 100%
\Hcase{\WMONO}{#2}{\Z}}%

\newcommand{\Wmono}{\@ifnextchar[{\basicWmono}%
{\hspace{\SOURCE\unitlength}\basicWmono[\ARROWLENGTH]}}%

\def\basicwmonO[#1]#2{%
\Z=#1%
\multiply \Z by 100%
\hcasE{\WMONO}{#2}{\Z}}%

\newcommand{\wmonO}{\@ifnextchar[{\basicwmonO}%
{\hspace{\SOURCE\unitlength}\basicwmonO[\ARROWLENGTH]}}%

\def\basicwepi[#1]{%
\Z=#1%
\multiply \Z by 100%
\hcase{\WEPI}{\Z}}%

\newcommand{\wepi}{\@ifnextchar[{\basicwepi}%
{\hspace{\SOURCE\unitlength}\basicwepi[\ARROWLENGTH]}}%

\def\basicWepi[#1]#2{%
\Z=#1%
\multiply \Z by 100%
\Hcase{\WEPI}{#2}{\Z}}%

\newcommand{\Wepi}{\@ifnextchar[{\basicWepi}%
{\hspace{\SOURCE\unitlength}\basicWepi[\ARROWLENGTH]}}%

\def\basicwepI[#1]#2{%
\Z=#1%
\multiply \Z by 100%
\hcasE{\WEPI}{#2}{\Z}}%

\newcommand{\wepI}{\@ifnextchar[{\basicwepI}%
{\hspace{\SOURCE\unitlength}\basicwepI[\ARROWLENGTH]}}%

\def\basicwbimo[#1]{%
\Z=#1%
\multiply \Z by 100%
\hcase{\WBIMO}{\Z}}%

\newcommand{\wbimo}{\@ifnextchar[{\basicwbimo}%
{\hspace{\SOURCE\unitlength}\basicwbimo[\ARROWLENGTH]}}%

\def\basicWbimo[#1]#2{%
\Z=#1%
\multiply \Z by 100%
\Hcase{\WBIMO}{#2}{\Z}}%

\newcommand{\Wbimo}{\@ifnextchar[{\basicWbimo}%
{\hspace{\SOURCE\unitlength}\basicWbimo[\ARROWLENGTH]}}%

\def\basicwbimO[#1]#2{%
\Z=#1%
\multiply \Z by 100%
\hcasE{\WBIMO}{#2}{\Z}}%

\newcommand{\wbimO}{\@ifnextchar[{\basicwbimO}%
{\hspace{\SOURCE\unitlength}\basicwbimO[\ARROWLENGTH]}}%

\def\basicwiso[#1]{%
\Z=#1%
\multiply \Z by 100%
\Hisocase{\WAR}{\isosign}{}{\Z}}%

\newcommand{\wiso}{\@ifnextchar[{\basicwiso}%
{\hspace{\SOURCE\unitlength}\basicwiso[\ARROWLENGTH]}}%

\def\basicWiso[#1]#2{%
\Z=#1%
\multiply \Z by 100%
\Hisocase{\WAR}{#2}{\isosign}{\Z}}%

\newcommand{\Wiso}{\@ifnextchar[{\basicWiso}%
{\hspace{\SOURCE\unitlength}\basicWiso[\ARROWLENGTH]}}%

\def\basicwisO[#1]#2{%
\Z=#1%
\multiply \Z by 100%
\Hisocase{\WAR}{\isosign}{#2}{\Z}}%

\newcommand{\wisO}{\@ifnextchar[{\basicwisO}%
{\hspace{\SOURCE\unitlength}\basicwisO[\ARROWLENGTH]}}%

\def\basicwbiar[#1]{%
\Z=#1%
\multiply \Z by 100%
\hbicase{\WBIAR}{\Z}}%

\newcommand{\wbiar}{\@ifnextchar[{\basicwbiar}%
{\hspace{\SOURCE\unitlength}\basicwbiar[\ARROWLENGTH]}}%

\def\basicWbiar[#1]#2#3{%
\Z=#1%
\multiply \Z by 100%
\Hbicase{\WBIAR}{#2}{#3}{\Z}}%

\newcommand{\Wbiar}{\@ifnextchar[{\basicWbiar}%
{\hspace{\SOURCE\unitlength}\basicWbiar[\ARROWLENGTH]}}%

\def\basicwbidist[#1]{%
\Z=#1%
\multiply \Z by 100%
\hbicase{\WBIDIST}{\Z}}%

\newcommand{\wbidist}{\@ifnextchar[{\basicwbidist}%
{\hspace{\SOURCE\unitlength}\basicwbidist[\ARROWLENGTH]}}%

\def\basicWbidist[#1]#2#3{%
\Z=#1%
\multiply \Z by 100%
\Hbicase{\WBIDIST}{\DUP{#2}}{\DDOWN{#3}}{\Z}}%

\newcommand{\Wbidist}{\@ifnextchar[{\basicWbidist}%
{\hspace{\SOURCE\unitlength}\basicWbidist[\ARROWLENGTH]}}%

\def\basicwadjar[#1]{%
\Z=#1%
\multiply \Z by 100%
\hbicase{\WADJAR}{\Z}}%

\newcommand{\wadjar}{\@ifnextchar[{\basicwadjar}%
{\hspace{\SOURCE\unitlength}\basicwadjar[\ARROWLENGTH]}}%

\def\basicWadjar[#1]#2#3{%
\Z=#1%
\multiply \Z by 100%
\Hbicase{\WADJAR}{#2}{#3}{\Z}}%

\newcommand{\Wadjar}{\@ifnextchar[{\basicWadjar}%
{\hspace{\SOURCE\unitlength}\basicWadjar[\ARROWLENGTH]}}%

\def\basicwadjdist[#1]{%
\Z=#1%
\multiply \Z by 100%
\hbicase{\WADJDIST}{\Z}}%

\newcommand{\wadjdist}{\@ifnextchar[{\basicwadjdist}%
{\hspace{\SOURCE\unitlength}\basicwadjdist[\ARROWLENGTH]}}%

\def\basicWadjdist[#1]#2#3{%
\Z=#1%
\multiply \Z by 100%
\Hbicase{\WADJDIST}{\DUP{#2}}{\DDOWN{#3}}{\Z}}%

\newcommand{\Wadjdist}{\@ifnextchar[{\basicWadjdist}%
{\hspace{\SOURCE\unitlength}\basicWadjdist[\ARROWLENGTH]}}%

\def\basicwtriar[#1]{%
\Z=#1%
\multiply \Z by 100%
\htricase{\WTRIAR}{\Z}}%

\newcommand{\wtriar}{\@ifnextchar[{\basicwtriar}%
{\hspace{\SOURCE\unitlength}\basicwtriar[\ARROWLENGTH]}}%

\def\basicWtriar[#1]#2#3#4{%
\Z=#1%
\multiply \Z by 100%
\Htricase{\WTRIAR}{#2}{#3}{#4}{}{\Z}}%

\newcommand{\Wtriar}{\@ifnextchar[{\basicWtriar}%
{\hspace{\SOURCE\unitlength}\basicWtriar[\ARROWLENGTH]}}%

\def\basicwtriaR[#1]#2#3#4{%
\Z=#1%
\multiply \Z by 100%
\Htricase{\WTRIAR}{}{#2}{#3}{#4}{\Z}}%

\newcommand{\wtriaR}{\@ifnextchar[{\basicwtriaR}%
{\hspace{\SOURCE\unitlength}\basicwtriaR[\ARROWLENGTH]}}%

\def\basicwtriadjar[#1]{%
\Z=#1%
\multiply \Z by 100%
\htricase{\WTRIADJAR}{\Z}}%

\newcommand{\wtriadjar}{\@ifnextchar[{\basicwtriadjar}%
{\hspace{\SOURCE\unitlength}\basicwtriadjar[\ARROWLENGTH]}}%

\def\basicWtriadjar[#1]#2#3#4{%
\Z=#1%
\multiply \Z by 100%
\Htricase{\WTRIADJAR}{#2}{#3}{#4}{}{\Z}}%

\newcommand{\Wtriadjar}{\@ifnextchar[{\basicWtriadjar}%
{\hspace{\SOURCE\unitlength}\basicWtriadjar[\ARROWLENGTH]}}%

\def\basicwtriadjaR[#1]#2#3#4{%
\Z=#1%
\multiply \Z by 100%
\Htricase{\WTRIADJAR}{}{#2}{#3}{#4}{\Z}}%

\newcommand{\wtriadjaR}{\@ifnextchar[{\basicwtriadjaR}%
{\hspace{\SOURCE\unitlength}\basicwtriadjaR[\ARROWLENGTH]}}%

\def\basicwquadriar[#1]{%
\Z=#1%
\multiply \Z by 100%
\hmulticase{\WQUADRIAR}{\Z}}%

\newcommand{\wquadriar}{\@ifnextchar[{\basicwquadriar}%
{\hspace{\SOURCE\unitlength}\basicwquadriar[\ARROWLENGTH]}}%

\def\basicwquadriadjar[#1]{%
\Z=#1%
\multiply \Z by 100%
\hmulticase{\WQUADRIADJAR}{\Z}}%

\newcommand{\wquadriadjar}{\@ifnextchar[{\basicwquadriadjar}%
{\hspace{\SOURCE\unitlength}\basicwquadriadjar[\ARROWLENGTH]}}%

\def\basicwquintiar[#1]{%
\Z=#1%
\multiply \Z by 100%
\hmulticase{\WQUINTIAR}{\Z}}%

\newcommand{\wquintiar}{\@ifnextchar[{\basicwquintiar}%
{\hspace{\SOURCE\unitlength}\basicwquintiar[\ARROWLENGTH]}}%

\def\basicwquintiadjar[#1]{%
\Z=#1%
\multiply \Z by 100%
\hmulticase{\WQUINTIADJAR}{\Z}}%

\newcommand{\wquintiadjar}{\@ifnextchar[{\basicwquintiadjar}%
{\hspace{\SOURCE\unitlength}\basicwquintiadjar[\ARROWLENGTH]}}%

\newcommand{\vcase}[2]{#1{#2}}%

\newcommand{\Vcase}[3]{\makebox[0pt]%
{\makebox[0pt][r]{\raisebox{0pt}[0pt][0pt]{${#2}\hspace{2pt}$}}}#1{#3}}%

\newcommand{\vcasE}[3]{\makebox[0pt]%
{#1{#3}\makebox[0pt][l]{\raisebox{0pt}[0pt][0pt]{\hspace{2pt}$#2$}}}}%

\newcommand{\Visocase}[4]{\makebox[0pt]%
{\makebox[0pt][r]{\raisebox{0pt}[0pt][0pt]{$#2$\hspace{2pt}}}#1{#4}%
\makebox[0pt][l]{\raisebox{0pt}[0pt][0pt]{\hspace{2pt}$#3$}}}}%

\newcommand{\vbicase}[2]{\makebox[0pt]{{#1{#2}}}}%

\newcommand{\Vbicase}[4]{\makebox[0pt]%
{\makebox[0pt][r]{\raisebox{0pt}[0pt][0pt]{$#2$\hspace{5.5pt}}}#1{#4}%
\makebox[0pt][l]{\raisebox{0pt}[0pt][0pt]{\hspace{6.5pt}$#3$}}}}%

\newcommand{\vtricase}[2]%
{\makebox[0pt]{#1{}{}{}{}{#2}}}%

\newcommand{\Vtricase}[6]%
{\makebox[0pt]{#1{#2}{#3}{#4}{#5}{#6}}}%

\newcommand{\vmulticase}[2]{\makebox[0pt]{#1{#2}}}%

\newcommand{\SAR}[1]%
{\begin{picture}(0,0)%
\put(0,0){\makebox(0,0)%
{\begin{picture}(0,#1)%
\put(0,#1){\line(0,-1){#1}}%
\put(0,0){\shead}%
\end{picture}}}\end{picture}}%

\newcommand{\SDIST}[1]%
{\begin{picture}(0,0)%
\put(0,0){\makebox(0,0)%
{\begin{picture}(0,#1)%
\put(0,#1){\line(0,-1){#1}}%
\put(0,0){\shead}%
\end{picture}}}%
\put(0,0){\distsign}%
\end{picture}}%

\newcommand{\SDOTAR}[1]%
{\truex{100}\truey{300}%
\NUMBEROFDOTS=#1%
\divide\NUMBEROFDOTS by \value{y}%
\advance\NUMBEROFDOTS by 1%
\begin{picture}(0,0)%
\put(0,0){\makebox(0,0)%
{\begin{picture}(0,#1)%
\multiput(0,#1)(0,-\value{y}){\NUMBEROFDOTS}%
{\circle*{\value{x}}}%
\put(0,0){\shead}%
\end{picture}}}\end{picture}}%

\newcommand{\SMONO}[1]%
{\truetail%
\monolength=#1%
\advance\monolength by -\truemonotail%
\begin{picture}(0,0)%
\put(0,0){\makebox(0,0)%
{\begin{picture}(0,#1)%
\put(0,\monolength){\line(0,-1){\monolength}}%
\put(0,\monolength){\shead}%
\put(0,0){\shead}%
\end{picture}}}\end{picture}}%

\newcommand{\SEPI}[1]%
{\truehead%
\epilength=#1%
\advance\epilength by -\trueepihead%
\begin{picture}(0,0)%
\put(0,0){\makebox(0,0)%
{\begin{picture}(0,#1)%
\put(0,#1){\line(0,-1){\epilength}}%
\put(0,\trueepihead){\shead}%
\put(0,0){\shead}%
\end{picture}}}\end{picture}}%

\newcommand{\SBIMO}[1]%
{\truehead\truetail%
\monolength=#1%
\advance\monolength by -\truemonotail%
\epilength=\monolength%
\advance\epilength by -\trueepihead%
\begin{picture}(0,0)%
\put(0,0){\makebox(0,0)%
{\begin{picture}(0,#1)%
\put(0,\monolength){\line(0,-1){\epilength}}%
\put(0,\monolength){\shead}%
\put(0,\trueepihead){\shead}%
\put(0,0){\shead}%
\end{picture}}}\end{picture}}%

\newcommand{\SBIAR}[1]%
{\begin{picture}(0,0)%
\truex{350}%
\put(0,0){\makebox(0,0)%
{\begin{picture}(0,#1)%
\put(-\value{x},#1){\line(0,-1){#1}}%
\put(-\value{x},0){\shead}%
\put(\value{x},#1){\line(0,-1){#1}}%
\put(\value{x},0){\shead}%
\end{picture}}}\end{picture}}%

\newcommand{\SBIDIST}[1]%
{\begin{picture}(0,0)%
\truex{350}%
\put(0,0){\makebox(0,0)%
{\begin{picture}(0,#1)%
\put(-\value{x},#1){\line(0,-1){#1}}%
\put(-\value{x},0){\shead}%
\put(\value{x},#1){\line(0,-1){#1}}%
\put(\value{x},0){\shead}%
\end{picture}}}%
\put(-\value{x},0){\distsign}%
\put(\value{x},0){\distsign}%
\end{picture}}%

\newcommand{\SEQL}[1]%
{\begin{picture}(0,0)%
\truex{100}%
\put(0,0){\makebox(0,0)%
{\begin{picture}(0,#1)\put(-\value{x},#1){\line(0,-1){#1}}%
\put(\value{x},#1){\line(0,-1){#1}}%
\end{picture}}}\end{picture}}%

\newcommand{\SLINE}[1]%
{\begin{picture}(0,0)%
\put(0,0){\makebox(0,0)%
{\begin{picture}(0,#1)\put(0,#1){\line(0,-1){#1}}%
\end{picture}}}\end{picture}}%

\newcommand{\SADJAR}[1]{\begin{picture}(0,0)%
\truex{350}%
\put(0,0){\makebox(0,0)%
{\begin{picture}(0,#1)%
\put(-\value{x},#1){\line(0,-1){#1}}%
\put(-\value{x},0){\shead}%
\put(\value{x},0){\line(0,1){#1}}%
\put(\value{x},#1){\nhead}%
\end{picture}}}\end{picture}}%

\newcommand{\SADJDIST}[1]{\begin{picture}(0,0)%
\truex{350}%
\put(0,0){\makebox(0,0)%
{\begin{picture}(0,#1)%
\put(-\value{x},#1){\line(0,-1){#1}}%
\put(-\value{x},0){\shead}%
\put(\value{x},0){\line(0,1){#1}}%
\put(\value{x},#1){\nhead}%
\end{picture}}}%
\put(-\value{x},0){\distsign}%
\put(\value{x},0){\distsign}%
\end{picture}}%

\newcommand{\STRIAR}[5]%
{\truex{1200}%
\truey{1800}%
\truez{600}%
\begin{picture}(0,0)%
\put(0,0){\makebox(0,0)%
{\begin{picture}(0,#5)%
\X=#5\divide\X by 2%
\put(-\value{x},#5){\line(0,-1){#5}}%
\put(-\value{x},0){\shead}%
\put(0,#5){\line(0,-1){#5}}%
\put(0,0){\shead}%
\put(\value{x},#5){\line(0,-1){#5}}%
\put(\value{x},0){\shead}%
\put(-\value{y},\X){\makebox(0,0){$#1$}}%
\put(\value{y},\X){\makebox(0,0){$#4$}}%
\put(-\value{z},\X){\makebox(0,0){$#2$}}%
\put(\value{z},\X){\makebox(0,0){$#3$}}%
\end{picture}}}\end{picture}}%

\newcommand{\STRIADJAR}[5]%
{\truex{1200}%
\truey{1800}%
\truez{600}%
\begin{picture}(0,0)%
\put(0,0){\makebox(0,0)%
{\begin{picture}(0,#5)%
\X=#5\divide\X by 2%
\put(-\value{x},0){\line(0,1){#5}}%
\put(-\value{x},#5){\nhead}%
\put(0,#5){\line(0,-1){#5}}%
\put(0,0){\shead}%
\put(\value{x},0){\line(0,1){#5}}%
\put(\value{x},#5){\nhead}%
\put(-\value{y},\X){\makebox(0,0){$#1$}}%
\put(\value{y},\X){\makebox(0,0){$#4$}}%
\put(-\value{z},\X){\makebox(0,0){$#2$}}%
\put(\value{z},\X){\makebox(0,0){$#3$}}%
\end{picture}}}\end{picture}}%

\newcommand{\SQUADRIAR}[1]%
{\truex{400}%
\truey{1200}%
\begin{picture}(0,0)%
\put(0,0){\makebox(0,0)%
{\begin{picture}(0,#1)%
\put(-\value{y},#1){\line(0,-1){#1}}%
\put(-\value{y},0){\shead}%
\put(-\value{x},#1){\line(0,-1){#1}}%
\put(-\value{x},0){\shead}%
\put(\value{x},#1){\line(0,-1){#1}}%
\put(\value{x},0){\shead}%
\put(\value{y},#1){\line(0,-1){#1}}%
\put(\value{y},0){\shead}%
\end{picture}}}\end{picture}}%

\newcommand{\SQUADRIADJAR}[1]%
{\truex{400}%
\truey{1200}%
\begin{picture}(0,0)%
\put(0,0){\makebox(0,0)%
{\begin{picture}(0,#1)%
\put(-\value{y},#1){\line(0,-1){#1}}%
\put(-\value{y},0){\shead}%
\put(-\value{x},0){\line(0,1){#1}}%
\put(-\value{x},#1){\nhead}%
\put(\value{x},#1){\line(0,-1){#1}}%
\put(\value{x},0){\shead}%
\put(\value{y},0){\line(0,1){#1}}%
\put(\value{y},#1){\nhead}%
\end{picture}}}\end{picture}}%

\newcommand{\SQUINTIAR}[1]%
{\truex{600}%
\truey{1200}%
\begin{picture}(0,0)%
\put(0,0){\makebox(0,0)%
{\begin{picture}(0,#1)%
\put(-\value{y},#1){\line(0,-1){#1}}%
\put(-\value{y},0){\shead}%
\put(-\value{x},#1){\line(0,-1){#1}}%
\put(-\value{x},0){\shead}%
\put(0,#1){\line(0,-1){#1}}%
\put(0,0){\shead}%
\put(\value{x},#1){\line(0,-1){#1}}%
\put(\value{x},0){\shead}%
\put(\value{y},#1){\line(0,-1){#1}}%
\put(\value{y},0){\shead}%
\end{picture}}}\end{picture}}%

\newcommand{\SQUINTIADJAR}[1]%
{\truex{600}%
\truey{1200}%
\begin{picture}(0,0)%
\put(0,0){\makebox(0,0)%
{\begin{picture}(0,#1)%
\put(-\value{y},#1){\line(0,-1){#1}}%
\put(-\value{y},0){\shead}%
\put(-\value{x},0){\line(0,1){#1}}%
\put(-\value{x},#1){\nhead}%
\put(0,#1){\line(0,-1){#1}}%
\put(0,0){\shead}%
\put(\value{x},0){\line(0,1){#1}}%
\put(\value{x},#1){\nhead}%
\put(\value{y},#1){\line(0,-1){#1}}%
\put(\value{y},0){\shead}%
\end{picture}}}\end{picture}}%

\def\basicsar[#1]{\vcase{\SAR}{#100}}%

\newcommand{\sar}{\@ifnextchar[{\basicsar}{\basicsar[50]}}%

\def\basicSar[#1]#2{\Vcase{\SAR}{#2}{#100}}%

\newcommand{\Sar}{\@ifnextchar[{\basicSar}{\basicSar[50]}}%

\def\basicsaR[#1]#2{\vcasE{\SAR}{#2}{#100}}%

\newcommand{\saR}{\@ifnextchar[{\basicsaR}{\basicsaR[50]}}%

\def\basicsdist[#1]{\vcase{\SDIST}{#100}}%

\newcommand{\sdist}{\@ifnextchar[{\basicsdist}{\basicsdist[50]}}%

\def\basicSdist[#1]#2{\Vcase{\SDIST}{#2\hspace*{2pt}}{#100}}%

\newcommand{\Sdist}{\@ifnextchar[{\basicSdist}{\basicSdist[50]}}%

\def\basicsdisT[#1]#2{\vcasE{\SDIST}{\hspace*{2pt}#2}{#100}}%

\newcommand{\sdisT}{\@ifnextchar[{\basicsdisT}{\basicsdisT[50]}}%

\def\basicsdotar[#1]{\vcase{\SDOTAR}{#100}}%

\newcommand{\sdotar}{\@ifnextchar[{\basicsdotar}{\basicsdotar[50]}}%

\def\basicSdotar[#1]#2{\Vcase{\SDOTAR}{#2}{#100}}%

\newcommand{\Sdotar}{\@ifnextchar[{\basicSdotar}{\basicSdotar[50]}}%

\def\basicsdotaR[#1]#2{\vcasE{\SDOTAR}{#2}{#100}}%

\newcommand{\sdotaR}{\@ifnextchar[{\basicsdotaR}{\basicsdotaR[50]}}%

\def\basicsmono[#1]{\vcase{\SMONO}{#100}}%

\newcommand{\smono}{\@ifnextchar[{\basicsmono}{\basicsmono[50]}}%

\def\basicSmono[#1]#2{\Vcase{\SMONO}{#2}{#100}}%

\newcommand{\Smono}{\@ifnextchar[{\basicSmono}{\basicSmono[50]}}%

\def\basicsmonO[#1]#2{\vcasE{\SMONO}{#2}{#100}}%

\newcommand{\smonO}{\@ifnextchar[{\basicsmonO}{\basicsmonO[50]}}%

\def\basicsepi[#1]{\vcase{\SEPI}{#100}}%

\newcommand{\sepi}{\@ifnextchar[{\basicsepi}{\basicsepi[50]}}%

\def\basicSepi[#1]#2{\Vcase{\SEPI}{#2}{#100}}%

\newcommand{\Sepi}{\@ifnextchar[{\basicSepi}{\basicSepi[50]}}%

\def\basicsepI[#1]#2{\vcasE{\SEPI}{#2}{#100}}%

\newcommand{\sepI}{\@ifnextchar[{\basicsepI}{\basicsepI[50]}}%

\def\basicsbimo[#1]{\vcase{\SBIMO}{#100}}%

\newcommand{\sbimo}{\@ifnextchar[{\basicsbimo}{\basicsbimo[50]}}%

\def\basicSbimo[#1]#2{\Vcase{\SBIMO}{#2}{#100}}%

\newcommand{\Sbimo}{\@ifnextchar[{\basicSbimo}{\basicSbimo[50]}}%

\def\basicsbimO[#1]#2{\vcasE{\SBIMO}{#2}{#100}}%

\newcommand{\sbimO}{\@ifnextchar[{\basicsbimO}{\basicsbimO[50]}}%

\def\basicsiso[#1]{\Visocase{\SAR}{\isosign}{}{#100}}%

\newcommand{\siso}{\@ifnextchar[{\basicsiso}{\basicsiso[50]}}%

\def\basicSiso[#1]#2{\Visocase{\SAR}{#2}{\isosign}{#100}}%

\newcommand{\Siso}{\@ifnextchar[{\basicSiso}{\basicSiso[50]}}%

\def\basicsisO[#1]#2{\Visocase{\SAR}{\isosign}{#2}{#100}}%

\newcommand{\sisO}{\@ifnextchar[{\basicsisO}{\basicsisO[50]}}%

\def\basicseql[#1]{\vcase{\SEQL}{#100}}%

\newcommand{\seql}{\@ifnextchar[{\basicseql}{\basicseql[50]}}%

\def\basicSeql[#1]#2{\Vcase{\SEQL}{#2\hspace*{2pt}}{#100}}%

\newcommand{\Seql}{\@ifnextchar[{\basicSeql}{\basicSeql[50]}}%

\def\basicseqL[#1]#2{\vcasE{\SEQL}{\hspace*{2pt}#2}{#100}}%

\newcommand{\seqL}{\@ifnextchar[{\basicseqL}{\basicseqL[50]}}%

\def\basicsline[#1]{\vcase{\SLINE}{#100}}%

\newcommand{\sline}{\@ifnextchar[{\basicsline}{\basicsline[50]}}%

\def\basicSline[#1]#2{\Vcase{\SLINE}{#2\hspace*{2pt}}{#100}}%

\newcommand{\Sline}{\@ifnextchar[{\basicSline}{\basicSline[50]}}%

\def\basicslinE[#1]#2{\vcasE{\SLINE}{\hspace*{2pt}#2}{#100}}%

\newcommand{\slinE}{\@ifnextchar[{\basicslinE}{\basicslinE[50]}}%

\def\basicsbiar[#1]{\vbicase{\SBIAR}{#100}}%

\newcommand{\sbiar}{\@ifnextchar[{\basicsbiar}{\basicsbiar[50]}}%

\def\basicSbiar[#1]#2#3{\Vbicase{\SBIAR}{#2}{#3}{#100}}%

\newcommand{\Sbiar}{\@ifnextchar[{\basicSbiar}{\basicSbiar[50]}}%

\def\basicsbidist[#1]{\vbicase{\SBIDIST}{#100}}%

\newcommand{\sbidist}{\@ifnextchar[{\basicsbidist}{\basicsbidist[50]}}%

\def\basicSbidist[#1]#2#3%
{\Vbicase{\SBIDIST}{#2\hspace*{2pt}}{\hspace*{2pt}#3}{#100}}%

\newcommand{\Sbidist}{\@ifnextchar[{\basicSbidist}{\basicSbidist[50]}}%

\def\basicsadjar[#1]{\vbicase{\SADJAR}{#100}}%

\newcommand{\sadjar}{\@ifnextchar[{\basicsadjar}{\basicsadjar[50]}}%

\def\basicSadjar[#1]#2#3{\Vbicase{\SADJAR}{#2}{#3}{#100}}%

\newcommand{\Sadjar}{\@ifnextchar[{\basicSadjar}{\basicSadjar[50]}}%

\def\basicsadjdist[#1]{\vbicase{\SADJDIST}{#100}}%

\newcommand{\sadjdist}{\@ifnextchar[{\basicsadjdist}{\basicsadjdist[50]}}%

\def\basicSadjdist[#1]#2#3%
{\Vbicase{\SADJDIST}{#2\hspace*{2pt}}{\hspace*{2pt}#3}{#100}}%

\newcommand{\Sadjdist}{\@ifnextchar[{\basicSadjdist}{\basicSadjdist[50]}}%

\def\basicstriar[#1]{\vtricase{\STRIAR}{#100}}%

\newcommand{\striar}{\@ifnextchar[{\basicstriar}{\basicstriar[50]}}%

\def\basicStriar[#1]#2#3#4{\Vtricase{\STRIAR}{#2}{#3}{#4}{}{#100}}%

\newcommand{\Striar}{\@ifnextchar[{\basicStriar}{\basicStriar[50]}}%

\def\basicstriaR[#1]#2#3#4{\Vtricase{\STRIAR}{}{#2}{#3}{#4}{#100}}%

\newcommand{\striaR}{\@ifnextchar[{\basicstriaR}{\basicstriaR[50]}}%

\def\basicstriadjar[#1]{\vtricase{\STRIADJAR}{#100}}%

\newcommand{\striadjar}{\@ifnextchar[{\basicstriadjar}{\basicstriadjar[50]}}%

\def\basicStriadjar[#1]#2#3#4{\Vtricase{\STRIADJAR}{#2}{#3}{#4}{}{#100}}%

\newcommand{\Striadjar}{\@ifnextchar[{\basicStriadjar}{\basicStriadjar[50]}}%

\def\basicstriadjaR[#1]#2#3#4{\Vtricase{\STRIADJAR}{}{#2}{#3}{#4}{#100}}%

\newcommand{\striadjaR}{\@ifnextchar[{\basicstriadjaR}{\basicstriadjaR[50]}}%

\def\basicsquadriar[#1]{%
\Z=#1%
\multiply \Z by 100%
\vmulticase{\SQUADRIAR}{\Z}}%

\newcommand{\squadriar}{\@ifnextchar[{\basicsquadriar}%
{\hspace{\SOURCE\unitlength}\basicsquadriar[\ARROWLENGTH]}}%

\def\basicsquadriadjar[#1]{%
\Z=#1%
\multiply \Z by 100%
\vmulticase{\SQUADRIADJAR}{\Z}}%

\newcommand{\squadriadjar}{\@ifnextchar[{\basicsquadriadjar}%
{\hspace{\SOURCE\unitlength}\basicsquadriadjar[\ARROWLENGTH]}}%

\def\basicsquintiar[#1]{%
\Z=#1%
\multiply \Z by 100%
\vmulticase{\SQUINTIAR}{\Z}}%

\newcommand{\squintiar}{\@ifnextchar[{\basicsquintiar}%
{\hspace{\SOURCE\unitlength}\basicsquintiar[\ARROWLENGTH]}}%

\def\basicsquintiadjar[#1]{%
\Z=#1%
\multiply \Z by 100%
\vmulticase{\SQUINTIADJAR}{\Z}}%

\newcommand{\squintiadjar}{\@ifnextchar[{\basicsquintiadjar}%
{\hspace{\SOURCE\unitlength}\basicsquintiadjar[\ARROWLENGTH]}}%

\newcommand{\NAR}[1]%
{\begin{picture}(0,0)%
\put(0,0){\makebox(0,0)%
{\begin{picture}(0,#1)%
\put(0,0){\line(0,1){#1}}%
\put(0,#1){\nhead}%
\end{picture}}}\end{picture}}%

\newcommand{\NDIST}[1]%
{\begin{picture}(0,0)%
\put(0,0){\makebox(0,0)%
{\begin{picture}(0,#1)%
\put(0,0){\line(0,1){#1}}%
\put(0,#1){\nhead}%
\end{picture}}}
\put(0,0){\distsign}%
\end{picture}}%

\newcommand{\NDOTAR}[1]%
{\truex{100}\truey{300}%
\NUMBEROFDOTS=#1%
\divide\NUMBEROFDOTS by \value{y}%
\advance\NUMBEROFDOTS by 1%
\begin{picture}(0,0)%
\put(0,0){\makebox(0,0)%
{\begin{picture}(0,#1)%
\multiput(0,0)(0,\value{y}){\NUMBEROFDOTS}%
{\circle*{\value{x}}}%
\put(0,#1){\nhead}%
\end{picture}}}\end{picture}}%

\newcommand{\NMONO}[1]%
{\truetail%
\monolength=#1%
\advance\monolength by -\truemonotail%
\begin{picture}(0,0)%
\put(0,0){\makebox(0,0)%
{\begin{picture}(0,#1)%
\put(0,\truemonotail){\line(0,1){\monolength}}%
\put(0,#1){\nhead}%
\put(0,\truemonotail){\nhead}%
\end{picture}}}\end{picture}}%

\newcommand{\NEPI}[1]%
{\truehead%
\epilength=#1%
\advance\epilength by -\trueepihead%
\begin{picture}(0,0)%
\put(0,0){\makebox(0,0)%
{\begin{picture}(0,#1)%
\put(0,0){\line(0,1){\epilength}}%
\put(0,#1){\nhead}%
\put(0,\epilength){\nhead}%
\end{picture}}}\end{picture}}%

\newcommand{\NBIMO}[1]%
{\truehead\truetail%
\epilength=#1%
\advance\epilength by -\trueepihead%
\monolength=\epilength%
\advance\monolength by -\truemonotail%
\begin{picture}(0,0)%
\put(0,0){\makebox(0,0)%
{\begin{picture}(0,#1)%
\put(0,\truemonotail){\line(0,1){\monolength}}%
\put(0,#1){\nhead}%
\put(0,\truemonotail){\nhead}%
\put(0,\epilength){\nhead}%
\end{picture}}}\end{picture}}%

\newcommand{\NBIAR}[1]%
{\begin{picture}(0,0)%
\truex{350}%
\put(0,0){\makebox(0,0)%
{\begin{picture}(0,#1)%
\put(-\value{x},0){\line(0,1){#1}}%
\put(-\value{x},#1){\nhead}%
\put(\value{x},0){\line(0,1){#1}}%
\put(\value{x},#1){\nhead}%
\end{picture}}}\end{picture}}%

\newcommand{\NBIDIST}[1]%
{\begin{picture}(0,0)%
\truex{350}%
\put(0,0){\makebox(0,0)%
{\begin{picture}(0,#1)%
\put(-\value{x},0){\line(0,1){#1}}%
\put(-\value{x},#1){\nhead}%
\put(\value{x},0){\line(0,1){#1}}%
\put(\value{x},#1){\nhead}%
\end{picture}}}
\put(-\value{x},0){\distsign}%
\put(\value{x},0){\distsign}%
\end{picture}}%

\newcommand{\NADJAR}[1]{\begin{picture}(0,0)%
\truex{350}%
\put(0,0){\makebox(0,0)%
{\begin{picture}(0,#1)%
\put(\value{x},#1){\line(0,-1){#1}}%
\put(\value{x},0){\shead}%
\put(-\value{x},0){\line(0,1){#1}}%
\put(-\value{x},#1){\nhead}%
\end{picture}}}\end{picture}}%

\newcommand{\NADJDIST}[1]{\begin{picture}(0,0)%
\truex{350}%
\put(0,0){\makebox(0,0)%
{\begin{picture}(0,#1)%
\put(\value{x},#1){\line(0,-1){#1}}%
\put(\value{x},0){\shead}%
\put(-\value{x},0){\line(0,1){#1}}%
\put(-\value{x},#1){\nhead}%
\end{picture}}}
\put(-\value{x},0){\distsign}%
\put(\value{x},0){\distsign}%
\end{picture}}%

\newcommand{\NTRIAR}[5]%
{\truex{1200}%
\truey{1800}%
\truez{600}%
\begin{picture}(0,0)%
\put(0,0){\makebox(0,0)%
{\begin{picture}(0,#5)%
\X=#5\divide\X by 2%
\put(-\value{x},0){\line(0,1){#5}}%
\put(-\value{x},#5){\nhead}%
\put(0,0){\line(0,1){#5}}%
\put(0,#5){\nhead}%
\put(\value{x},0){\line(0,1){#5}}%
\put(\value{x},#5){\nhead}%
\put(-\value{y},\X){\makebox(0,0){$#1$}}%
\put(\value{y},\X){\makebox(0,0){$#4$}}%
\put(-\value{z},\X){\makebox(0,0){$#2$}}%
\put(\value{z},\X){\makebox(0,0){$#3$}}%
\end{picture}}}\end{picture}}%

\newcommand{\NTRIADJAR}[5]%
{\truex{1200}%
\truey{1800}%
\truez{600}%
\begin{picture}(0,0)%
\put(0,0){\makebox(0,0)%
{\begin{picture}(0,#5)%
\X=#5\divide\X by 2%
\put(-\value{x},#5){\line(0,-1){#5}}%
\put(-\value{x},0){\shead}%
\put(0,0){\line(0,1){#5}}%
\put(0,#5){\nhead}%
\put(\value{x},#5){\line(0,-1){#5}}%
\put(\value{x},0){\shead}%
\put(-\value{y},\X){\makebox(0,0){$#1$}}%
\put(\value{y},\X){\makebox(0,0){$#4$}}%
\put(-\value{z},\X){\makebox(0,0){$#2$}}%
\put(\value{z},\X){\makebox(0,0){$#3$}}%
\end{picture}}}\end{picture}}%

\newcommand{\NQUADRIAR}[1]%
{\truex{400}%
\truey{1200}%
\begin{picture}(0,0)%
\put(0,0){\makebox(0,0)%
{\begin{picture}(0,#1)%
\put(-\value{y},0){\line(0,1){#1}}%
\put(-\value{y},#1){\nhead}%
\put(-\value{x},0){\line(0,1){#1}}%
\put(-\value{x},#1){\nhead}%
\put(\value{x},0){\line(0,1){#1}}%
\put(\value{x},#1){\nhead}%
\put(\value{y},0){\line(0,1){#1}}%
\put(\value{y},#1){\nhead}%
\end{picture}}}\end{picture}}%

\newcommand{\NQUADRIADJAR}[1]%
{\truex{400}%
\truey{1200}%
\begin{picture}(0,0)%
\put(0,0){\makebox(0,0)%
{\begin{picture}(0,#1)%
\put(-\value{y},0){\line(0,1){#1}}%
\put(-\value{y},#1){\nhead}%
\put(-\value{x},#1){\line(0,-1){#1}}%
\put(-\value{x},0){\shead}%
\put(\value{x},0){\line(0,1){#1}}%
\put(\value{x},#1){\nhead}%
\put(\value{y},#1){\line(0,-1){#1}}%
\put(\value{y},0){\shead}%
\end{picture}}}\end{picture}}%

\newcommand{\NQUINTIAR}[1]%
{\truex{600}%
\truey{1200}%
\begin{picture}(0,0)%
\put(0,0){\makebox(0,0)%
{\begin{picture}(0,#1)%
\put(-\value{y},0){\line(0,1){#1}}%
\put(-\value{y},#1){\nhead}%
\put(-\value{x},0){\line(0,1){#1}}%
\put(-\value{x},#1){\nhead}%
\put(0,0){\line(0,1){#1}}%
\put(0,#1){\nhead}%
\put(\value{x},0){\line(0,1){#1}}%
\put(\value{x},#1){\nhead}%
\put(\value{y},0){\line(0,1){#1}}%
\put(\value{y},#1){\nhead}%
\end{picture}}}\end{picture}}%

\newcommand{\NQUINTIADJAR}[1]%
{\truex{600}%
\truey{1200}%
\begin{picture}(0,0)%
\put(0,0){\makebox(0,0)%
{\begin{picture}(0,#1)%
\put(-\value{y},0){\line(0,1){#1}}%
\put(-\value{y},#1){\nhead}%
\put(-\value{x},#1){\line(0,-1){#1}}%
\put(-\value{x},0){\shead}%
\put(0,0){\line(0,1){#1}}%
\put(0,#1){\nhead}%
\put(\value{x},#1){\line(0,-1){#1}}%
\put(\value{x},0){\shead}%
\put(\value{y},0){\line(0,1){#1}}%
\put(\value{y},#1){\nhead}%
\end{picture}}}\end{picture}}%

\def\basicnar[#1]{\vcase{\NAR}{#100}}%

\newcommand{\nar}{\@ifnextchar[{\basicnar}{\basicnar[50]}}%

\def\basicNar[#1]#2{\Vcase{\NAR}{#2}{#100}}%

\newcommand{\Nar}{\@ifnextchar[{\basicNar}{\basicNar[50]}}%

\def\basicnaR[#1]#2{\vcasE{\NAR}{#2}{#100}}%

\newcommand{\naR}{\@ifnextchar[{\basicnaR}{\basicnaR[50]}}%

\def\basicndist[#1]{\vcase{\NDIST}{#100}}%

\newcommand{\ndist}{\@ifnextchar[{\basicndist}{\basicndist[50]}}%

\def\basicNdist[#1]#2{\Vcase{\NDIST}{#2\hspace*{2pt}}{#100}}%

\newcommand{\Ndist}{\@ifnextchar[{\basicNdist}{\basicNdist[50]}}%

\def\basicndisT[#1]#2{\vcasE{\NDIST}{\hspace*{2pt}#2}{#100}}%

\newcommand{\ndisT}{\@ifnextchar[{\basicndisT}{\basicndisT[50]}}%

\def\basicndotar[#1]{\vcase{\NDOTAR}{#100}}%

\newcommand{\ndotar}{\@ifnextchar[{\basicndotar}{\basicndotar[50]}}%

\def\basicNdotar[#1]#2{\Vcase{\NDOTAR}{#2}{#100}}%

\newcommand{\Ndotar}{\@ifnextchar[{\basicNdotar}{\basicNdotar[50]}}%

\def\basicndotaR[#1]#2{\vcasE{\NDOTAR}{#2}{#100}}%

\newcommand{\ndotaR}{\@ifnextchar[{\basicndotaR}{\basicndotaR[50]}}%

\def\basicnmono[#1]{\vcase{\NMONO}{#100}}%

\newcommand{\nmono}{\@ifnextchar[{\basicnmono}%
{\basicnmono[50]}}%

\def\basicNmono[#1]#2{\Vcase{\NMONO}{#2}{#100}}%

\newcommand{\Nmono}{\@ifnextchar[{\basicNmono}{\basicNmono[50]}}%

\def\basicnmonO[#1]#2{\vcasE{\NMONO}{#2}{#100}}%

\newcommand{\nmonO}{\@ifnextchar[{\basicnmonO}{\basicnmonO[50]}}%

\def\basicnepi[#1]{\vcase{\NEPI}{#100}}%

\newcommand{\nepi}{\@ifnextchar[{\basicnepi}{\basicnepi[50]}}%

\def\basicNepi[#1]#2{\Vcase{\NEPI}{#2}{#100}}%

\newcommand{\Nepi}{\@ifnextchar[{\basicNepi}{\basicNepi[50]}}%

\def\basicnepI[#1]#2{\vcasE{\NEPI}{#2}{#100}}%

\newcommand{\nepI}{\@ifnextchar[{\basicnepI}{\basicnepI[50]}}%

\def\basicnbimo[#1]{\vcase{\NBIMO}{#100}}%

\newcommand{\nbimo}{\@ifnextchar[{\basicnbimo}{\basicnbimo[50]}}%

\def\basicNbimo[#1]#2{\Vcase{\NBIMO}{#2}{#100}}%

\newcommand{\Nbimo}{\@ifnextchar[{\basicNbimo}{\basicNbimo[50]}}%

\def\basicnbimO[#1]#2{\vcasE{\NBIMO}{#2}{#100}}%

\newcommand{\nbimO}{\@ifnextchar[{\basicnbimO}{\basicnbimO[50]}}%

\def\basicniso[#1]{\Visocase{\NAR}{\isosign}{}{#100}}%

\newcommand{\niso}{\@ifnextchar[{\basicniso}{\basicniso[50]}}%

\def\basicNiso[#1]#2{\Visocase{\NAR}{#2}{\isosign}{#100}}%

\newcommand{\Niso}{\@ifnextchar[{\basicNiso}{\basicNiso[50]}}%

\def\basicnisO[#1]#2{\Visocase{\NAR}{\isosign}{#2}{#100}}%

\newcommand{\nisO}{\@ifnextchar[{\basicnisO}{\basicnisO[50]}}%

\def\basicnbiar[#1]{\vbicase{\NBIAR}{#100}}%

\newcommand{\nbiar}{\@ifnextchar[{\basicnbiar}{\basicnbiar[50]}}%

\def\basicNbiar[#1]#2#3{\Vbicase{\NBIAR}{#2}{#3}{#100}}%

\newcommand{\Nbiar}{\@ifnextchar[{\basicNbiar}{\basicNbiar[50]}}%

\def\basicnbidist[#1]{\vbicase{\NBIDIST}{#100}}%

\newcommand{\nbidist}{\@ifnextchar[{\basicnbidist}{\basicnbidist[50]}}%

\def\basicNbidist[#1]#2#3%
{\Vbicase{\NBIDIST}{#2\hspace*{2pt}}{\hspace*{2pt}#3}{#100}}%

\newcommand{\Nbidist}{\@ifnextchar[{\basicNbidist}{\basicNbidist[50]}}%

\def\basicnadjar[#1]{\vbicase{\NADJAR}{#100}}%

\newcommand{\nadjar}{\@ifnextchar[{\basicnadjar}{\basicnadjar[50]}}%

\def\basicNadjar[#1]#2#3{\Vbicase{\NADJAR}{#2}{#3}{#100}}%

\newcommand{\Nadjar}{\@ifnextchar[{\basicNadjar}{\basicNadjar[50]}}%

\def\basicnadjdist[#1]{\vbicase{\NADJDIST}{#100}}%

\newcommand{\nadjdist}{\@ifnextchar[{\basicnadjdist}{\basicnadjdist[50]}}%

\def\basicNadjdist[#1]#2#3%
{\Vbicase{\NADJDIST}{#2\hspace*{2pt}}{\hspace*{2pt}#3}{#100}}%

\newcommand{\Nadjdist}{\@ifnextchar[{\basicNadjdist}{\basicNadjdist[50]}}%

\def\basicntriar[#1]{\vtricase{\NTRIAR}{#100}}%

\newcommand{\ntriar}{\@ifnextchar[{\basicntriar}{\basicntriar[50]}}%

\def\basicNtriar[#1]#2#3#4{\Vtricase{\NTRIAR}{#2}{#3}{#4}{}{#100}}%

\newcommand{\Ntriar}{\@ifnextchar[{\basicNtriar}{\basicNtriar[50]}}%

\def\basicntriaR[#1]#2#3#4{\Vtricase{\NTRIAR}{}{#2}{#3}{#4}{#100}}%

\newcommand{\ntriaR}{\@ifnextchar[{\basicntriaR}{\basicntriaR[50]}}%

\def\basicntriadjar[#1]{\vtricase{\NTRIADJAR}{#100}}%

\newcommand{\ntriadjar}{\@ifnextchar[{\basicntriadjar}{\basicntriadjar[50]}}%

\def\basicNtriadjar[#1]#2#3#4{\Vtricase{\NTRIADJAR}{#2}{#3}{#4}{}{#100}}%

\newcommand{\Ntriadjar}{\@ifnextchar[{\basicNtriadjar}{\basicNtriadjar[50]}}%

\def\basicntriadjaR[#1]#2#3#4{\Vtricase{\NTRIADJAR}{}{#2}{#3}{#4}{#100}}%

\newcommand{\ntriadjaR}{\@ifnextchar[{\basicntriadjaR}{\basicntriadjaR[50]}}%

\def\basicnquadriar[#1]{%
\Z=#1%
\multiply \Z by 100%
\vmulticase{\NQUADRIAR}{\Z}}%

\newcommand{\nquadriar}{\@ifnextchar[{\basicnquadriar}%
{\hspace{\SOURCE\unitlength}\basicnquadriar[\ARROWLENGTH]}}%

\def\basicnquadriadjar[#1]{%
\Z=#1%
\multiply \Z by 100%
\vmulticase{\NQUADRIADJAR}{\Z}}%

\newcommand{\nquadriadjar}{\@ifnextchar[{\basicnquadriadjar}%
{\hspace{\SOURCE\unitlength}\basicnquadriadjar[\ARROWLENGTH]}}%

\def\basicnquintiar[#1]{%
\Z=#1%
\multiply \Z by 100%
\vmulticase{\NQUINTIAR}{\Z}}%

\newcommand{\nquintiar}{\@ifnextchar[{\basicnquintiar}%
{\hspace{\SOURCE\unitlength}\basicnquintiar[\ARROWLENGTH]}}%

\def\basicnquintiadjar[#1]{%
\Z=#1%
\multiply \Z by 100%
\vmulticase{\NQUINTIADJAR}{\Z}}%

\newcommand{\nquintiadjar}{\@ifnextchar[{\basicnquintiadjar}%
{\hspace{\SOURCE\unitlength}\basicnquintiadjar[\ARROWLENGTH]}}%

\newcommand{\fdcase}[4]{\begin{picture}(0,0)%
\put(0,0){#1{#4}}%
\truex{200}\truey{600}\truez{600}%
\put(-\value{x},-\value{x}){\makebox(0,\value{z})[r]{${#2}$}}%
\put(\value{x},-\value{y}){\makebox(0,\value{z})[l]{${#3}$}}%
\end{picture}}%

\newcommand{\fdbicase}[4]{\begin{picture}(0,0)%
\put(0,0){#1{#4}}%
\truex{600}\truey{150}\truez{600}%
\put(-\value{x},\value{y}){\makebox(0,\value{z})[r]{${#2}$}}%
\truex{300}\truey{1050}\truez{600}%
\put(\value{x},-\value{y}){\makebox(0,\value{z})[l]{${#3}$}}%
\end{picture}}%

\newcommand{\NEAR}[1]{%
\Y=#1%
\divide\Y by 2%
\begin{picture}(0,0)%
\put(-\Y,-\Y){\line(1,1){#1}}%
\put(\Y,\Y){\nehead}%
\end{picture}}%

\newcommand{\NEDIST}[1]{%
\Y=#1%
\divide\Y by 2%
\begin{picture}(0,0)%
\put(-\Y,-\Y){\line(1,1){#1}}%
\put(\Y,\Y){\nehead}%
\put(0,0){\distsign}%
\end{picture}}%

\newcommand{\NEDOTAR}[1]%
{\truex{100}\truey{212}%
\Y=#1%
\divide\Y by 2%
\NUMBEROFDOTS=#1%
\divide\NUMBEROFDOTS by \value{y}%
\advance\NUMBEROFDOTS by 1%
\begin{picture}(0,0)%
\multiput(-\Y,-\Y)(\value{y},\value{y}){\NUMBEROFDOTS}%
{\circle*{\value{x}}}%
\put(\Y,\Y){\nehead}%
\end{picture}}%

\newcommand{\NEMONO}[1]{%
\Y=#1%
\divide \Y by 2%
\Truetail%
\bimolength=#1%
\advance\bimolength by -\Truemonotail%
\monolength=\bimolength%
\advance\monolength by -\Y%
\begin{picture}(0,0)%
\put(-\monolength,-\monolength){\line(1,1){\bimolength}}%
\put(-\monolength,-\monolength){\nehead}%
\put(\Y,\Y){\nehead}%
\end{picture}}%

\newcommand{\NEEPI}[1]{%
\Y=#1%
\divide\Y by 2%
\Truehead%
\bimolength=#1%
\advance\bimolength by -\Trueepihead%
\epilength=\bimolength%
\advance\epilength by -\Y%
\begin{picture}(0,0)%
\put(-\Y,-\Y){\line(1,1){\bimolength}}%
\put(\epilength,\epilength){\nehead}%
\put(\Y,\Y){\nehead}%
\end{picture}}%

\newcommand{\NEBIMO}[1]{%
\Y=#1%
\divide\Y by 2%
\Truetail\Truehead%
\bimolength=#1%
\advance\bimolength by -\Truemonotail%
\monolength=\bimolength%
\advance\monolength by -\Y%
\advance\bimolength by -\Trueepihead%
\epilength=\bimolength%
\advance\epilength by -\monolength%
\begin{picture}(0,0)%
\put(-\monolength,-\monolength){\line(1,1){\bimolength}}%
\put(-\monolength,-\monolength){\nehead}%
\put(\epilength,\epilength){\nehead}%
\put(\Y,\Y){\nehead}%
\end{picture}}%

\newcommand{\NEBIAR}[1]{%
\Y=#1%
\divide\Y by 2%
\begin{picture}(0,0)%
\put(-\Y,-\Y){\begin{picture}(0,0)%
\truex{247}%
\put(-\value{x},\value{x}){\line(1,1){#1}}%
\put(\value{x},-\value{x}){\line(1,1){#1}}%
\monolength=#1%
\advance\monolength by -\value{x}%
\epilength=#1%
\advance\epilength by \value{x}%
\put(\monolength,\epilength){\nehead}%
\put(\epilength,\monolength){\nehead}%
\end{picture}}\end{picture}}%

\newcommand{\NEBIDIST}[1]{%
\Y=#1%
\divide\Y by 2%
\begin{picture}(0,0)%
\put(-\Y,-\Y){\begin{picture}(0,0)%
\truex{247}%
\monolength=#1%
\advance\monolength by -\value{x}%
\epilength=#1%
\advance\epilength by \value{x}%
\put(\value{x},-\value{x}){\line(1,1){#1}}%
\put(\epilength,\monolength){\nehead}%
\end{picture}}%
\put(-\Y,-\Y){\begin{picture}(0,0)%
\truex{247}%
\monolength=#1%
\advance\monolength by \value{x}%
\epilength=#1%
\advance\epilength by -\value{x}%
\put(-\value{x},\value{x}){\line(1,1){#1}}%
\put(\epilength,\monolength){\nehead}%
\end{picture}}%
\put(-\value{x},\value{x}){\distsign}%
\put(\value{x},-\value{x}){\distsign}%
\end{picture}}%

\newcommand{\NEEQL}[1]{%
\Y=#1%
\divide\Y by 2%
\begin{picture}(0,0)%
\put(-\Y,-\Y){\begin{picture}(0,0)%
\truex{70}%
\put(-\value{x},\value{x}){\line(1,1){#1}}%
\put(\value{x},-\value{x}){\line(1,1){#1}}%
\end{picture}}\end{picture}}%

\newcommand{\NELINE}[1]{%
\Y=#1%
\divide\Y by 2%
\begin{picture}(0,0)%
\put(-\Y,-\Y){\begin{picture}(0,0)%
\put(0,0){\line(1,1){#1}}%
\end{picture}}\end{picture}}%

\newcommand{\NEADJAR}[1]{%
\Y=#1%
\divide\Y by 2%
\begin{picture}(0,0)%
\put(-\Y,-\Y){\begin{picture}(0,0)%
\truex{247}%
\monolength=#1%
\advance\monolength by -\value{x}%
\epilength=#1%
\advance\epilength by \value{x}%
\put(\value{x},-\value{x}){\line(1,1){#1}}%
\put(\epilength,\monolength){\nehead}%
\end{picture}}%
\put(\Y,\Y){\begin{picture}(0,0)%
\truex{247}%
\monolength=#1%
\advance\monolength by -\value{x}%
\epilength=#1%
\advance\epilength by \value{x}%
\put(-\value{x},\value{x}){\line(-1,-1){#1}}%
\put(-\epilength,-\monolength){\swhead}%
\end{picture}}\end{picture}}%

\newcommand{\NEADJDIST}[1]{%
\Y=#1%
\divide\Y by 2%
\begin{picture}(0,0)%
\put(-\Y,-\Y){\begin{picture}(0,0)%
\truex{247}%
\monolength=#1%
\advance\monolength by -\value{x}%
\epilength=#1%
\advance\epilength by \value{x}%
\put(\value{x},-\value{x}){\line(1,1){#1}}%
\put(\epilength,\monolength){\nehead}%
\end{picture}}%
\put(\Y,\Y){\begin{picture}(0,0)%
\truex{247}%
\monolength=#1%
\advance\monolength by -\value{x}%
\epilength=#1%
\advance\epilength by \value{x}%
\put(-\value{x},\value{x}){\line(-1,-1){#1}}%
\put(-\epilength,-\monolength){\swhead}%
\end{picture}}%
\put(-\value{x},\value{x}){\distsign}%
\put(\value{x},-\value{x}){\distsign}%
\end{picture}}%

\def\basicnear[#1]{\fdcase{\NEAR}{}{}{#100}}%

\newcommand{\near}{\@ifnextchar[{\basicnear}{\basicnear[59]}}%

\def\basicNear[#1]#2{\fdcase{\NEAR}{#2}{}{#100}}%

\newcommand{\Near}{\@ifnextchar[{\basicNear}{\basicNear[59]}}%

\def\basicneaR[#1]#2{\fdcase{\NEAR}{}{#2}{#100}}%

\newcommand{\neaR}{\@ifnextchar[{\basicneaR}{\basicneaR[59]}}%

\def\basicnedist[#1]{\fdcase{\NEDIST}{}{}{#100}}%

\newcommand{\nedist}{\@ifnextchar[{\basicnedist}{\basicnedist[59]}}%

\def\basicNedist[#1]#2{\fdcase{\NEDIST}{#2}{}{#100}}%

\newcommand{\Nedist}{\@ifnextchar[{\basicNedist}{\basicNedist[59]}}%

\def\basicnedisT[#1]#2{\fdcase{\NEDIST}{}{#2}{#100}}%

\newcommand{\nedisT}{\@ifnextchar[{\basicnedisT}{\basicnedisT[59]}}%

\def\basicnedotar[#1]{\fdcase{\NEDOTAR}{}{}{#100}}%

\newcommand{\nedotar}{\@ifnextchar[{\basicnedotar}{\basicnedotar[59]}}%

\def\basicNedotar[#1]#2{\fdcase{\NEDOTAR}{#2}{}{#100}}%

\newcommand{\Nedotar}{\@ifnextchar[{\basicNedotar}{\basicNedotar[59]}}%

\def\basicnedotaR[#1]#2{\fdcase{\NEDOTAR}{}{#2}{#100}}%

\newcommand{\nedotaR}{\@ifnextchar[{\basicnedotaR}{\basicnedotaR[59]}}%

\def\basicnemono[#1]{\fdcase{\NEMONO}{}{}{#100}}%

\newcommand{\nemono}{\@ifnextchar[{\basicnemono}{\basicnemono[59]}}%

\def\basicNemono[#1]#2{\fdcase{\NEMONO}{#2}{}{#100}}%

\newcommand{\Nemono}{\@ifnextchar[{\basicNemono}{\basicNemono[59]}}%

\def\basicnemonO[#1]#2{\fdcase{\NEMONO}{}{#2}{#100}}%

\newcommand{\nemonO}{\@ifnextchar[{\basicnemonO}{\basicnemonO[59]}}%

\def\basicneepi[#1]{\fdcase{\NEEPI}{}{}{#100}}%

\newcommand{\neepi}{\@ifnextchar[{\basicneepi}{\basicneepi[59]}}%

\def\basicNeepi[#1]#2{\fdcase{\NEEPI}{#2}{}{#100}}%

\newcommand{\Neepi}{\@ifnextchar[{\basicNeepi}{\basicNeepi[59]}}%

\def\basicneepI[#1]#2{\fdcase{\NEEPI}{}{#2}{#100}}%

\newcommand{\neepI}{\@ifnextchar[{\basicneepI}{\basicneepI[59]}}%

\def\basicnebimo[#1]{\fdcase{\NEBIMO}{}{}{#100}}%

\newcommand{\nebimo}{\@ifnextchar[{\basicnebimo}{\basicnebimo[59]}}%

\def\basicNebimo[#1]#2{\fdcase{\NEBIMO}{#2}{}{#100}}%

\newcommand{\Nebimo}{\@ifnextchar[{\basicNebimo}{\basicNebimo[59]}}%

\def\basicnebimO[#1]#2{\fdcase{\NEBIMO}{}{#2}{#100}}%

\newcommand{\nebimO}{\@ifnextchar[{\basicnebimO}{\basicnebimO[59]}}%

\def\basicneiso[#1]{\fdcase{\NEAR}{\hspace{-2pt}\isosign}{}{#100}}%

\newcommand{\neiso}{\@ifnextchar[{\basicneiso}{\basicneiso[59]}}%

\def\basicNeiso[#1]#2{\fdcase{\NEAR}{#2}{\isosign}{#100}}%

\newcommand{\Neiso}{\@ifnextchar[{\basicNeiso}{\basicNeiso[59]}}%

\def\basicneisO[#1]#2{\fdcase{\NEAR}{\hspace{-2pt}\isosign}{#2}{#100}}%

\newcommand{\neisO}{\@ifnextchar[{\basicneisO}{\basicneisO[59]}}%

\def\basicneeql[#1]{\fdcase{\NEEQL}{}{}{#100}}%

\newcommand{\neeql}{\@ifnextchar[{\basicneeql}{\basicneeql[59]}}%

\def\basicNeeql[#1]#2{\fdcase{\NEEQL}{#2}{}{#100}}%

\newcommand{\Neeql}{\@ifnextchar[{\basicNeeql}{\basicNeeql[59]}}%

\def\basicneeqL[#1]#2{\fdcase{\NEEQL}{}{#2}{#100}}%

\newcommand{\neeqL}{\@ifnextchar[{\basicneeqL}{\basicneeqL[59]}}%

\def\basicneline[#1]{\fdcase{\NELINE}{}{}{#100}}%

\newcommand{\neline}{\@ifnextchar[{\basicneline}{\basicneline[59]}}%

\def\basicNeline[#1]#2{\fdcase{\NELINE}{#2}{}{#100}}%

\newcommand{\Neline}{\@ifnextchar[{\basicNeline}{\basicNeline[59]}}%

\def\basicnelinE[#1]#2{\fdcase{\NELINE}{}{#2}{#100}}%

\newcommand{\nelinE}{\@ifnextchar[{\basicnelinE}{\basicneeqL[59]}}%

\def\basicnebiar[#1]{\fdbicase{\NEBIAR}{}{}{#100}}%

\newcommand{\nebiar}{\@ifnextchar[{\basicnebiar}{\basicnebiar[59]}}%

\def\basicNebiar[#1]#2#3{\fdbicase{\NEBIAR}{#2}{#3}{#100}}%

\newcommand{\Nebiar}{\@ifnextchar[{\basicNebiar}{\basicNebiar[59]}}%

\def\basicneadjar[#1]{\fdbicase{\NEADJAR}{}{}{#100}}%

\newcommand{\neadjar}{\@ifnextchar[{\basicneadjar}{\basicneadjar[59]}}%

\def\basicNeadjar[#1]#2#3{\fdbicase{\NEADJAR}{#2}{#3}{#100}}%

\newcommand{\Neadjar}{\@ifnextchar[{\basicNeadjar}{\basicNeadjar[59]}}%

\def\basicnebidist[#1]{\fdbicase{\NEBIDIST}{}{}{#100}}%

\newcommand{\nebidist}{\@ifnextchar[{\basicnebidist}{\basicnebidist[59]}}%

\def\basicNebidist[#1]#2#3{\fdbicase{\NEBIDIST}{#2}{#3}{#100}}%

\newcommand{\Nebidist}{\@ifnextchar[{\basicNebidist}{\basicNebidist[59]}}%

\def\basicneadjdist[#1]{\fdbicase{\NEADJDIST}{}{}{#100}}%

\newcommand{\neadjdist}{\@ifnextchar[{\basicneadjdist}{\basicneadjdist[59]}}%

\def\basicNeadjdist[#1]#2#3{\fdbicase{\NEADJDIST}{#2}{#3}{#100}}%

\newcommand{\Neadjdist}{\@ifnextchar[{\basicNeadjdist}{\basicNeadjdist[59]}}%

\newcommand{\SWAR}[1]{%
\Y=#1%
\divide\Y by 2%
\begin{picture}(0,0)%
\put(\Y,\Y){\line(-1,-1){#1}}%
\put(-\Y,-\Y){\swhead}%
\end{picture}}%

\newcommand{\SWDIST}[1]{%
\Y=#1%
\divide\Y by 2%
\begin{picture}(0,0)%
\put(\Y,\Y){\line(-1,-1){#1}}%
\put(-\Y,-\Y){\swhead}%
\put(0,0){\distsign}%
\end{picture}}%

\newcommand{\SWDOTAR}[1]%
{\truex{100}\truey{212}%
\Y=#1%
\divide\Y by 2%
\NUMBEROFDOTS=#1%
\divide\NUMBEROFDOTS by \value{y}%
\advance\NUMBEROFDOTS by 1%
\begin{picture}(0,0)%
\multiput(\Y,\Y)(-\value{y},-\value{y}){\NUMBEROFDOTS}%
{\circle*{\value{x}}}%
\put(-\Y,-\Y){\swhead}%
\end{picture}}%

\newcommand{\SWMONO}[1]{%
\Y=#1%
\divide \Y by 2%
\Truetail%
\bimolength=#1%
\advance\bimolength by -\Truemonotail%
\monolength=\bimolength%
\advance\monolength by -\Y%
\begin{picture}(0,0)%
\put(\monolength,\monolength){\line(-1,-1){\bimolength}}%
\put(\monolength,\monolength){\swhead}%
\put(-\Y,-\Y){\swhead}%
\end{picture}}%

\newcommand{\SWEPI}[1]{%
\Y=#1%
\divide\Y by 2%
\Truehead%
\bimolength=#1%
\advance\bimolength by -\Trueepihead%
\epilength=\bimolength%
\advance\epilength by -\Y%
\begin{picture}(0,0)%
\put(\Y,\Y){\line(-1,-1){\bimolength}}%
\put(-\epilength,-\epilength){\swhead}%
\put(-\Y,-\Y){\swhead}%
\end{picture}}%

\newcommand{\SWBIMO}[1]{%
\Y=#1%
\divide\Y by 2%
\Truetail\Truehead%
\bimolength=#1%
\advance\bimolength by -\Truemonotail%
\monolength=\bimolength%
\advance\monolength by -\Y%
\advance\bimolength by -\Trueepihead%
\epilength=\bimolength%
\advance\epilength by -\monolength%
\begin{picture}(0,0)%
\put(\monolength,\monolength){\line(-1,-1){\bimolength}}%
\put(\monolength,\monolength){\swhead}%
\put(-\epilength,-\epilength){\swhead}%
\put(-\Y,-\Y){\swhead}%
\end{picture}}%

\newcommand{\SWBIAR}[1]{%
\Y=#1%
\divide\Y by 2%
\begin{picture}(0,0)%
\put(\Y,\Y){\begin{picture}(0,0)%
\truex{247}%
\put(\value{x},-\value{x}){\line(-1,-1){#1}}%
\put(-\value{x},\value{x}){\line(-1,-1){#1}}%
\monolength=#1%
\advance\monolength by -\value{x}%
\epilength=#1%
\advance\epilength by \value{x}%
\put(-\monolength,-\epilength){\swhead}%
\put(-\epilength,-\monolength){\swhead}%
\end{picture}}\end{picture}}%

\newcommand{\SWBIDIST}[1]{%
\Y=#1%
\divide\Y by 2%
\begin{picture}(0,0)%
\put(\Y,\Y){\begin{picture}(0,0)%
\truex{247}%
\monolength=#1%
\advance\monolength by -\value{x}%
\epilength=#1%
\advance\epilength by \value{x}%
\put(-\value{x},\value{x}){\line(-1,-1){#1}}%
\put(-\epilength,-\monolength){\swhead}%
\end{picture}}%
\put(\Y,\Y){\begin{picture}(0,0)%
\truex{247}%
\monolength=#1%
\advance\monolength by \value{x}%
\epilength=#1%
\advance\epilength by -\value{x}%
\put(\value{x},-\value{x}){\line(-1,-1){#1}}%
\put(-\epilength,-\monolength){\swhead}%
\end{picture}}%
\put(\value{x},-\value{x}){\distsign}%
\put(-\value{x},\value{x}){\distsign}%
\end{picture}}%

\newcommand{\SWADJAR}[1]{%
\Y=#1%
\divide\Y by 2%
\begin{picture}(0,0)%
\put(\Y,\Y){\begin{picture}(0,0)%
\truex{247}%
\monolength=#1%
\advance\monolength by -\value{x}%
\epilength=#1%
\advance\epilength by \value{x}%
\put(\value{x},-\value{x}){\line(-1,-1){#1}}%
\put(-\monolength,-\epilength){\swhead}%
\end{picture}}%
\put(-\Y,-\Y){\begin{picture}(0,0)%
\truex{247}%
\monolength=#1%
\advance\monolength by -\value{x}%
\epilength=#1%
\advance\epilength by \value{x}%
\put(-\value{x},\value{x}){\line(1,1){#1}}%
\put(\monolength,\epilength){\nehead}%
\end{picture}}\end{picture}}%

\newcommand{\SWADJDIST}[1]{%
\Y=#1%
\divide\Y by 2%
\begin{picture}(0,0)%
\put(\Y,\Y){\begin{picture}(0,0)%
\truex{247}%
\monolength=#1%
\advance\monolength by -\value{x}%
\epilength=#1%
\advance\epilength by \value{x}%
\put(\value{x},-\value{x}){\line(-1,-1){#1}}%
\put(-\monolength,-\epilength){\swhead}%
\end{picture}}%
\put(-\Y,-\Y){\begin{picture}(0,0)%
\truex{247}%
\monolength=#1%
\advance\monolength by -\value{x}%
\epilength=#1%
\advance\epilength by \value{x}%
\put(-\value{x},\value{x}){\line(1,1){#1}}%
\put(\monolength,\epilength){\nehead}%
\end{picture}}%
\put(-\value{x},\value{x}){\distsign}%
\put(\value{x},-\value{x}){\distsign}%
\end{picture}}%

\def\basicswar[#1]{\fdcase{\SWAR}{}{}{#100}}%

\newcommand{\swar}{\@ifnextchar[{\basicswar}{\basicswar[59]}}%

\def\basicSwar[#1]#2{\fdcase{\SWAR}{#2}{}{#100}}%

\newcommand{\Swar}{\@ifnextchar[{\basicSwar}{\basicSwar[59]}}%

\def\basicswaR[#1]#2{\fdcase{\SWAR}{}{#2}{#100}}%

\newcommand{\swaR}{\@ifnextchar[{\basicswaR}{\basicswaR[59]}}%

\def\basicswdist[#1]{\fdcase{\SWDIST}{}{}{#100}}%

\newcommand{\swdist}{\@ifnextchar[{\basicswdist}{\basicswdist[59]}}%

\def\basicSwdist[#1]#2{\fdcase{\SWDIST}{#2}{}{#100}}%

\newcommand{\Swdist}{\@ifnextchar[{\basicSwdist}{\basicSwdist[59]}}%

\def\basicswdisT[#1]#2{\fdcase{\SWDIST}{}{#2}{#100}}%

\newcommand{\swdisT}{\@ifnextchar[{\basicswdisT}{\basicswdisT[59]}}%

\def\basicswdotar[#1]{\fdcase{\SWDOTAR}{}{}{#100}}%

\newcommand{\swdotar}{\@ifnextchar[{\basicswdotar}{\basicswdotar[59]}}%

\def\basicSwdotar[#1]#2{\fdcase{\SWDOTAR}{#2}{}{#100}}%

\newcommand{\Swdotar}{\@ifnextchar[{\basicSwdotar}{\basicSwdotar[59]}}%

\def\basicswdotaR[#1]#2{\fdcase{\SWDOTAR}{}{#2}{#100}}%

\newcommand{\swdotaR}{\@ifnextchar[{\basicswdotaR}{\basicswdotaR[59]}}%

\def\basicswmono[#1]{\fdcase{\SWMONO}{}{}{#100}}%

\newcommand{\swmono}{\@ifnextchar[{\basicswmono}{\basicswmono[59]}}%

\def\basicSwmono[#1]#2{\fdcase{\SWMONO}{#2}{}{#100}}%

\newcommand{\Swmono}{\@ifnextchar[{\basicSwmono}{\basicSwmono[59]}}%

\def\basicswmonO[#1]#2{\fdcase{\SWMONO}{}{#2}{#100}}%

\newcommand{\swmonO}{\@ifnextchar[{\basicswmonO}{\basicswmonO[59]}}%

\def\basicswepi[#1]{\fdcase{\SWEPI}{}{}{#100}}%

\newcommand{\swepi}{\@ifnextchar[{\basicswepi}{\basicswepi[59]}}%

\def\basicSwepi[#1]#2{\fdcase{\SWEPI}{#2}{}{#100}}%

\newcommand{\Swepi}{\@ifnextchar[{\basicSwepi}{\basicSwepi[59]}}%

\def\basicswepI[#1]#2{\fdcase{\SWEPI}{}{#2}{#100}}%

\newcommand{\swepI}{\@ifnextchar[{\basicswepI}{\basicswepI[59]}}%

\def\basicswbimo[#1]{\fdcase{\SWBIMO}{}{}{#100}}%

\newcommand{\swbimo}{\@ifnextchar[{\basicswbimo}{\basicswbimo[59]}}%

\def\basicSwbimo[#1]#2{\fdcase{\SWBIMO}{#2}{}{#100}}%

\newcommand{\Swbimo}{\@ifnextchar[{\basicSwbimo}{\basicSwbimo[59]}}%

\def\basicswbimO[#1]#2{\fdcase{\SWBIMO}{}{#2}{#100}}%

\newcommand{\swbimO}{\@ifnextchar[{\basicswbimO}{\basicswbimO[59]}}%

\def\basicswiso[#1]{\fdcase{\SWAR}{\hspace{-2pt}\isosign}{}{#100}}%

\newcommand{\swiso}{\@ifnextchar[{\basicswiso}{\basicswiso[59]}}%

\def\basicSwiso[#1]#2{\fdcase{\SWAR}{#2}{\isosign}{#100}}%

\newcommand{\Swiso}{\@ifnextchar[{\basicSwiso}{\basicSwiso[59]}}%

\def\basicswisO[#1]#2{\fdcase{\SWAR}{\hspace{-2pt}\isosign}{#2}{#100}}%

\newcommand{\swisO}{\@ifnextchar[{\basicswisO}{\basicswisO[59]}}%

\def\basicswbiar[#1]{\fdbicase{\SWBIAR}{}{}{#100}}%

\newcommand{\swbiar}{\@ifnextchar[{\basicswbiar}{\basicswbiar[59]}}%

\def\basicSwbiar[#1]#2#3{\fdbicase{\SWBIAR}{#2}{#3}{#100}}%

\newcommand{\Swbiar}{\@ifnextchar[{\basicSwbiar}{\basicSwbiar[59]}}%

\def\basicswadjar[#1]{\fdbicase{\SWADJAR}{}{}{#100}}%

\newcommand{\swadjar}{\@ifnextchar[{\basicswadjar}{\basicswadjar[59]}}%

\def\basicSwadjar[#1]#2#3{\fdbicase{\SWADJAR}{#2}{#3}{#100}}%

\newcommand{\Swadjar}{\@ifnextchar[{\basicSwadjar}{\basicSwadjar[59]}}%

\def\basicswbidist[#1]{\fdbicase{\SWBIDIST}{}{}{#100}}%

\newcommand{\swbidist}{\@ifnextchar[{\basicswbidist}{\basicswbidist[59]}}%

\def\basicSwbidist[#1]#2#3{\fdbicase{\SWBIDIST}{#2}{#3}{#100}}%

\newcommand{\Swbidist}{\@ifnextchar[{\basicSwbidist}{\basicSwbidist[59]}}%

\def\basicswadjdist[#1]{\fdbicase{\SWADJDIST}{}{}{#100}}%

\newcommand{\swadjdist}{\@ifnextchar[{\basicswadjdist}{\basicswadjdist[59]}}%

\def\basicSwadjdist[#1]#2#3{\fdbicase{\SWADJDIST}{#2}{#3}{#100}}%

\newcommand{\Swadjdist}{\@ifnextchar[{\basicSwadjdist}{\basicSwadjdist[59]}}%

\newcommand{\sdcase}[4]{\begin{picture}(0,0)%
\put(0,0){#1{#4}}%
\truex{100}\truez{600}%
\put(\value{x},\value{x}){\makebox(0,\value{z})[l]{${#2}$}}%
\truex{300}\truey{800}%
\put(-\value{x},-\value{y}){\makebox(0,\value{z})[r]{${#3}$}}%
\end{picture}}%

\newcommand{\sdbicase}[4]{\begin{picture}(0,0)%
\put(0,0){#1{#4}}%
\truex{350}\truey{600}\truez{950}%
\put(\value{x},\value{x}){\makebox(0,\value{y})[l]{${#2}$}}%
\truex{450}\truey{600}\truez{1050}%
\put(-\value{x},-\value{z}){\makebox(0,\value{y})[r]{${#3}$}}%
\end{picture}}%

\newcommand{\SEAR}[1]{%
\Y=#1%
\divide\Y by 2%
\begin{picture}(0,0)%
\put(-\Y,\Y){\line(1,-1){#1}}%
\put(\Y,-\Y){\sehead}%
\end{picture}}%

\newcommand{\SEDIST}[1]{%
\Y=#1%
\divide\Y by 2%
\begin{picture}(0,0)%
\put(-\Y,\Y){\line(1,-1){#1}}%
\put(\Y,-\Y){\sehead}%
\put(0,0){\distsign}%
\end{picture}}%

\newcommand{\SEDOTAR}[1]%
{\truex{100}\truey{212}%
\Y=#1%
\divide\Y by 2%
\NUMBEROFDOTS=#1%
\divide\NUMBEROFDOTS by \value{y}%
\advance\NUMBEROFDOTS by 1%
\begin{picture}(0,0)%
\multiput(-\Y,\Y)(\value{y},-\value{y}){\NUMBEROFDOTS}%
{\circle*{\value{x}}}%
\put(\Y,-\Y){\sehead}%
\end{picture}}%

\newcommand{\SEMONO}[1]{%
\Y=#1%
\divide \Y by 2%
\Truetail%
\bimolength=#1%
\advance\bimolength by -\Truemonotail%
\monolength=\bimolength%
\advance\monolength by -\Y%
\begin{picture}(0,0)%
\put(-\monolength,\monolength){\line(1,-1){\bimolength}}%
\put(-\monolength,\monolength){\sehead}%
\put(\Y,-\Y){\sehead}%
\end{picture}}%

\newcommand{\SEEPI}[1]{%
\Y=#1%
\divide\Y by 2%
\Truehead%
\bimolength=#1%
\advance\bimolength by -\Trueepihead%
\epilength=\bimolength%
\advance\epilength by -\Y%
\begin{picture}(0,0)%
\put(-\Y,\Y){\line(1,-1){\bimolength}}%
\put(\epilength,-\epilength){\sehead}%
\put(\Y,-\Y){\sehead}%
\end{picture}}%

\newcommand{\SEBIMO}[1]{%
\Y=#1%
\divide\Y by 2%
\Truetail\Truehead%
\bimolength=#1%
\advance\bimolength by -\Truemonotail%
\monolength=\bimolength%
\advance\monolength by -\Y%
\advance\bimolength by -\Trueepihead%
\epilength=\bimolength%
\advance\epilength by -\monolength%
\begin{picture}(0,0)%
\put(-\monolength,\monolength){\line(1,-1){\bimolength}}%
\put(-\monolength,\monolength){\sehead}%
\put(\epilength,-\epilength){\sehead}%
\put(\Y,-\Y){\sehead}%
\end{picture}}%

\newcommand{\SEBIAR}[1]{%
\Y=#1%
\divide\Y by 2%
\begin{picture}(0,0)%
\put(-\Y,\Y){\begin{picture}(0,0)%
\truex{247}%
\put(-\value{x},-\value{x}){\line(1,-1){#1}}%
\put(\value{x},\value{x}){\line(1,-1){#1}}%
\monolength=#1%
\advance\monolength by -\value{x}%
\epilength=#1%
\advance\epilength by \value{x}%
\put(\monolength,-\epilength){\sehead}%
\put(\epilength,-\monolength){\sehead}%
\end{picture}}\end{picture}}%

\newcommand{\SEBIDIST}[1]{%
\Y=#1%
\divide\Y by 2%
\begin{picture}(0,0)%
\put(-\Y,\Y){\begin{picture}(0,0)%
\truex{247}%
\monolength=#1%
\advance\monolength by -\value{x}%
\epilength=#1%
\advance\epilength by \value{x}%
\put(\value{x},\value{x}){\line(1,-1){#1}}%
\put(\epilength,-\monolength){\sehead}%
\end{picture}}%
\put(-\Y,\Y){\begin{picture}(0,0)%
\truex{247}%
\monolength=#1%
\advance\monolength by \value{x}%
\epilength=#1%
\advance\epilength by -\value{x}%
\put(-\value{x},-\value{x}){\line(1,-1){#1}}%
\put(\epilength,-\monolength){\sehead}%
\end{picture}}%
\put(-\value{x},-\value{x}){\distsign}%
\put(\value{x},\value{x}){\distsign}%
\end{picture}}%

\newcommand{\SEEQL}[1]{%
\Y=#1%
\divide\Y by 2%
\begin{picture}(0,0)%
\put(-\Y,\Y){\begin{picture}(0,0)%
\truex{70}%
\put(-\value{x},-\value{x}){\line(1,-1){#1}}%
\put(\value{x},\value{x}){\line(1,-1){#1}}%
\end{picture}}\end{picture}}%

\newcommand{\SELINE}[1]{%
\Y=#1%
\divide\Y by 2%
\begin{picture}(0,0)%
\put(-\Y,\Y){\begin{picture}(0,0)%
\put(0,0){\line(1,-1){#1}}%
\end{picture}}\end{picture}}%

\newcommand{\SEADJAR}[1]{%
\Y=#1%
\divide\Y by 2%
\begin{picture}(0,0)%
\put(-\Y,\Y){\begin{picture}(0,0)%
\truex{247}%
\monolength=#1%
\advance\monolength by -\value{x}%
\epilength=#1%
\advance\epilength by \value{x}%
\put(-\value{x},-\value{x}){\line(1,-1){#1}}%
\put(\monolength,-\epilength){\sehead}%
\end{picture}}%
\put(\Y,-\Y){\begin{picture}(0,0)%
\truex{247}%
\monolength=#1%
\advance\monolength by -\value{x}%
\epilength=#1%
\advance\epilength by \value{x}%
\put(\value{x},\value{x}){\line(-1,1){#1}}%
\put(-\monolength,\epilength){\nwhead}%
\end{picture}}\end{picture}}%

\newcommand{\SEADJDIST}[1]{%
\Y=#1%
\divide\Y by 2%
\begin{picture}(0,0)%
\put(-\Y,\Y){\begin{picture}(0,0)%
\truex{247}%
\monolength=#1%
\advance\monolength by -\value{x}%
\epilength=#1%
\advance\epilength by \value{x}%
\put(-\value{x},-\value{x}){\line(1,-1){#1}}%
\put(\monolength,-\epilength){\sehead}%
\end{picture}}%
\put(\Y,-\Y){\begin{picture}(0,0)%
\truex{247}%
\monolength=#1%
\advance\monolength by -\value{x}%
\epilength=#1%
\advance\epilength by \value{x}%
\put(\value{x},\value{x}){\line(-1,1){#1}}%
\put(-\monolength,\epilength){\nwhead}%
\end{picture}}%
\put(-\value{x},-\value{x}){\distsign}%
\put(\value{x},\value{x}){\distsign}%
\end{picture}}%

\def\basicsear[#1]{\sdcase{\SEAR}{}{}{#100}}%

\newcommand{\sear}{\@ifnextchar[{\basicsear}{\basicsear[59]}}%

\def\basicSear[#1]#2{\sdcase{\SEAR}{#2}{}{#100}}%

\newcommand{\Sear}{\@ifnextchar[{\basicSear}{\basicSear[59]}}%

\def\basicseaR[#1]#2{\sdcase{\SEAR}{}{#2}{#100}}%

\newcommand{\seaR}{\@ifnextchar[{\basicseaR}{\basicseaR[59]}}%

\def\basicsedist[#1]{\sdcase{\SEDIST}{}{}{#100}}%

\newcommand{\sedist}{\@ifnextchar[{\basicsedist}{\basicsedist[59]}}%

\def\basicSedist[#1]#2{\sdcase{\SEDIST}{#2}{}{#100}}%

\newcommand{\Sedist}{\@ifnextchar[{\basicSedist}{\basicSedist[59]}}%

\def\basicsedisT[#1]#2{\sdcase{\SEDIST}{}{#2}{#100}}%

\newcommand{\sedisT}{\@ifnextchar[{\basicsedisT}{\basicsedisT[59]}}%

\def\basicsedotar[#1]{\sdcase{\SEDOTAR}{}{}{#100}}%

\newcommand{\sedotar}{\@ifnextchar[{\basicsedotar}{\basicsedotar[59]}}%

\def\basicSedotar[#1]#2{\sdcase{\SEDOTAR}{#2}{}{#100}}%

\newcommand{\Sedotar}{\@ifnextchar[{\basicSedotar}{\basicSedotar[59]}}%

\def\basicsedotaR[#1]#2{\sdcase{\SEDOTAR}{}{#2}{#100}}%

\newcommand{\sedotaR}{\@ifnextchar[{\basicsedotaR}{\basicsedotaR[59]}}%

\def\basicsemono[#1]{\sdcase{\SEMONO}{}{}{#100}}%

\newcommand{\semono}{\@ifnextchar[{\basicsemono}{\basicsemono[59]}}%

\def\basicSemono[#1]#2{\sdcase{\SEMONO}{#2}{}{#100}}%

\newcommand{\Semono}{\@ifnextchar[{\basicSemono}{\basicSemono[59]}}%

\def\basicsemonO[#1]#2{\sdcase{\SEMONO}{}{#2}{#100}}%

\newcommand{\semonO}{\@ifnextchar[{\basicsemonO}{\basicsemonO[59]}}%

\def\basicseepi[#1]{\sdcase{\SEEPI}{}{}{#100}}%

\newcommand{\seepi}{\@ifnextchar[{\basicseepi}{\basicseepi[59]}}%

\def\basicSeepi[#1]#2{\sdcase{\SEEPI}{#2}{}{#100}}%

\newcommand{\Seepi}{\@ifnextchar[{\basicSeepi}{\basicSeepi[59]}}%

\def\basicseepI[#1]#2{\sdcase{\SEEPI}{}{#2}{#100}}%

\newcommand{\seepI}{\@ifnextchar[{\basicseepI}{\basicseepI[59]}}%

\def\basicsebimo[#1]{\sdcase{\SEBIMO}{}{}{#100}}%

\newcommand{\sebimo}{\@ifnextchar[{\basicsebimo}{\basicsebimo[59]}}%

\def\basicSebimo[#1]#2{\sdcase{\SEBIMO}{#2}{}{#100}}%

\newcommand{\Sebimo}{\@ifnextchar[{\basicSebimo}{\basicSebimo[59]}}%

\def\basicsebimO[#1]#2{\sdcase{\SEBIMO}{}{#2}{#100}}%

\newcommand{\sebimO}{\@ifnextchar[{\basicsebimO}{\basicsebimO[59]}}%

\def\basicseiso[#1]{\sdcase{\SEAR}{\hspace{-2pt}\isosign}{}{#100}}%

\newcommand{\seiso}{\@ifnextchar[{\basicseiso}{\basicseiso[59]}}%

\def\basicSeiso[#1]#2{\sdcase{\SEAR}{#2}{\isosign}{#100}}%

\newcommand{\Seiso}{\@ifnextchar[{\basicSeiso}{\basicSeiso[59]}}%

\def\basicseisO[#1]#2{\sdcase{\SEAR}{\hspace{-2pt}\isosign}{#2}{#100}}%

\newcommand{\seisO}{\@ifnextchar[{\basicseisO}{\basicseisO[59]}}%

\def\basicseeql[#1]{\sdcase{\SEEQL}{}{}{#100}}%

\newcommand{\seeql}{\@ifnextchar[{\basicseeql}{\basicseeql[59]}}%

\def\basicSeeql[#1]#2{\sdcase{\SEEQL}{#2}{}{#100}}%

\newcommand{\Seeql}{\@ifnextchar[{\basicSeeql}{\basicSeeql[59]}}%

\def\basicseeqL[#1]#2{\sdcase{\SEEQL}{}{#2}{#100}}%

\newcommand{\seeqL}{\@ifnextchar[{\basicseeqL}{\basicseeqL[59]}}%

\def\basicseline[#1]{\sdcase{\SELINE}{}{}{#100}}%

\newcommand{\seline}{\@ifnextchar[{\basicseline}{\basicseline[59]}}%

\def\basicSeline[#1]#2{\sdcase{\SELINE}{#2}{}{#100}}%

\newcommand{\Seline}{\@ifnextchar[{\basicSeline}{\basicSeline[59]}}%

\def\basicselinE[#1]#2{\sdcase{\SELINE}{}{#2}{#100}}%

\newcommand{\selinE}{\@ifnextchar[{\basicselinE}{\basicselinE[59]}}%

\def\basicsebiar[#1]{\sdbicase{\SEBIAR}{}{}{#100}}%

\newcommand{\sebiar}{\@ifnextchar[{\basicsebiar}{\basicsebiar[59]}}%

\def\basicSebiar[#1]#2#3{\sdbicase{\SEBIAR}{#2}{#3}{#100}}%

\newcommand{\Sebiar}{\@ifnextchar[{\basicSebiar}{\basicSebiar[59]}}%

\def\basicseadjar[#1]{\sdbicase{\SEADJAR}{}{}{#100}}%

\newcommand{\seadjar}{\@ifnextchar[{\basicseadjar}{\basicseadjar[59]}}%

\def\basicSeadjar[#1]#2#3{\sdbicase{\SEADJAR}{#2}{#3}{#100}}%

\newcommand{\Seadjar}{\@ifnextchar[{\basicSeadjar}{\basicSeadjar[59]}}%

\def\basicsebidist[#1]{\sdbicase{\SEBIDIST}{}{}{#100}}%

\newcommand{\sebidist}{\@ifnextchar[{\basicsebidist}{\basicsebidist[59]}}%

\def\basicSebidist[#1]#2#3{\sdbicase{\SEBIDIST}{#2}{#3}{#100}}%

\newcommand{\Sebidist}{\@ifnextchar[{\basicSebidist}{\basicSebidist[59]}}%

\def\basicseadjdist[#1]{\sdbicase{\SEADJDIST}{}{}{#100}}%

\newcommand{\seadjdist}{\@ifnextchar[{\basicseadjdist}{\basicseadjdist[59]}}%

\def\basicSeadjdist[#1]#2#3{\sdbicase{\SEADJDIST}{#2}{#3}{#100}}%

\newcommand{\Seadjdist}{\@ifnextchar[{\basicSeadjdist}{\basicSeadjdist[59]}}%

\newcommand{\NWAR}[1]{%
\Y=#1%
\divide\Y by 2%
\begin{picture}(0,0)%
\put(\Y,-\Y){\line(-1,1){#1}}%
\put(-\Y,\Y){\nwhead}%
\end{picture}}%

\newcommand{\NWDIST}[1]{%
\Y=#1%
\divide\Y by 2%
\begin{picture}(0,0)%
\put(\Y,-\Y){\line(-1,1){#1}}%
\put(-\Y,\Y){\nwhead}%
\put(0,0){\distsign}%
\end{picture}}%

\newcommand{\NWDOTAR}[1]%
{\truex{100}\truey{212}%
\Y=#1%
\divide\Y by 2%
\NUMBEROFDOTS=#1%
\divide\NUMBEROFDOTS by \value{y}%
\advance\NUMBEROFDOTS by 1%
\begin{picture}(0,0)%
\multiput(\Y,-\Y)(-\value{y},\value{y}){\NUMBEROFDOTS}%
{\circle*{\value{x}}}%
\put(-\Y,\Y){\nwhead}%
\end{picture}}%

\newcommand{\NWMONO}[1]{%
\Y=#1%
\divide \Y by 2%
\Truetail%
\bimolength=#1%
\advance\bimolength by -\Truemonotail%
\monolength=\bimolength%
\advance\monolength by -\Y%
\begin{picture}(0,0)%
\put(\monolength,-\monolength){\line(-1,1){\bimolength}}%
\put(\monolength,-\monolength){\nwhead}%
\put(-\Y,\Y){\nwhead}%
\end{picture}}%

\newcommand{\NWEPI}[1]{%
\Y=#1%
\divide\Y by 2%
\Truehead%
\bimolength=#1%
\advance\bimolength by -\Trueepihead%
\epilength=\bimolength%
\advance\epilength by -\Y%
\begin{picture}(0,0)%
\put(\Y,-\Y){\line(-1,1){\bimolength}}%
\put(-\epilength,\epilength){\nwhead}%
\put(-\Y,\Y){\nwhead}%
\end{picture}}%

\newcommand{\NWBIMO}[1]{%
\Y=#1%
\divide\Y by 2%
\Truetail\Truehead%
\bimolength=#1%
\advance\bimolength by -\Truemonotail%
\monolength=\bimolength%
\advance\monolength by -\Y%
\advance\bimolength by -\Trueepihead%
\epilength=\bimolength%
\advance\epilength by -\monolength%
\begin{picture}(0,0)%
\put(\monolength,-\monolength){\line(-1,1){\bimolength}}%
\put(\monolength,-\monolength){\nwhead}%
\put(-\epilength,\epilength){\nwhead}%
\put(-\Y,\Y){\nwhead}%
\end{picture}}%

\newcommand{\NWBIAR}[1]{%
\Y=#1%
\divide\Y by 2%
\begin{picture}(0,0)%
\put(\Y,-\Y){\begin{picture}(0,0)%
\truex{247}%
\put(-\value{x},-\value{x}){\line(-1,1){#1}}%
\put(\value{x},\value{x}){\line(-1,1){#1}}%
\monolength=#1%
\advance\monolength by -\value{x}%
\epilength=#1%
\advance\epilength by \value{x}%
\put(-\monolength,\epilength){\nwhead}%
\put(-\epilength,\monolength){\nwhead}%
\end{picture}}\end{picture}}%

\newcommand{\NWBIDIST}[1]{%
\Y=#1%
\divide\Y by 2%
\begin{picture}(0,0)%
\put(\Y,-\Y){\begin{picture}(0,0)%
\truex{247}%
\monolength=#1%
\advance\monolength by -\value{x}%
\epilength=#1%
\advance\epilength by \value{x}%
\put(-\value{x},-\value{x}){\line(-1,1){#1}}%
\put(-\epilength,\monolength){\nwhead}%
\end{picture}}%
\put(\Y,-\Y){\begin{picture}(0,0)%
\truex{247}%
\monolength=#1%
\advance\monolength by \value{x}%
\epilength=#1%
\advance\epilength by -\value{x}%
\put(\value{x},\value{x}){\line(-1,1){#1}}%
\put(-\epilength,\monolength){\nwhead}%
\end{picture}}%
\put(-\value{x},-\value{x}){\distsign}%
\put(\value{x},\value{x}){\distsign}%
\end{picture}}%

\newcommand{\NWADJAR}[1]{%
\Y=#1%
\divide\Y by 2%
\begin{picture}(0,0)%
\put(\Y,-\Y){\begin{picture}(0,0)%
\truex{247}%
\monolength=#1%
\advance\monolength by -\value{x}%
\epilength=#1%
\advance\epilength by \value{x}%
\put(-\value{x},-\value{x}){\line(-1,1){#1}}%
\put(-\epilength,\monolength){\nwhead}%
\end{picture}}%
\put(-\Y,\Y){\begin{picture}(0,0)%
\truex{247}%
\monolength=#1%
\advance\monolength by -\value{x}%
\epilength=#1%
\advance\epilength by \value{x}%
\put(\value{x},\value{x}){\line(1,-1){#1}}%
\put(\epilength,-\monolength){\sehead}%
\end{picture}}\end{picture}}%

\newcommand{\NWADJDIST}[1]{%
\Y=#1%
\divide\Y by 2%
\begin{picture}(0,0)%
\put(\Y,-\Y){\begin{picture}(0,0)%
\truex{247}%
\monolength=#1%
\advance\monolength by -\value{x}%
\epilength=#1%
\advance\epilength by \value{x}%
\put(-\value{x},-\value{x}){\line(-1,1){#1}}%
\put(-\epilength,\monolength){\nwhead}%
\end{picture}}%
\put(-\Y,\Y){\begin{picture}(0,0)%
\truex{247}%
\monolength=#1%
\advance\monolength by -\value{x}%
\epilength=#1%
\advance\epilength by \value{x}%
\put(\value{x},\value{x}){\line(1,-1){#1}}%
\put(\epilength,-\monolength){\sehead}%
\end{picture}}%
\put(-\value{x},-\value{x}){\distsign}%
\put(\value{x},\value{x}){\distsign}%
\end{picture}}%

\def\basicnwar[#1]{\sdcase{\NWAR}{}{}{#100}}%

\newcommand{\nwar}{\@ifnextchar[{\basicnwar}{\basicnwar[59]}}%

\def\basicNwar[#1]#2{\sdcase{\NWAR}{#2}{}{#100}}%

\newcommand{\Nwar}{\@ifnextchar[{\basicNwar}{\basicNwar[59]}}%

\def\basicnwaR[#1]#2{\sdcase{\NWAR}{}{#2}{#100}}%

\newcommand{\nwaR}{\@ifnextchar[{\basicnwaR}{\basicnwaR[59]}}%

\def\basicnwdist[#1]{\sdcase{\NWDIST}{}{}{#100}}%

\newcommand{\nwdist}{\@ifnextchar[{\basicnwdist}{\basicnwdist[59]}}%

\def\basicNwdist[#1]#2{\sdcase{\NWDIST}{#2}{}{#100}}%

\newcommand{\Nwdist}{\@ifnextchar[{\basicNwdist}{\basicNwdist[59]}}%

\def\basicnwdisT[#1]#2{\sdcase{\NWDIST}{}{#2}{#100}}%

\newcommand{\nwdisT}{\@ifnextchar[{\basicnwdisT}{\basicnwdisT[59]}}%

\def\basicnwdotar[#1]{\sdcase{\NWDOTAR}{}{}{#100}}%

\newcommand{\nwdotar}{\@ifnextchar[{\basicnwdotar}{\basicnwdotar[59]}}%

\def\basicNwdotar[#1]#2{\sdcase{\NWDOTAR}{#2}{}{#100}}%

\newcommand{\Nwdotar}{\@ifnextchar[{\basicNwdotar}{\basicNwdotar[59]}}%

\def\basicnwdotaR[#1]#2{\sdcase{\NWDOTAR}{}{#2}{#100}}%

\newcommand{\nwdotaR}{\@ifnextchar[{\basicnwdotaR}{\basicnwdotaR[59]}}%

\def\basicnwmono[#1]{\sdcase{\NWMONO}{}{}{#100}}%

\newcommand{\nwmono}{\@ifnextchar[{\basicnwmono}{\basicnwmono[59]}}%

\def\basicNwmono[#1]#2{\sdcase{\NWMONO}{#2}{}{#100}}%

\newcommand{\Nwmono}{\@ifnextchar[{\basicNwmono}{\basicNwmono[59]}}%

\def\basicnwmonO[#1]#2{\sdcase{\NWMONO}{}{#2}{#100}}%

\newcommand{\nwmonO}{\@ifnextchar[{\basicnwmonO}{\basicnwmonO[59]}}%

\def\basicnwepi[#1]{\sdcase{\NWEPI}{}{}{#100}}%

\newcommand{\nwepi}{\@ifnextchar[{\basicnwepi}{\basicnwepi[59]}}%

\def\basicNwepi[#1]#2{\sdcase{\NWEPI}{#2}{}{#100}}%

\newcommand{\Nwepi}{\@ifnextchar[{\basicNwepi}{\basicNwepi[59]}}%

\def\basicnwepI[#1]#2{\sdcase{\NWEPI}{}{#2}{#100}}%

\newcommand{\nwepI}{\@ifnextchar[{\basicnwepI}{\basicnwepI[59]}}%

\def\basicnwbimo[#1]{\sdcase{\NWBIMO}{}{}{#100}}%

\newcommand{\nwbimo}{\@ifnextchar[{\basicnwbimo}{\basicnwbimo[59]}}%

\def\basicNwbimo[#1]#2{\sdcase{\NWBIMO}{#2}{}{#100}}%

\newcommand{\Nwbimo}{\@ifnextchar[{\basicNwbimo}{\basicNwbimo[59]}}%

\def\basicnwbimO[#1]#2{\sdcase{\NWBIMO}{}{#2}{#100}}%

\newcommand{\nwbimO}{\@ifnextchar[{\basicnwbimO}{\basicnwbimO[59]}}%

\def\basicnwiso[#1]{\sdcase{\NWAR}{\hspace{-2pt}\isosign}{}{#100}}%

\newcommand{\nwiso}{\@ifnextchar[{\basicnwiso}{\basicnwiso[59]}}%

\def\basicNwiso[#1]#2{\sdcase{\NWAR}{#2}{\isosign}{#100}}%

\newcommand{\Nwiso}{\@ifnextchar[{\basicNwiso}{\basicNwiso[59]}}%

\def\basicnwisO[#1]#2{\sdcase{\NWAR}{\hspace{-2pt}\isosign}{#2}{#100}}%

\newcommand{\nwisO}{\@ifnextchar[{\basicnwisO}{\basicnwisO[59]}}%

\def\basicnwbiar[#1]{\sdbicase{\NWBIAR}{}{}{#100}}%

\newcommand{\nwbiar}{\@ifnextchar[{\basicnwbiar}{\basicnwbiar[59]}}%

\def\basicNwbiar[#1]#2#3{\sdbicase{\NWBIAR}{#2}{#3}{#100}}%

\newcommand{\Nwbiar}{\@ifnextchar[{\basicNwbiar}{\basicNwbiar[59]}}%

\def\basicnwadjar[#1]{\sdbicase{\NWADJAR}{}{}{#100}}%

\newcommand{\nwadjar}{\@ifnextchar[{\basicnwadjar}{\basicnwadjar[59]}}%

\def\basicNwadjar[#1]#2#3{\sdbicase{\NWADJAR}{#2}{#3}{#100}}%

\newcommand{\Nwadjar}{\@ifnextchar[{\basicNwadjar}{\basicNwadjar[59]}}%

\def\basicnwbidist[#1]{\sdbicase{\NWBIDIST}{}{}{#100}}%

\newcommand{\nwbidist}{\@ifnextchar[{\basicnwbidist}{\basicnwbidist[59]}}%

\def\basicNwbidist[#1]#2#3{\sdbicase{\NWBIDIST}{#2}{#3}{#100}}%

\newcommand{\Nwbidist}{\@ifnextchar[{\basicNwbidist}{\basicNwbidist[59]}}%

\def\basicnwadjdist[#1]{\sdbicase{\NWADJDIST}{}{}{#100}}%

\newcommand{\nwadjdist}{\@ifnextchar[{\basicnwadjdist}{\basicnwadjdist[59]}}%

\def\basicNwadjdist[#1]#2#3{\sdbicase{\NWADJDIST}{#2}{#3}{#100}}%

\newcommand{\Nwadjdist}{\@ifnextchar[{\basicNwadjdist}{\basicNwadjdist[59]}}%

\newcommand{\ENEAR}[3]{%
\Y=#3%
\divide\Y by 2%
\Z=\Y%
\divide\Z by 2%
\begin{picture}(0,0)%
\put(-\Y,-\Z){\line(2,1){#3}}%
\put(\Y,\Z){\enehead}%
\truex{200}\truey{800}\truez{600}%
\put(-\value{x},\value{x}){\makebox(0,\value{z})[r]{${#1}$}}%
\put(\value{x},-\value{y}){\makebox(0,\value{z})[l]{${#2}$}}%
\end{picture}}%

\newcommand{\ENEDIST}[3]{%
\Y=#3%
\divide\Y by 2%
\Z=\Y%
\divide\Z by 2%
\begin{picture}(0,0)%
\put(-\Y,-\Z){\line(2,1){#3}}%
\put(\Y,\Z){\enehead}%
\put(0,0){\distsign}%
\truex{200}\truey{800}\truez{600}%
\put(-\value{x},\value{x}){\makebox(0,\value{z})[r]{${#1}$}}%
\put(\value{x},-\value{y}){\makebox(0,\value{z})[l]{${#2}$}}%
\end{picture}}%

\newcommand{\ENEDOTAR}[3]{%
\truex{100}\truey{268}\truez{134}%
\Y=#3%
\divide\Y by 2%
\Z=\Y%
\divide\Z by 2%
\NUMBEROFDOTS=#3%
\divide\NUMBEROFDOTS by \value{y}%
\advance\NUMBEROFDOTS by 1%
\begin{picture}(0,0)%
\multiput(-\Y,-\Z)(\value{y},\value{z}){\NUMBEROFDOTS}%
{\circle*{\value{x}}}%
\put(\Y,\Z){\enehead}%
\truex{200}\truey{800}\truez{600}%
\put(-\value{x},\value{x}){\makebox(0,\value{z})[r]{${#1}$}}%
\put(\value{x},-\value{y}){\makebox(0,\value{z})[l]{${#2}$}}%
\end{picture}}%

\newcommand{\ENEMONO}[3]{%
\Y=#3%
\divide\Y by 2%
\Z=\Y%
\divide\Z by 2%
\TrueTail%
\bimolength=#3%
\advance\bimolength by -\TrueMonoTail%
\monolength=\bimolength%
\advance\monolength by -\Y%
\secondmonolength=\monolength%
\divide\secondmonolength by 2%
\begin{picture}(0,0)%
\put(-\monolength,-\secondmonolength){\line(2,1){\bimolength}}%
\put(-\monolength,-\secondmonolength){\enehead}%
\put(\Y,\Z){\enehead}%
\truex{200}\truey{800}\truez{600}%
\put(-\value{x},\value{x}){\makebox(0,\value{z})[r]{${#1}$}}%
\put(\value{x},-\value{y}){\makebox(0,\value{z})[l]{${#2}$}}%
\end{picture}}%

\newcommand{\ENEEPI}[3]{%
\Y=#3%
\divide\Y by 2%
\Z=\Y%
\divide\Z by 2%
\TrueHead%
\bimolength=#3%
\advance\bimolength by -\TrueEpiHead%
\epilength=\bimolength%
\advance\epilength by -\Y%
\secondepilength=\epilength%
\divide\secondepilength by 2%
\begin{picture}(0,0)%
\put(-\Y,-\Z){\line(2,1){\bimolength}}%
\put(\epilength,\secondepilength){\enehead}%
\put(\Y,\Z){\enehead}%
\truex{200}\truey{800}\truez{600}%
\put(-\value{x},\value{x}){\makebox(0,\value{z})[r]{${#1}$}}%
\put(\value{x},-\value{y}){\makebox(0,\value{z})[l]{${#2}$}}%
\end{picture}}%

\newcommand{\ENEBIMO}[3]{%
\Y=#3%
\divide\Y by 2%
\Z=\Y%
\divide\Z by 2%
\TrueTail\TrueHead%
\bimolength=#3%
\advance\bimolength by -\TrueMonoTail%
\monolength=\bimolength%
\advance\monolength by -\Y%
\advance\bimolength by -\TrueEpiHead%
\epilength=\bimolength%
\advance\epilength by -\monolength%
\secondmonolength=\monolength%
\divide\secondmonolength by 2%
\secondepilength=\epilength%
\divide\secondepilength by 2%
\begin{picture}(0,0)%
\put(-\monolength,-\secondmonolength){\line(2,1){\bimolength}}%
\put(-\monolength,-\secondmonolength){\enehead}%
\put(\epilength,\secondepilength){\enehead}%
\put(\Y,\Z){\enehead}%
\truex{200}\truey{800}\truez{600}%
\put(-\value{x},\value{x}){\makebox(0,\value{z})[r]{${#1}$}}%
\put(\value{x},-\value{y}){\makebox(0,\value{z})[l]{${#2}$}}%
\end{picture}}%

\newcommand{\ENEEQL}[3]{%
\Y=#3%
\divide\Y by 2%
\Z=\Y%
\divide\Z by 2%
\begin{picture}(0,0)%
\put(-\Y,-\Z){\begin{picture}(0,0)%
\truex{44}\truey{89}%
\put(-\value{x},\value{y}){\line(2,1){#3}}%
\put(\value{x},-\value{y}){\line(2,1){#3}}%
\end{picture}}%
\truex{200}\truey{800}\truez{600}%
\put(-\value{x},\value{x}){\makebox(0,\value{z})[r]{${#1}$}}%
\put(\value{x},-\value{y}){\makebox(0,\value{z})[l]{${#2}$}}%
\end{picture}}%

\newcommand{\ENEBIAR}[3]{%
\Y=#3%
\divide\Y by 2%
\Z=\Y%
\divide\Z by 2%
\begin{picture}(0,0)%
\put(-\Y,-\Z){\begin{picture}(0,0)%
\truex{156}\truey{313}%
\put(-\value{x},\value{y}){\line(2,1){#3}}%
\put(\value{x},-\value{y}){\line(2,1){#3}}%
\monolength=#3%
\advance\monolength by -\value{x}%
\epilength=#3%
\advance\epilength by \value{x}%
\secondmonolength=\Y%
\advance\secondmonolength by -\value{y}%
\secondepilength=\Y%
\advance\secondepilength by \value{y}%
\put(\monolength,\secondepilength){\enehead}%
\put(\epilength,\secondmonolength){\enehead}%
\end{picture}}
\truex{300}\truey{1000}\truez{600}%
\put(-\value{x},\value{x}){\makebox(0,\value{z})[r]{${#1}$}}%
\put(\value{x},-\value{y}){\makebox(0,\value{z})[l]{${#2}$}}%
\end{picture}}%

\newcommand{\ENEBIDIST}[3]{%
\Y=#3%
\divide\Y by 2%
\Z=\Y%
\divide\Z by 2%
\begin{picture}(0,0)%
\truex{156}\truey{313}\truez{400}%
\put(-\Y,-\Z){\begin{picture}(0,0)%
\put(-\value{x},\value{y}){\line(2,1){#3}}%
\put(\value{x},-\value{y}){\line(2,1){#3}}%
\monolength=#3%
\advance\monolength by -\value{x}%
\epilength=#3%
\advance\epilength by \value{x}%
\secondmonolength=\Y%
\advance\secondmonolength by -\value{y}%
\secondepilength=\Y%
\advance\secondepilength by \value{y}%
\put(\monolength,\secondepilength){\enehead}%
\put(\epilength,\secondmonolength){\enehead}%
\end{picture}}
\put(-\value{x},\value{y}){\distsign}%
\put(\value{x},-\value{y}){\distsign}%
\truex{300}\truey{1000}\truez{600}%
\put(-\value{x},\value{x}){\makebox(0,\value{z})[r]{${#1}$}}%
\put(\value{x},-\value{y}){\makebox(0,\value{z})[l]{${#2}$}}%
\end{picture}}%

\newcommand{\ENEADJAR}[3]{%
\Y=#3%
\divide\Y by 2%
\Z=\Y%
\divide\Z by 2%
\begin{picture}(0,0)%
\put(-\Y,-\Z){\begin{picture}(0,0)%
\truex{156}\truey{313}%
\monolength=#3%
\advance\monolength by -\value{x}%
\epilength=#3%
\advance\epilength by \value{x}%
\secondmonolength=\Y%
\advance\secondmonolength by -\value{y}%
\secondepilength=\Y%
\advance\secondepilength by \value{y}%
\put(\value{x},-\value{y}){\line(2,1){#3}}%
\put(\epilength,\secondmonolength){\enehead}%
\put(\monolength,\secondepilength){\line(-2,-1){#3}}%
\put(-\value{x},\value{y}){\wswhead}%
\end{picture}}
\truex{300}\truey{1000}\truez{600}%
\put(-\value{x},\value{x}){\makebox(0,\value{z})[r]{${#1}$}}%
\put(\value{x},-\value{y}){\makebox(0,\value{z})[l]{${#2}$}}%
\end{picture}}%

\newcommand{\ENEADJDIST}[3]{%
\Y=#3%
\divide\Y by 2%
\Z=\Y%
\divide\Z by 2%
\begin{picture}(0,0)%
\truex{156}\truey{313}\truez{400}%
\put(-\Y,-\Z){\begin{picture}(0,0)%
\monolength=#3%
\advance\monolength by -\value{x}%
\epilength=#3%
\advance\epilength by \value{x}%
\secondmonolength=\Y%
\advance\secondmonolength by -\value{y}%
\secondepilength=\Y%
\advance\secondepilength by \value{y}%
\put(\value{x},-\value{y}){\line(2,1){#3}}%
\put(\epilength,\secondmonolength){\enehead}%
\put(\monolength,\secondepilength){\line(-2,-1){#3}}%
\put(-\value{x},\value{y}){\wswhead}%
\end{picture}}
\put(-\value{x},\value{y}){\distsign}%
\put(\value{x},-\value{y}){\distsign}%
\truex{300}\truey{1000}\truez{600}%
\put(-\value{x},\value{x}){\makebox(0,\value{z})[r]{${#1}$}}%
\put(\value{x},-\value{y}){\makebox(0,\value{z})[l]{${#2}$}}%
\end{picture}}%

\def\basicenear[#1]{\ENEAR{}{}{#100}}%

\newcommand{\enear}{\@ifnextchar[{\basicenear}{\basicenear[133]}}%

\def\basicEnear[#1]#2{\ENEAR{#2}{}{#100}}%

\newcommand{\Enear}{\@ifnextchar[{\basicEnear}{\basicEnear[133]}}%

\def\basiceneaR[#1]#2{\ENEAR{}{#2}{#100}}%

\newcommand{\eneaR}{\@ifnextchar[{\basiceneaR}{\basiceneaR[133]}}%

\def\basicenedist[#1]{\ENEDIST{}{}{#100}}%

\newcommand{\enedist}{\@ifnextchar[{\basicenedist}{\basicenedist[133]}}%

\def\basicEnedist[#1]#2{\ENEDIST{#2}{}{#100}}%

\newcommand{\Enedist}{\@ifnextchar[{\basicEnedist}{\basicEnedist[133]}}%

\def\basicenedisT[#1]#2{\ENEDIST{}{#2}{#100}}%

\newcommand{\enedisT}{\@ifnextchar[{\basicenedisT}{\basicenedisT[133]}}%

\def\basicenedotar[#1]{\ENEDOTAR{}{}{#100}}%

\newcommand{\enedotar}{\@ifnextchar[{\basicenedotar}{\basicenedotar[133]}}%

\def\basicEnedotar[#1]#2{\ENEDOTAR{#2}{}{#100}}%

\newcommand{\Enedotar}{\@ifnextchar[{\basicEnedotar}{\basicEnedotar[133]}}%

\def\basicenedotaR[#1]#2{\ENEDOTAR{}{#2}{#100}}%

\newcommand{\enedotaR}{\@ifnextchar[{\basicenedotaR}{\basicenedotaR[133]}}%

\def\basicenemono[#1]{\ENEMONO{}{}{#100}}%

\newcommand{\enemono}{\@ifnextchar[{\basicenemono}{\basicenemono[133]}}%

\def\basicEnemono[#1]#2{\ENEMONO{#2}{}{#100}}%

\newcommand{\Enemono}{\@ifnextchar[{\basicEnemono}{\basicEnemono[133]}}%

\def\basicenemonO[#1]#2{\ENEMONO{}{#2}{#100}}%

\newcommand{\enemonO}{\@ifnextchar[{\basicenemonO}{\basicenemonO[133]}}%

\def\basiceneepi[#1]{\ENEEPI{}{}{#100}}%

\newcommand{\eneepi}{\@ifnextchar[{\basiceneepi}{\basiceneepi[133]}}%

\def\basicEneepi[#1]#2{\ENEEPI{#2}{}{#100}}%

\newcommand{\Eneepi}{\@ifnextchar[{\basicEneepi}{\basicEneepi[133]}}%

\def\basiceneepI[#1]#2{\ENEEPI{}{#2}{#100}}%

\newcommand{\eneepI}{\@ifnextchar[{\basiceneepI}{\basiceneepI[133]}}%

\def\basicenebimo[#1]{\ENEBIMO{}{}{#100}}%

\newcommand{\enebimo}{\@ifnextchar[{\basicenebimo}{\basicenebimo[133]}}%

\def\basicEnebimo[#1]#2{\ENEBIMO{#2}{}{#100}}%

\newcommand{\Enebimo}{\@ifnextchar[{\basicEnebimo}{\basicEnebimo[133]}}%

\def\basicenebimO[#1]#2{\ENEBIMO{}{#2}{#100}}%

\newcommand{\enebimO}{\@ifnextchar[{\basicenebimO}{\basicenebimO[133]}}%

\def\basiceneiso[#1]{\ENEAR{\isosign}{}{#100}}%

\newcommand{\eneiso}{\@ifnextchar[{\basiceneiso}{\basiceneiso[133]}}%

\def\basicEneiso[#1]#2{\ENEAR{#2}{\isosign}{#100}}%

\newcommand{\Eneiso}{\@ifnextchar[{\basicEneiso}{\basicEneiso[133]}}%

\def\basiceneisO[#1]#2{\ENEAR{\isosign}{#2}{#100}}%

\newcommand{\eneisO}{\@ifnextchar[{\basiceneisO}{\basiceneisO[133]}}%

\def\basiceneeql[#1]{\ENEEQL{}{}{#100}}%

\newcommand{\eneeql}{\@ifnextchar[{\basiceneeql}{\basiceneeql[133]}}%

\def\basicEneeql[#1]#2{\ENEEQL{#2}{}{#100}}%

\newcommand{\Eneeql}{\@ifnextchar[{\basicEneeql}{\basicEneeql[133]}}%

\def\basiceneeqL[#1]#2{\ENEEQL{}{#2}{#100}}%

\newcommand{\eneeqL}{\@ifnextchar[{\basiceneeqL}{\basiceneeqL[133]}}%

\def\basicenebiar[#1]{\ENEBIAR{}{}{#100}}%

\newcommand{\enebiar}{\@ifnextchar[{\basicenebiar}{\basicenebiar[133]}}%

\def\basicEnebiar[#1]#2#3{\ENEBIAR{#2}{#3}{#100}}%

\newcommand{\Enebiar}{\@ifnextchar[{\basicEnebiar}{\basicEnebiar[133]}}%

\def\basicenebidist[#1]{\ENEBIDIST{}{}{#100}}%

\newcommand{\enebidist}{\@ifnextchar[{\basicenebidist}{\basicenebidist[133]}}%

\def\basicEnebidist[#1]#2#3{\ENEBIDIST{#2}{#3}{#100}}%

\newcommand{\Enebidist}{\@ifnextchar[{\basicEnebidist}{\basicEnebidist[133]}}%

\def\basiceneadjar[#1]{\ENEADJAR{}{}{#100}}%

\newcommand{\eneadjar}{\@ifnextchar[{\basiceneadjar}{\basiceneadjar[133]}}%

\def\basicEneadjar[#1]#2#3{\ENEADJAR{#2}{#3}{#100}}%

\newcommand{\Eneadjar}{\@ifnextchar[{\basicEneadjar}{\basicEneadjar[133]}}%

\def\basiceneadjdist[#1]{\ENEADJDIST{}{}{#100}}%

\newcommand{\eneadjdist}{\@ifnextchar[{\basiceneadjdist}{\basiceneadjdist[13
3]}}%

\def\basicEneadjdist[#1]#2#3{\ENEADJDIST{#2}{#3}{#100}}%

\newcommand{\Eneadjdist}{\@ifnextchar[{\basicEneadjdist}{\basicEneadjdist[13
3]}}%

\newcommand{\ESEAR}[3]{%
\Y=#3%
\divide\Y by 2%
\Z=\Y%
\divide\Z by 2%
\begin{picture}(0,0)%
\put(-\Y,\Z){\line(2,-1){#3}}%
\put(\Y,-\Z){\esehead}%
\truex{200}\truey{800}\truez{600}%
\put(\value{x},\value{x}){\makebox(0,\value{z})[l]{${#1}$}}%
\put(-\value{x},-\value{y}){\makebox(0,\value{z})[r]{${#2}$}}%
\end{picture}}%

\newcommand{\ESEDIST}[3]{%
\Y=#3%
\divide\Y by 2%
\Z=\Y%
\divide\Z by 2%
\begin{picture}(0,0)%
\put(-\Y,\Z){\line(2,-1){#3}}%
\put(\Y,-\Z){\esehead}%
\put(0,0){\distsign}%
\truex{200}\truey{800}\truez{600}%
\put(\value{x},\value{x}){\makebox(0,\value{z})[l]{${#1}$}}%
\put(-\value{x},-\value{y}){\makebox(0,\value{z})[r]{${#2}$}}%
\end{picture}}%

\newcommand{\ESEDOTAR}[3]{%
\truex{100}\truey{268}\truez{134}%
\Y=#3%
\divide\Y by 2%
\Z=\Y%
\divide\Z by 2%
\NUMBEROFDOTS=#3%
\divide\NUMBEROFDOTS by \value{y}%
\advance\NUMBEROFDOTS by 1%
\begin{picture}(0,0)%
\multiput(-\Y,\Z)(\value{y},-\value{z}){\NUMBEROFDOTS}%
{\circle*{\value{x}}}%
\put(\Y,-\Z){\esehead}%
\truex{200}\truey{800}\truez{600}%
\put(\value{x},\value{x}){\makebox(0,\value{z})[l]{${#1}$}}%
\put(-\value{x},-\value{y}){\makebox(0,\value{z})[r]{${#2}$}}%
\end{picture}}%

\newcommand{\ESEMONO}[3]{%
\Y=#3%
\divide\Y by 2%
\Z=\Y%
\divide\Z by 2%
\TrueTail%
\bimolength=#3%
\advance\bimolength by -\TrueMonoTail%
\monolength=\bimolength%
\advance\monolength by -\Y%
\secondmonolength=\monolength%
\divide\secondmonolength by 2%
\begin{picture}(0,0)%
\put(-\monolength,\secondmonolength){\line(2,-1){\bimolength}}%
\put(-\monolength,\secondmonolength){\esehead}%
\put(\Y,-\Z){\esehead}%
\truex{200}\truey{800}\truez{600}%
\put(\value{x},\value{x}){\makebox(0,\value{z})[l]{${#1}$}}%
\put(-\value{x},-\value{y}){\makebox(0,\value{z})[r]{${#2}$}}%
\end{picture}}%

\newcommand{\ESEEPI}[3]{%
\Y=#3%
\divide\Y by 2%
\Z=\Y%
\divide\Z by 2%
\TrueHead%
\bimolength=#3%
\advance\bimolength by -\TrueEpiHead%
\epilength=\bimolength%
\advance\epilength by -\Y%
\secondepilength=\epilength%
\divide\secondepilength by 2%
\begin{picture}(0,0)%
\put(-\Y,\Z){\line(2,-1){\bimolength}}%
\put(\epilength,-\secondepilength){\esehead}%
\put(\Y,-\Z){\esehead}%
\truex{200}\truey{800}\truez{600}%
\put(\value{x},\value{x}){\makebox(0,\value{z})[l]{${#1}$}}%
\put(-\value{x},-\value{y}){\makebox(0,\value{z})[r]{${#2}$}}%
\end{picture}}%

\newcommand{\ESEBIMO}[3]{%
\Y=#3%
\divide\Y by 2%
\Z=\Y%
\divide\Z by 2%
\TrueTail\TrueHead%
\bimolength=#3%
\advance\bimolength by -\TrueMonoTail%
\monolength=\bimolength%
\advance\monolength by -\Y%
\advance\bimolength by -\TrueEpiHead%
\epilength=\bimolength%
\advance\epilength by -\monolength%
\secondmonolength=\monolength%
\divide\secondmonolength by 2%
\secondepilength=\epilength%
\divide\secondepilength by 2%
\begin{picture}(0,0)%
\put(-\monolength,\secondmonolength){\line(2,-1){\bimolength}}%
\put(-\monolength,\secondmonolength){\esehead}%
\put(\epilength,-\secondepilength){\esehead}%
\put(\Y,-\Z){\esehead}%
\truex{200}\truey{800}\truez{600}%
\put(\value{x},\value{x}){\makebox(0,\value{z})[l]{${#1}$}}%
\put(-\value{x},-\value{y}){\makebox(0,\value{z})[r]{${#2}$}}%
\end{picture}}%

\newcommand{\ESEEQL}[3]{%
\Y=#3%
\divide\Y by 2%
\Z=\Y%
\divide\Z by 2%
\begin{picture}(0,0)%
\put(-\Y,\Z){\begin{picture}(0,0)%
\truex{44}\truey{89}%
\put(-\value{x},-\value{y}){\line(2,-1){#3}}%
\put(\value{x},\value{y}){\line(2,-1){#3}}%
\end{picture}}%
\truex{200}\truey{800}\truez{600}%
\put(\value{x},\value{x}){\makebox(0,\value{z})[l]{${#1}$}}%
\put(-\value{x},-\value{y}){\makebox(0,\value{z})[r]{${#2}$}}%
\end{picture}}%

\newcommand{\ESEBIAR}[3]{%
\Y=#3%
\divide\Y by 2%
\Z=\Y%
\divide\Z by 2%
\begin{picture}(0,0)%
\put(-\Y,\Z){\begin{picture}(0,0)%
\truex{156}\truey{313}%
\put(-\value{x},-\value{y}){\line(2,-1){#3}}%
\put(\value{x},\value{y}){\line(2,-1){#3}}%
\monolength=#3%
\advance\monolength by -\value{x}%
\epilength=#3%
\advance\epilength by \value{x}%
\secondmonolength=\Y%
\advance\secondmonolength by -\value{y}%
\secondepilength=\Y%
\advance\secondepilength by \value{y}%
\put(\monolength,-\secondepilength){\esehead}%
\put(\epilength,-\secondmonolength){\esehead}%
\end{picture}}
\truex{400}\truey{1000}\truez{600}%
\put(\value{x},\value{x}){\makebox(0,\value{z})[l]{${#1}$}}%
\put(-\value{x},-\value{y}){\makebox(0,\value{z})[r]{${#2}$}}%
\end{picture}}%

\newcommand{\ESEBIDIST}[3]{%
\Y=#3%
\divide\Y by 2%
\Z=\Y%
\divide\Z by 2%
\begin{picture}(0,0)%
\truex{156}\truey{313}\truez{400}%
\put(-\Y,\Z){\begin{picture}(0,0)%
\put(-\value{x},-\value{y}){\line(2,-1){#3}}%
\put(\value{x},\value{y}){\line(2,-1){#3}}%
\monolength=#3%
\advance\monolength by -\value{x}%
\epilength=#3%
\advance\epilength by \value{x}%
\secondmonolength=\Y%
\advance\secondmonolength by -\value{y}%
\secondepilength=\Y%
\advance\secondepilength by \value{y}%
\put(\monolength,-\secondepilength){\esehead}%
\put(\epilength,-\secondmonolength){\esehead}%
\end{picture}}
\put(\value{x},\value{y}){\distsign}%
\put(-\value{x},-\value{y}){\distsign}%
\truex{400}\truey{1000}\truez{600}%
\put(\value{x},\value{x}){\makebox(0,\value{z})[l]{${#1}$}}%
\put(-\value{x},-\value{y}){\makebox(0,\value{z})[r]{${#2}$}}%
\end{picture}}%

\newcommand{\ESEADJAR}[3]{%
\Y=#3%
\divide\Y by 2%
\Z=\Y%
\divide\Z by 2%
\begin{picture}(0,0)%
\put(-\Y,\Z){\begin{picture}(0,0)%
\truex{156}\truey{313}%
\monolength=#3%
\advance\monolength by -\value{x}%
\epilength=#3%
\advance\epilength by \value{x}%
\secondmonolength=\Y%
\advance\secondmonolength by -\value{y}%
\secondepilength=\Y%
\advance\secondepilength by \value{y}%
\put(-\value{x},-\value{y}){\line(2,-1){#3}}%
\put(\monolength,-\secondepilength){\esehead}%
\put(\epilength,-\secondmonolength){\line(-2,1){#3}}%
\put(\value{x},\value{y}){\wnwhead}%
\end{picture}}
\truex{400}\truey{1000}\truez{600}%
\put(\value{x},\value{x}){\makebox(0,\value{z})[l]{${#1}$}}%
\put(-\value{x},-\value{y}){\makebox(0,\value{z})[r]{${#2}$}}%
\end{picture}}%

\newcommand{\ESEADJDIST}[3]{%
\Y=#3%
\divide\Y by 2%
\Z=\Y%
\divide\Z by 2%
\begin{picture}(0,0)%
\truex{156}\truey{313}\truez{400}%
\put(-\Y,\Z){\begin{picture}(0,0)%
\monolength=#3%
\advance\monolength by -\value{x}%
\epilength=#3%
\advance\epilength by \value{x}%
\secondmonolength=\Y%
\advance\secondmonolength by -\value{y}%
\secondepilength=\Y%
\advance\secondepilength by \value{y}%
\put(-\value{x},-\value{y}){\line(2,-1){#3}}%
\put(\monolength,-\secondepilength){\esehead}%
\put(\epilength,-\secondmonolength){\line(-2,1){#3}}%
\put(\value{x},\value{y}){\wnwhead}%
\end{picture}}
\put(\value{x},\value{y}){\distsign}%
\put(-\value{x},-\value{y}){\distsign}%
\truex{400}\truey{1000}\truez{600}%
\put(\value{x},\value{x}){\makebox(0,\value{z})[l]{${#1}$}}%
\put(-\value{x},-\value{y}){\makebox(0,\value{z})[r]{${#2}$}}%
\end{picture}}%

\def\basicesear[#1]{\ESEAR{}{}{#100}}%

\newcommand{\esear}{\@ifnextchar[{\basicesear}{\basicesear[133]}}%

\def\basicEsear[#1]#2{\ESEAR{#2}{}{#100}}%

\newcommand{\Esear}{\@ifnextchar[{\basicEsear}{\basicEsear[133]}}%

\def\basiceseaR[#1]#2{\ESEAR{}{#2}{#100}}%

\newcommand{\eseaR}{\@ifnextchar[{\basiceseaR}{\basiceseaR[133]}}%

\def\basicesedist[#1]{\ESEDIST{}{}{#100}}%

\newcommand{\esedist}{\@ifnextchar[{\basicesedist}{\basicesedist[133]}}%

\def\basicEsedist[#1]#2{\ESEDIST{#2}{}{#100}}%

\newcommand{\Esedist}{\@ifnextchar[{\basicEsedist}{\basicEsedist[133]}}%

\def\basicesedisT[#1]#2{\ESEDIST{}{#2}{#100}}%

\newcommand{\esedisT}{\@ifnextchar[{\basicesedisT}{\basicesedisT[133]}}%

\def\basicesedotar[#1]{\ESEDOTAR{}{}{#100}}%

\newcommand{\esedotar}{\@ifnextchar[{\basicesedotar}{\basicesedotar[133]}}%

\def\basicEsedotar[#1]#2{\ESEDOTAR{#2}{}{#100}}%

\newcommand{\Esedotar}{\@ifnextchar[{\basicEsedotar}{\basicEsedotar[133]}}%

\def\basicesedotaR[#1]#2{\ESEDOTAR{}{#2}{#100}}%

\newcommand{\esedotaR}{\@ifnextchar[{\basicesedotaR}{\basicesedotaR[133]}}%

\def\basicesemono[#1]{\ESEMONO{}{}{#100}}%

\newcommand{\esemono}{\@ifnextchar[{\basicesemono}{\basicesemono[133]}}%

\def\basicEsemono[#1]#2{\ESEMONO{#2}{}{#100}}%

\newcommand{\Esemono}{\@ifnextchar[{\basicEsemono}{\basicEsemono[133]}}%

\def\basicesemonO[#1]#2{\ESEMONO{}{#2}{#100}}%

\newcommand{\esemonO}{\@ifnextchar[{\basicesemonO}{\basicesemonO[133]}}%

\def\basiceseepi[#1]{\ESEEPI{}{}{#100}}%

\newcommand{\eseepi}{\@ifnextchar[{\basiceseepi}{\basiceseepi[133]}}%

\def\basicEseepi[#1]#2{\ESEEPI{#2}{}{#100}}%

\newcommand{\Eseepi}{\@ifnextchar[{\basicEseepi}{\basicEseepi[133]}}%

\def\basiceseepI[#1]#2{\ESEEPI{}{#2}{#100}}%

\newcommand{\eseepI}{\@ifnextchar[{\basiceseepI}{\basiceseepI[133]}}%

\def\basicesebimo[#1]{\ESEBIMO{}{}{#100}}%

\newcommand{\esebimo}{\@ifnextchar[{\basicesebimo}{\basicesebimo[133]}}%

\def\basicEsebimo[#1]#2{\ESEBIMO{#2}{}{#100}}%

\newcommand{\Esebimo}{\@ifnextchar[{\basicEsebimo}{\basicEsebimo[133]}}%

\def\basicesebimO[#1]#2{\ESEBIMO{}{#2}{#100}}%

\newcommand{\esebimO}{\@ifnextchar[{\basicesebimO}{\basicesebimO[133]}}%

\def\basiceseiso[#1]{\ESEAR{\isosign}{}{#100}}%

\newcommand{\eseiso}{\@ifnextchar[{\basiceseiso}{\basiceseiso[133]}}%

\def\basicEseiso[#1]#2{\ESEAR{#2}{\isosign}{#100}}%

\newcommand{\Eseiso}{\@ifnextchar[{\basicEseiso}{\basicEseiso[133]}}%

\def\basiceseisO[#1]#2{\ESEAR{\isosign}{#2}{#100}}%

\newcommand{\eseisO}{\@ifnextchar[{\basiceseisO}{\basiceseisO[133]}}%

\def\basiceseeql[#1]{\ESEEQL{}{}{#100}}%

\newcommand{\eseeql}{\@ifnextchar[{\basiceseeql}{\basiceseeql[133]}}%

\def\basicEseeql[#1]#2{\ESEEQL{#2}{}{#100}}%

\newcommand{\Eseeql}{\@ifnextchar[{\basicEseeql}{\basicEseeql[133]}}%

\def\basiceseeqL[#1]#2{\ESEEQL{}{#2}{#100}}%

\newcommand{\eseeqL}{\@ifnextchar[{\basiceseeqL}{\basiceseeqL[133]}}%

\def\basicesebiar[#1]{\ESEBIAR{}{}{#100}}%

\newcommand{\esebiar}{\@ifnextchar[{\basicesebiar}{\basicesebiar[133]}}%

\def\basicEsebiar[#1]#2#3{\ESEBIAR{#2}{#3}{#100}}%

\newcommand{\Esebiar}{\@ifnextchar[{\basicEsebiar}{\basicEsebiar[133]}}%

\def\basicesebidist[#1]{\ESEBIDIST{}{}{#100}}%

\newcommand{\esebidist}{\@ifnextchar[{\basicesebidist}{\basicesebidist[133]}}%

\def\basicEsebidist[#1]#2#3{\ESEBIDIST{#2}{#3}{#100}}%

\newcommand{\Esebidist}{\@ifnextchar[{\basicEsebidist}{\basicEsebidist[133]}}%

\def\basiceseadjar[#1]{\ESEADJAR{}{}{#100}}%

\newcommand{\eseadjar}{\@ifnextchar[{\basiceseadjar}{\basiceseadjar[133]}}%

\def\basicEseadjar[#1]#2#3{\ESEADJAR{#2}{#3}{#100}}%

\newcommand{\Eseadjar}{\@ifnextchar[{\basicEseadjar}{\basicEseadjar[133]}}%

\def\basiceseadjdist[#1]{\ESEADJDIST{}{}{#100}}%

\newcommand{\eseadjdist}{\@ifnextchar[{\basiceseadjdist}{\basiceseadjdist[13
3]}}%

\def\basicEseadjdist[#1]#2#3{\ESEADJDIST{#2}{#3}{#100}}%

\newcommand{\Eseadjdist}{\@ifnextchar[{\basicEseadjdist}{\basicEseadjdist[13
3]}}%

\newcommand{\WSWAR}[3]{%
\Y=#3%
\divide\Y by 2%
\Z=\Y%
\divide\Z by 2%
\begin{picture}(0,0)%
\put(\Y,\Z){\line(-2,-1){#3}}%
\put(-\Y,-\Z){\wswhead}%
\truex{200}\truey{800}\truez{600}%
\put(-\value{x},\value{x}){\makebox(0,\value{z})[r]{${#1}$}}%
\put(\value{x},-\value{y}){\makebox(0,\value{z})[l]{${#2}$}}%
\end{picture}}%

\newcommand{\WSWDIST}[3]{%
\Y=#3%
\divide\Y by 2%
\Z=\Y%
\divide\Z by 2%
\begin{picture}(0,0)%
\put(\Y,\Z){\line(-2,-1){#3}}%
\put(-\Y,-\Z){\wswhead}%
\truex{400}%
\put(0,0){\distsign}%
\truex{200}\truey{800}\truez{600}%
\put(-\value{x},\value{x}){\makebox(0,\value{z})[r]{${#1}$}}%
\put(\value{x},-\value{y}){\makebox(0,\value{z})[l]{${#2}$}}%
\end{picture}}%

\newcommand{\WSWDOTAR}[3]{%
\truex{100}\truey{268}\truez{134}%
\Y=#3%
\divide\Y by 2%
\Z=\Y%
\divide\Z by 2%
\NUMBEROFDOTS=#3%
\divide\NUMBEROFDOTS by \value{y}%
\advance\NUMBEROFDOTS by 1%
\begin{picture}(0,0)%
\multiput(\Y,\Z)(-\value{y},-\value{z}){\NUMBEROFDOTS}%
{\circle*{\value{x}}}%
\put(-\Y,-\Z){\wswhead}%
\truex{200}\truey{800}\truez{600}%
\put(-\value{x},\value{x}){\makebox(0,\value{z})[r]{${#1}$}}%
\put(\value{x},-\value{y}){\makebox(0,\value{z})[l]{${#2}$}}%
\end{picture}}%

\newcommand{\WSWMONO}[3]{%
\Y=#3%
\divide\Y by 2%
\Z=\Y%
\divide\Z by 2%
\TrueTail%
\bimolength=#3%
\advance\bimolength by -\TrueMonoTail%
\monolength=\bimolength%
\advance\monolength by -\Y%
\secondmonolength=\monolength%
\divide\secondmonolength by 2%
\begin{picture}(0,0)%
\put(\monolength,\secondmonolength){\line(-2,-1){\bimolength}}%
\put(\monolength,\secondmonolength){\wswhead}%
\put(-\Y,-\Z){\wswhead}%
\truex{200}\truey{800}\truez{600}%
\put(-\value{x},\value{x}){\makebox(0,\value{z})[r]{${#1}$}}%
\put(\value{x},-\value{y}){\makebox(0,\value{z})[l]{${#2}$}}%
\end{picture}}%

\newcommand{\WSWEPI}[3]{%
\Y=#3%
\divide\Y by 2%
\Z=\Y%
\divide\Z by 2%
\TrueHead%
\bimolength=#3%
\advance\bimolength by -\TrueEpiHead%
\epilength=\bimolength%
\advance\epilength by -\Y%
\secondepilength=\epilength%
\divide\secondepilength by 2%
\begin{picture}(0,0)%
\put(\Y,\Z){\line(-2,-1){\bimolength}}%
\put(-\epilength,-\secondepilength){\wswhead}%
\put(-\Y,-\Z){\wswhead}%
\truex{200}\truey{800}\truez{600}%
\put(-\value{x},\value{x}){\makebox(0,\value{z})[r]{${#1}$}}%
\put(\value{x},-\value{y}){\makebox(0,\value{z})[l]{${#2}$}}%
\end{picture}}%

\newcommand{\WSWBIMO}[3]{%
\Y=#3%
\divide\Y by 2%
\Z=\Y%
\divide\Z by 2%
\TrueTail\TrueHead%
\bimolength=#3%
\advance\bimolength by -\TrueMonoTail%
\monolength=\bimolength%
\advance\monolength by -\Y%
\advance\bimolength by -\TrueEpiHead%
\epilength=\bimolength%
\advance\epilength by -\monolength%
\secondmonolength=\monolength%
\divide\secondmonolength by 2%
\secondepilength=\epilength%
\divide\secondepilength by 2%
\begin{picture}(0,0)%
\put(\monolength,\secondmonolength){\line(-2,-1){\bimolength}}%
\put(\monolength,\secondmonolength){\wswhead}%
\put(-\epilength,-\secondepilength){\wswhead}%
\put(-\Y,-\Z){\wswhead}%
\truex{200}\truey{800}\truez{600}%
\put(-\value{x},\value{x}){\makebox(0,\value{z})[r]{${#1}$}}%
\put(\value{x},-\value{y}){\makebox(0,\value{z})[l]{${#2}$}}%
\end{picture}}%

\newcommand{\WSWBIAR}[3]{%
\Y=#3%
\divide\Y by 2%
\Z=\Y%
\divide\Z by 2%
\begin{picture}(0,0)%
\put(\Y,\Z){\begin{picture}(0,0)%
\truex{156}\truey{313}%
\put(-\value{x},\value{y}){\line(-2,-1){#3}}%
\put(\value{x},-\value{y}){\line(-2,-1){#3}}%
\monolength=#3%
\advance\monolength by -\value{x}%
\epilength=#3%
\advance\epilength by \value{x}%
\secondmonolength=\Y%
\advance\secondmonolength by -\value{y}%
\secondepilength=\Y%
\advance\secondepilength by \value{y}%
\put(-\monolength,-\secondepilength){\wswhead}%
\put(-\epilength,-\secondmonolength){\wswhead}%
\end{picture}}
\truex{300}\truey{1000}\truez{600}%
\put(-\value{x},\value{x}){\makebox(0,\value{z})[r]{${#1}$}}%
\put(\value{x},-\value{y}){\makebox(0,\value{z})[l]{${#2}$}}%
\end{picture}}%

\newcommand{\WSWBIDIST}[3]{%
\Y=#3%
\divide\Y by 2%
\Z=\Y%
\divide\Z by 2%
\begin{picture}(0,0)%
\truex{156}\truey{313}\truez{400}%
\put(\Y,\Z){\begin{picture}(0,0)%
\put(-\value{x},\value{y}){\line(-2,-1){#3}}%
\put(\value{x},-\value{y}){\line(-2,-1){#3}}%
\monolength=#3%
\advance\monolength by -\value{x}%
\epilength=#3%
\advance\epilength by \value{x}%
\secondmonolength=\Y%
\advance\secondmonolength by -\value{y}%
\secondepilength=\Y%
\advance\secondepilength by \value{y}%
\put(-\monolength,-\secondepilength){\wswhead}%
\put(-\epilength,-\secondmonolength){\wswhead}%
\end{picture}}
\put(-\value{x},\value{y}){\distsign}%
\put(\value{x},-\value{y}){\distsign}%
\truex{300}\truey{1000}\truez{600}%
\put(-\value{x},\value{x}){\makebox(0,\value{z})[r]{${#1}$}}%
\put(\value{x},-\value{y}){\makebox(0,\value{z})[l]{${#2}$}}%
\end{picture}}%

\newcommand{\WSWADJAR}[3]{%
\Y=#3%
\divide\Y by 2%
\Z=\Y%
\divide\Z by 2%
\begin{picture}(0,0)%
\put(\Y,\Z){\begin{picture}(0,0)%
\truex{156}\truey{313}%
\monolength=#3%
\advance\monolength by -\value{x}%
\epilength=#3%
\advance\epilength by \value{x}%
\secondmonolength=\Y%
\advance\secondmonolength by -\value{y}%
\secondepilength=\Y%
\advance\secondepilength by \value{y}%
\put(\value{x},-\value{y}){\line(-2,-1){#3}}%
\put(-\monolength,-\secondepilength){\wswhead}%
\put(-\epilength,-\secondmonolength){\line(2,1){#3}}%
\put(-\value{x},\value{y}){\enehead}%
\end{picture}}
\truex{300}\truey{1000}\truez{600}%
\put(-\value{x},\value{x}){\makebox(0,\value{z})[r]{${#1}$}}%
\put(\value{x},-\value{y}){\makebox(0,\value{z})[l]{${#2}$}}%
\end{picture}}%

\newcommand{\WSWADJDIST}[3]{%
\Y=#3%
\divide\Y by 2%
\Z=\Y%
\divide\Z by 2%
\begin{picture}(0,0)%
\truex{156}\truey{313}\truez{400}%
\put(\Y,\Z){\begin{picture}(0,0)%
\monolength=#3%
\advance\monolength by -\value{x}%
\epilength=#3%
\advance\epilength by \value{x}%
\secondmonolength=\Y%
\advance\secondmonolength by -\value{y}%
\secondepilength=\Y%
\advance\secondepilength by \value{y}%
\put(\value{x},-\value{y}){\line(-2,-1){#3}}%
\put(-\monolength,-\secondepilength){\wswhead}%
\put(-\epilength,-\secondmonolength){\line(2,1){#3}}%
\put(-\value{x},\value{y}){\enehead}%
\end{picture}}
\put(-\value{x},\value{y}){\distsign}%
\put(\value{x},-\value{y}){\distsign}%
\truex{300}\truey{1000}\truez{600}%
\put(-\value{x},\value{x}){\makebox(0,\value{z})[r]{${#1}$}}%
\put(\value{x},-\value{y}){\makebox(0,\value{z})[l]{${#2}$}}%
\end{picture}}%

\def\basicwswar[#1]{\WSWAR{}{}{#100}}%

\newcommand{\wswar}{\@ifnextchar[{\basicwswar}{\basicwswar[133]}}%

\def\basicWswar[#1]#2{\WSWAR{#2}{}{#100}}%

\newcommand{\Wswar}{\@ifnextchar[{\basicWswar}{\basicWswar[133]}}%

\def\basicwswaR[#1]#2{\WSWAR{}{#2}{#100}}%

\newcommand{\wswaR}{\@ifnextchar[{\basicwswaR}{\basicwswaR[133]}}%

\def\basicwswdist[#1]{\WSWDIST{}{}{#100}}%

\newcommand{\wswdist}{\@ifnextchar[{\basicwswdist}{\basicwswdist[133]}}%

\def\basicWswdist[#1]#2{\WSWDIST{#2}{}{#100}}%

\newcommand{\Wswdist}{\@ifnextchar[{\basicWswdist}{\basicWswdist[133]}}%

\def\basicwswdisT[#1]#2{\WSWDIST{}{#2}{#100}}%

\newcommand{\wswdisT}{\@ifnextchar[{\basicwswdisT}{\basicwswdisT[133]}}%

\def\basicwswdotar[#1]{\WSWDOTAR{}{}{#100}}%

\newcommand{\wswdotar}{\@ifnextchar[{\basicwswdotar}{\basicwswdotar[133]}}%

\def\basicWswdotar[#1]#2{\WSWDOTAR{#2}{}{#100}}%

\newcommand{\Wswdotar}{\@ifnextchar[{\basicWswdotar}{\basicWswdotar[133]}}%

\def\basicwswdotaR[#1]#2{\WSWDOTAR{}{#2}{#100}}%

\newcommand{\wswdotaR}{\@ifnextchar[{\basicwswdotaR}{\basicwswdotaR[133]}}%

\def\basicwswmono[#1]{\WSWMONO{}{}{#100}}%

\newcommand{\wswmono}{\@ifnextchar[{\basicwswmono}{\basicwswmono[133]}}%

\def\basicWswmono[#1]#2{\WSWMONO{#2}{}{#100}}%

\newcommand{\Wswmono}{\@ifnextchar[{\basicWswmono}{\basicWswmono[133]}}%

\def\basicwswmonO[#1]#2{\WSWMONO{}{#2}{#100}}%

\newcommand{\wswmonO}{\@ifnextchar[{\basicwswmonO}{\basicwswmonO[133]}}%

\def\basicwswepi[#1]{\WSWEPI{}{}{#100}}%

\newcommand{\wswepi}{\@ifnextchar[{\basicwswepi}{\basicwswepi[133]}}%

\def\basicWswepi[#1]#2{\WSWEPI{#2}{}{#100}}%

\newcommand{\Wswepi}{\@ifnextchar[{\basicWswepi}{\basicWswepi[133]}}%

\def\basicwswepI[#1]#2{\WSWEPI{}{#2}{#100}}%

\newcommand{\wswepI}{\@ifnextchar[{\basicwswepI}{\basicwswepI[133]}}%

\def\basicwswbimo[#1]{\WSWBIMO{}{}{#100}}%

\newcommand{\wswbimo}{\@ifnextchar[{\basicwswbimo}{\basicwswbimo[133]}}%

\def\basicWswbimo[#1]#2{\WSWBIMO{#2}{}{#100}}%

\newcommand{\Wswbimo}{\@ifnextchar[{\basicWswbimo}{\basicWswbimo[133]}}%

\def\basicwswbimO[#1]#2{\WSWBIMO{}{#2}{#100}}%

\newcommand{\wswbimO}{\@ifnextchar[{\basicwswbimO}{\basicwswbimO[133]}}%

\def\basicwswiso[#1]{\WSWAR{\isosign}{}{#100}}%

\newcommand{\wswiso}{\@ifnextchar[{\basicwswiso}{\basicwswiso[133]}}%

\def\basicWswiso[#1]#2{\WSWAR{#2}{\isosign}{#100}}%

\newcommand{\Wswiso}{\@ifnextchar[{\basicWswiso}{\basicWswiso[133]}}%

\def\basicwswisO[#1]#2{\WSWAR{\isosign}{#2}{#100}}%

\newcommand{\wswisO}{\@ifnextchar[{\basicwswisO}{\basicwswisO[133]}}%

\def\basicwswbiar[#1]{\WSWBIAR{}{}{#100}}%

\newcommand{\wswbiar}{\@ifnextchar[{\basicwswbiar}{\basicwswbiar[133]}}%

\def\basicWswbiar[#1]#2#3{\WSWBIAR{#2}{#3}{#100}}%

\newcommand{\Wswbiar}{\@ifnextchar[{\basicWswbiar}{\basicWswbiar[133]}}%

\def\basicwswbidist[#1]{\WSWBIDIST{}{}{#100}}%

\newcommand{\wswbidist}{\@ifnextchar[{\basicwswbidist}{\basicwswbidist[133]}}%

\def\basicWswbidist[#1]#2#3{\WSWBIDIST{#2}{#3}{#100}}%

\newcommand{\Wswbidist}{\@ifnextchar[{\basicWswbidist}{\basicWswbidist[133]}}%

\def\basicwswadjar[#1]{\WSWADJAR{}{}{#100}}%

\newcommand{\wswadjar}{\@ifnextchar[{\basicwswadjar}{\basicwswadjar[133]}}%

\def\basicWswadjar[#1]#2#3{\WSWADJAR{#2}{#3}{#100}}%

\newcommand{\Wswadjar}{\@ifnextchar[{\basicWswadjar}{\basicWswadjar[133]}}%

\def\basicwswadjdist[#1]{\WSWADJDIST{}{}{#100}}%

\newcommand{\wswadjdist}{\@ifnextchar[{\basicwswadjdist}{\basicwswadjdist[13
3]}}%

\def\basicWswadjdist[#1]#2#3{\WSWADJDIST{#2}{#3}{#100}}%

\newcommand{\Wswadjdist}{\@ifnextchar[{\basicWswadjdist}{\basicWswadjdist[13
3]}}%

\newcommand{\WNWAR}[3]{%
\Y=#3%
\divide\Y by 2%
\Z=\Y%
\divide\Z by 2%
\begin{picture}(0,0)%
\put(\Y,-\Z){\line(-2,1){#3}}%
\put(-\Y,\Z){\wnwhead}%
\truex{200}\truey{800}\truez{600}%
\put(\value{x},\value{x}){\makebox(0,\value{z})[l]{${#1}$}}%
\put(-\value{x},-\value{y}){\makebox(0,\value{z})[r]{${#2}$}}%
\end{picture}}%

\newcommand{\WNWDIST}[3]{%
\Y=#3%
\divide\Y by 2%
\Z=\Y%
\divide\Z by 2%
\begin{picture}(0,0)%
\put(\Y,-\Z){\line(-2,1){#3}}%
\put(-\Y,\Z){\wnwhead}%
\truex{400}%
\put(0,0){\distsign}%
\truex{200}\truey{800}\truez{600}%
\put(\value{x},\value{x}){\makebox(0,\value{z})[l]{${#1}$}}%
\put(-\value{x},-\value{y}){\makebox(0,\value{z})[r]{${#2}$}}%
\end{picture}}%

\newcommand{\WNWDOTAR}[3]{%
\truex{100}\truey{268}\truez{134}%
\Y=#3%
\divide\Y by 2%
\Z=\Y%
\divide\Z by 2%
\NUMBEROFDOTS=#3%
\divide\NUMBEROFDOTS by \value{y}%
\advance\NUMBEROFDOTS by 1%
\begin{picture}(0,0)%
\multiput(\Y,-\Z)(-\value{y},\value{z}){\NUMBEROFDOTS}%
{\circle*{\value{x}}}%
\put(-\Y,\Z){\wnwhead}%
\truex{200}\truey{800}\truez{600}%
\put(\value{x},\value{x}){\makebox(0,\value{z})[l]{${#1}$}}%
\put(-\value{x},-\value{y}){\makebox(0,\value{z})[r]{${#2}$}}%
\end{picture}}%

\newcommand{\WNWMONO}[3]{%
\Y=#3%
\divide\Y by 2%
\Z=\Y%
\divide\Z by 2%
\TrueTail%
\bimolength=#3%
\advance\bimolength by -\TrueMonoTail%
\monolength=\bimolength%
\advance\monolength by -\Y%
\secondmonolength=\monolength%
\divide\secondmonolength by 2%
\begin{picture}(0,0)%
\put(\monolength,-\secondmonolength){\line(-2,1){\bimolength}}%
\put(\monolength,-\secondmonolength){\wnwhead}%
\put(-\Y,\Z){\wnwhead}%
\truex{200}\truey{800}\truez{600}%
\put(\value{x},\value{x}){\makebox(0,\value{z})[l]{${#1}$}}%
\put(-\value{x},-\value{y}){\makebox(0,\value{z})[r]{${#2}$}}%
\end{picture}}%

\newcommand{\WNWEPI}[3]{%
\Y=#3%
\divide\Y by 2%
\Z=\Y%
\divide\Z by 2%
\TrueHead%
\bimolength=#3%
\advance\bimolength by -\TrueEpiHead%
\epilength=\bimolength%
\advance\epilength by -\Y%
\secondepilength=\epilength%
\divide\secondepilength by 2%
\begin{picture}(0,0)%
\put(\Y,-\Z){\line(-2,1){\bimolength}}%
\put(-\epilength,\secondepilength){\wnwhead}%
\put(-\Y,\Z){\wnwhead}%
\truex{200}\truey{800}\truez{600}%
\put(\value{x},\value{x}){\makebox(0,\value{z})[l]{${#1}$}}%
\put(-\value{x},-\value{y}){\makebox(0,\value{z})[r]{${#2}$}}%
\end{picture}}%

\newcommand{\WNWBIMO}[3]{%
\Y=#3%
\divide\Y by 2%
\Z=\Y%
\divide\Z by 2%
\TrueTail\TrueHead%
\bimolength=#3%
\advance\bimolength by -\TrueMonoTail%
\monolength=\bimolength%
\advance\monolength by -\Y%
\advance\bimolength by -\TrueEpiHead%
\epilength=\bimolength%
\advance\epilength by -\monolength%
\secondmonolength=\monolength%
\divide\secondmonolength by 2%
\secondepilength=\epilength%
\divide\secondepilength by 2%
\begin{picture}(0,0)%
\put(\monolength,-\secondmonolength){\line(-2,1){\bimolength}}%
\put(\monolength,-\secondmonolength){\wnwhead}%
\put(-\epilength,\secondepilength){\wnwhead}%
\put(-\Y,\Z){\wnwhead}%
\truex{200}\truey{800}\truez{600}%
\put(\value{x},\value{x}){\makebox(0,\value{z})[l]{${#1}$}}%
\put(-\value{x},-\value{y}){\makebox(0,\value{z})[r]{${#2}$}}%
\end{picture}}%

\newcommand{\WNWBIAR}[3]{%
\Y=#3%
\divide\Y by 2%
\Z=\Y%
\divide\Z by 2%
\begin{picture}(0,0)%
\put(\Y,-\Z){\begin{picture}(0,0)%
\truex{156}\truey{313}%
\put(-\value{x},-\value{y}){\line(-2,1){#3}}%
\put(\value{x},\value{y}){\line(-2,1){#3}}%
\monolength=#3%
\advance\monolength by -\value{x}%
\epilength=#3%
\advance\epilength by \value{x}%
\secondmonolength=\Y%
\advance\secondmonolength by -\value{y}%
\secondepilength=\Y%
\advance\secondepilength by \value{y}%
\put(-\monolength,\secondepilength){\wnwhead}%
\put(-\epilength,\secondmonolength){\wnwhead}%
\end{picture}}
\truex{400}\truey{1000}\truez{600}%
\put(\value{x},\value{x}){\makebox(0,\value{z})[l]{${#1}$}}%
\put(-\value{x},-\value{y}){\makebox(0,\value{z})[r]{${#2}$}}%
\end{picture}}%

\newcommand{\WNWBIDIST}[3]{%
\Y=#3%
\divide\Y by 2%
\Z=\Y%
\divide\Z by 2%
\begin{picture}(0,0)%
\truex{156}\truey{313}\truez{400}%
\put(\Y,-\Z){\begin{picture}(0,0)%
\put(-\value{x},-\value{y}){\line(-2,1){#3}}%
\put(\value{x},\value{y}){\line(-2,1){#3}}%
\monolength=#3%
\advance\monolength by -\value{x}%
\epilength=#3%
\advance\epilength by \value{x}%
\secondmonolength=\Y%
\advance\secondmonolength by -\value{y}%
\secondepilength=\Y%
\advance\secondepilength by \value{y}%
\put(-\monolength,\secondepilength){\wnwhead}%
\put(-\epilength,\secondmonolength){\wnwhead}%
\end{picture}}
\put(\value{x},\value{y}){\distsign}%
\put(-\value{x},-\value{y}){\distsign}%
\truex{400}\truey{1000}\truez{600}%
\put(\value{x},\value{x}){\makebox(0,\value{z})[l]{${#1}$}}%
\put(-\value{x},-\value{y}){\makebox(0,\value{z})[r]{${#2}$}}%
\end{picture}}%

\newcommand{\WNWADJAR}[3]{%
\Y=#3%
\divide\Y by 2%
\Z=\Y%
\divide\Z by 2%
\begin{picture}(0,0)%
\put(\Y,-\Z){\begin{picture}(0,0)%
\truex{156}\truey{313}%
\monolength=#3%
\advance\monolength by -\value{x}%
\epilength=#3%
\advance\epilength by \value{x}%
\secondmonolength=\Y%
\advance\secondmonolength by -\value{y}%
\secondepilength=\Y%
\advance\secondepilength by \value{y}%
\put(-\value{x},-\value{y}){\line(-2,1){#3}}%
\put(-\epilength,\secondmonolength){\wnwhead}%
\put(-\monolength,\secondepilength){\line(2,-1){#3}}%
\put(\value{x},\value{y}){\esehead}%
\end{picture}}
\truex{400}\truey{1000}\truez{600}%
\put(\value{x},\value{x}){\makebox(0,\value{z})[l]{${#1}$}}%
\put(-\value{x},-\value{y}){\makebox(0,\value{z})[r]{${#2}$}}%
\end{picture}}%

\newcommand{\WNWADJDIST}[3]{%
\Y=#3%
\divide\Y by 2%
\Z=\Y%
\divide\Z by 2%
\begin{picture}(0,0)%
\truex{156}\truey{313}\truez{400}%
\put(\Y,-\Z){\begin{picture}(0,0)%
\monolength=#3%
\advance\monolength by -\value{x}%
\epilength=#3%
\advance\epilength by \value{x}%
\secondmonolength=\Y%
\advance\secondmonolength by -\value{y}%
\secondepilength=\Y%
\advance\secondepilength by \value{y}%
\put(-\value{x},-\value{y}){\line(-2,1){#3}}%
\put(-\epilength,\secondmonolength){\wnwhead}%
\put(-\monolength,\secondepilength){\line(2,-1){#3}}%
\put(\value{x},\value{y}){\esehead}%
\end{picture}}
\put(\value{x},\value{y}){\distsign}%
\put(-\value{x},-\value{y}){\distsign}%
\truex{400}\truey{1000}\truez{600}%
\put(\value{x},\value{x}){\makebox(0,\value{z})[l]{${#1}$}}%
\put(-\value{x},-\value{y}){\makebox(0,\value{z})[r]{${#2}$}}%
\end{picture}}%

\def\basicwnwar[#1]{\WNWAR{}{}{#100}}%

\newcommand{\wnwar}{\@ifnextchar[{\basicwnwar}{\basicwnwar[133]}}%

\def\basicWnwar[#1]#2{\WNWAR{#2}{}{#100}}%

\newcommand{\Wnwar}{\@ifnextchar[{\basicWnwar}{\basicWnwar[133]}}%

\def\basicwnwaR[#1]#2{\WNWAR{}{#2}{#100}}%

\newcommand{\wnwaR}{\@ifnextchar[{\basicwnwaR}{\basicwnwaR[133]}}%

\def\basicwnwdist[#1]{\WNWDIST{}{}{#100}}%

\newcommand{\wnwdist}{\@ifnextchar[{\basicwnwdist}{\basicwnwdist[133]}}%

\def\basicWnwdist[#1]#2{\WNWDIST{#2}{}{#100}}%

\newcommand{\Wnwdist}{\@ifnextchar[{\basicWnwdist}{\basicWnwdist[133]}}%

\def\basicwnwdisT[#1]#2{\WNWDIST{}{#2}{#100}}%

\newcommand{\wnwdisT}{\@ifnextchar[{\basicwnwdisT}{\basicwnwdisT[133]}}%

\def\basicwnwdotar[#1]{\WNWDOTAR{}{}{#100}}%

\newcommand{\wnwdotar}{\@ifnextchar[{\basicwnwdotar}{\basicwnwdotar[133]}}%

\def\basicWnwdotar[#1]#2{\WNWDOTAR{#2}{}{#100}}%

\newcommand{\Wnwdotar}{\@ifnextchar[{\basicWnwdotar}{\basicWnwdotar[133]}}%

\def\basicwnwdotaR[#1]#2{\WNWDOTAR{}{#2}{#100}}%

\newcommand{\wnwdotaR}{\@ifnextchar[{\basicwnwdotaR}{\basicwnwdotaR[133]}}%

\def\basicwnwmono[#1]{\WNWMONO{}{}{#100}}%

\newcommand{\wnwmono}{\@ifnextchar[{\basicwnwmono}{\basicwnwmono[133]}}%

\def\basicWnwmono[#1]#2{\WNWMONO{#2}{}{#100}}%

\newcommand{\Wnwmono}{\@ifnextchar[{\basicWnwmono}{\basicWnwmono[133]}}%

\def\basicwnwmonO[#1]#2{\WNWMONO{}{#2}{#100}}%

\newcommand{\wnwmonO}{\@ifnextchar[{\basicwnwmonO}{\basicwnwmonO[133]}}%

\def\basicwnwepi[#1]{\WNWEPI{}{}{#100}}%

\newcommand{\wnwepi}{\@ifnextchar[{\basicwnwepi}{\basicwnwepi[133]}}%

\def\basicWnwepi[#1]#2{\WNWEPI{#2}{}{#100}}%

\newcommand{\Wnwepi}{\@ifnextchar[{\basicWnwepi}{\basicWnwepi[133]}}%

\def\basicwnwepI[#1]#2{\WNWEPI{}{#2}{#100}}%

\newcommand{\wnwepI}{\@ifnextchar[{\basicwnwepI}{\basicwnwepI[133]}}%

\def\basicwnwbimo[#1]{\WNWBIMO{}{}{#100}}%

\newcommand{\wnwbimo}{\@ifnextchar[{\basicwnwbimo}{\basicwnwbimo[133]}}%

\def\basicWnwbimo[#1]#2{\WNWBIMO{#2}{}{#100}}%

\newcommand{\Wnwbimo}{\@ifnextchar[{\basicWnwbimo}{\basicWnwbimo[133]}}%

\def\basicwnwbimO[#1]#2{\WNWBIMO{}{#2}{#100}}%

\newcommand{\wnwbimO}{\@ifnextchar[{\basicwnwbimO}{\basicwnwbimO[133]}}%

\def\basicwnwiso[#1]{\WNWAR{\isosign}{}{#100}}%

\newcommand{\wnwiso}{\@ifnextchar[{\basicwnwiso}{\basicwnwiso[133]}}%

\def\basicWnwiso[#1]#2{\WNWAR{#2}{\isosign}{#100}}%

\newcommand{\Wnwiso}{\@ifnextchar[{\basicWnwiso}{\basicWnwiso[133]}}%

\def\basicwnwisO[#1]#2{\WNWAR{\isosign}{#2}{#100}}%

\newcommand{\wnwisO}{\@ifnextchar[{\basicwnwisO}{\basicwnwisO[133]}}%

\def\basicwnwbiar[#1]{\WNWBIAR{}{}{#100}}%

\newcommand{\wnwbiar}{\@ifnextchar[{\basicwnwbiar}{\basicwnwbiar[133]}}%

\def\basicWnwbiar[#1]#2#3{\WNWBIAR{#2}{#3}{#100}}%

\newcommand{\Wnwbiar}{\@ifnextchar[{\basicWnwbiar}{\basicWnwbiar[133]}}%

\def\basicwnwbidist[#1]{\WNWBIDIST{}{}{#100}}%

\newcommand{\wnwbidist}{\@ifnextchar[{\basicwnwbidist}{\basicwnwbidist[133]}}%

\def\basicWnwbidist[#1]#2#3{\WNWBIDIST{#2}{#3}{#100}}%

\newcommand{\Wnwbidist}{\@ifnextchar[{\basicWnwbidist}{\basicWnwbidist[133]}}%

\def\basicwnwadjar[#1]{\WNWADJAR{}{}{#100}}%

\newcommand{\wnwadjar}{\@ifnextchar[{\basicwnwadjar}{\basicwnwadjar[133]}}%

\def\basicWnwadjar[#1]#2#3{\WNWADJAR{#2}{#3}{#100}}%

\newcommand{\Wnwadjar}{\@ifnextchar[{\basicWnwadjar}{\basicWnwadjar[133]}}%

\def\basicwnwadjdist[#1]{\WNWADJDIST{}{}{#100}}%

\newcommand{\wnwadjdist}{\@ifnextchar[{\basicwnwadjdist}{\basicwnwadjdist[13
3]}}%

\def\basicWnwadjdist[#1]#2#3{\WNWADJDIST{#2}{#3}{#100}}%

\newcommand{\Wnwadjdist}{\@ifnextchar[{\basicWnwadjdist}{\basicWnwadjdist[13
3]}}%

\newcommand{\NNEAR}[3]{%
\Z=#3%
\divide\Z by 2%
\begin{picture}(0,0)%
\put(-\Z,-#3){\line(1,2){#3}}%
\put(\Z,#3){\nnehead}%
\truex{200}\truey{800}\truez{600}%
\put(-\value{x},\value{x}){\makebox(0,\value{z})[r]{${#1}$}}%
\put(\value{x},-\value{y}){\makebox(0,\value{z})[l]{${#2}$}}%
\end{picture}}%

\newcommand{\NNEDIST}[3]{%
\Z=#3%
\divide\Z by 2%
\begin{picture}(0,0)%
\put(-\Z,-#3){\line(1,2){#3}}%
\put(\Z,#3){\nnehead}%
\truex{400}%
\put(0,0){\distsign}%
\truex{200}\truey{800}\truez{600}%
\put(-\value{x},\value{x}){\makebox(0,\value{z})[r]{${#1}$}}%
\put(\value{x},-\value{y}){\makebox(0,\value{z})[l]{${#2}$}}%
\end{picture}}%

\newcommand{\NNEDOTAR}[3]{%
\truex{100}\truey{268}\truez{134}%
\Z=#3%
\divide\Z by 2%
\NUMBEROFDOTS=#3%
\divide\NUMBEROFDOTS by \value{z}%
\advance\NUMBEROFDOTS by 1%
\begin{picture}(0,0)%
\multiput(-\Z,-#3)(\value{z},\value{y}){\NUMBEROFDOTS}%
{\circle*{\value{x}}}%
\put(\Z,#3){\nnehead}%
\truex{200}\truey{800}\truez{600}%
\put(-\value{x},\value{x}){\makebox(0,\value{z})[r]{${#1}$}}%
\put(\value{x},-\value{y}){\makebox(0,\value{z})[l]{${#2}$}}%
\end{picture}}%

\newcommand{\NNEMONO}[3]{%
\Z=#3%
\divide\Z by 2%
\truetaiL%
\bimolength=#3%
\advance\bimolength by -\truemonotaiL%
\monolength=\bimolength%
\advance\monolength by -\Z%
\secondmonolength=\monolength%
\multiply\secondmonolength by 2%
\begin{picture}(0,0)%
\put(-\monolength,-\secondmonolength){\line(1,2){\bimolength}}%
\put(-\monolength,-\secondmonolength){\nnehead}%
\put(\Z,#3){\nnehead}%
\truex{200}\truey{800}\truez{600}%
\put(-\value{x},\value{x}){\makebox(0,\value{z})[r]{${#1}$}}%
\put(\value{x},-\value{y}){\makebox(0,\value{z})[l]{${#2}$}}%
\end{picture}}%

\newcommand{\NNEEPI}[3]{%
\Z=#3%
\divide\Z by 2%
\trueheaD%
\bimolength=#3%
\advance\bimolength by -\trueepiheaD%
\epilength=\bimolength%
\advance\epilength by -\Z%
\secondepilength=\epilength%
\multiply\secondepilength by 2%
\begin{picture}(0,0)%
\put(-\Z,-#3){\line(1,2){\bimolength}}%
\put(\epilength,\secondepilength){\nnehead}%
\put(\Z,#3){\nnehead}%
\truex{200}\truey{800}\truez{600}%
\put(-\value{x},\value{x}){\makebox(0,\value{z})[r]{${#1}$}}%
\put(\value{x},-\value{y}){\makebox(0,\value{z})[l]{${#2}$}}%
\end{picture}}%

\newcommand{\NNEBIMO}[3]{%
\Z=#3%
\divide\Z by 2%
\truetaiL\trueheaD%
\bimolength=#3%
\advance\bimolength by -\truemonotaiL%
\monolength=\bimolength%
\advance\monolength by -\Z%
\advance\bimolength by -\trueepiheaD%
\epilength=\bimolength%
\advance\epilength by -\monolength%
\secondmonolength=\monolength%
\multiply\secondmonolength by 2%
\secondepilength=\epilength%
\multiply\secondepilength by 2%
\begin{picture}(0,0)%
\put(-\monolength,-\secondmonolength){\line(1,2){\bimolength}}%
\put(-\monolength,-\secondmonolength){\nnehead}%
\put(\epilength,\secondepilength){\nnehead}%
\put(\Z,#3){\nnehead}%
\truex{200}\truey{800}\truez{600}%
\put(-\value{x},\value{x}){\makebox(0,\value{z})[r]{${#1}$}}%
\put(\value{x},-\value{y}){\makebox(0,\value{z})[l]{${#2}$}}%
\end{picture}}%

\newcommand{\NNEEQL}[3]{%
\Z=#3%
\divide\Z by 2%
\begin{picture}(0,0)%
\put(-\Z,-#3){\begin{picture}(0,0)%
\truex{44}\truey{89}%
\put(-\value{y},\value{x}){\line(1,2){#3}}%
\put(\value{y},-\value{x}){\line(1,2){#3}}%
\end{picture}}%
\truex{200}\truey{800}\truez{600}%
\put(-\value{x},\value{x}){\makebox(0,\value{z})[r]{${#1}$}}%
\put(\value{x},-\value{y}){\makebox(0,\value{z})[l]{${#2}$}}%
\end{picture}}%

\newcommand{\NNEBIAR}[3]{%
\Y=#3%
\divide\Y by 2%
\Z=#3%
\multiply \Z by 2%
\begin{picture}(0,0)%
\put(-\Y,-#3){\begin{picture}(0,0)%
\truex{313}\truey{156}%
\put(-\value{x},\value{y}){\line(1,2){#3}}%
\put(\value{x},-\value{y}){\line(1,2){#3}}%
\monolength=#3%
\advance\monolength by -\value{x}%
\epilength=#3%
\advance\epilength by \value{x}%
\secondmonolength=\Z%
\advance\secondmonolength by -\value{y}%
\secondepilength=\Z%
\advance\secondepilength by \value{y}%
\put(\monolength,\secondepilength){\nnehead}%
\put(\epilength,\secondmonolength){\nnehead}%
\end{picture}}
\truex{300}\truey{1000}\truez{600}%
\put(-\value{x},\value{x}){\makebox(0,\value{z})[r]{${#1}$}}%
\put(\value{x},-\value{y}){\makebox(0,\value{z})[l]{${#2}$}}%
\end{picture}}%

\newcommand{\NNEBIDIST}[3]{%
\Y=#3%
\divide\Y by 2%
\Z=#3%
\multiply \Z by 2%
\begin{picture}(0,0)%
\truex{313}\truey{156}\truez{400}%
\put(-\Y,-#3){\begin{picture}(0,0)%
\put(-\value{x},\value{y}){\line(1,2){#3}}%
\put(\value{x},-\value{y}){\line(1,2){#3}}%
\monolength=#3%
\advance\monolength by -\value{x}%
\epilength=#3%
\advance\epilength by \value{x}%
\secondmonolength=\Z%
\advance\secondmonolength by -\value{y}%
\secondepilength=\Z%
\advance\secondepilength by \value{y}%
\put(\monolength,\secondepilength){\nnehead}%
\put(\epilength,\secondmonolength){\nnehead}%
\end{picture}}
\put(-\value{x},\value{y}){\distsign}%
\put(\value{x},-\value{y}){\distsign}%
\truex{300}\truey{1000}\truez{600}%
\put(-\value{x},\value{x}){\makebox(0,\value{z})[r]{${#1}$}}%
\put(\value{x},-\value{y}){\makebox(0,\value{z})[l]{${#2}$}}%
\end{picture}}%

\newcommand{\NNEADJAR}[3]{%
\Y=#3%
\divide\Y by 2%
\Z=#3%
\multiply \Z by 2%
\begin{picture}(0,0)%
\put(-\Y,-#3){\begin{picture}(0,0)%
\truex{313}\truey{156}%
\monolength=#3%
\advance\monolength by -\value{x}%
\epilength=#3%
\advance\epilength by \value{x}%
\secondmonolength=\Z%
\advance\secondmonolength by -\value{y}%
\secondepilength=\Z%
\advance\secondepilength by \value{y}%
\put(\value{x},-\value{y}){\line(1,2){#3}}%
\put(\epilength,\secondmonolength){\nnehead}%
\put(\monolength,\secondepilength){\line(-1,-2){#3}}%
\put(-\value{x},\value{y}){\sswhead}%
\end{picture}}
\truex{300}\truey{1000}\truez{600}%
\put(-\value{x},\value{x}){\makebox(0,\value{z})[r]{${#1}$}}%
\put(\value{x},-\value{y}){\makebox(0,\value{z})[l]{${#2}$}}%
\end{picture}}%

\newcommand{\NNEADJDIST}[3]{%
\Y=#3%
\divide\Y by 2%
\Z=#3%
\multiply \Z by 2%
\begin{picture}(0,0)%
\truex{313}\truey{156}\truez{400}%
\put(-\Y,-#3){\begin{picture}(0,0)%
\monolength=#3%
\advance\monolength by -\value{x}%
\epilength=#3%
\advance\epilength by \value{x}%
\secondmonolength=\Z%
\advance\secondmonolength by -\value{y}%
\secondepilength=\Z%
\advance\secondepilength by \value{y}%
\put(\value{x},-\value{y}){\line(1,2){#3}}%
\put(\epilength,\secondmonolength){\nnehead}%
\put(\monolength,\secondepilength){\line(-1,-2){#3}}%
\put(-\value{x},\value{y}){\sswhead}%
\end{picture}}
\put(-\value{x},\value{y}){\distsign}%
\put(\value{x},-\value{y}){\distsign}%
\truex{300}\truey{1000}\truez{600}%
\put(-\value{x},\value{x}){\makebox(0,\value{z})[r]{${#1}$}}%
\put(\value{x},-\value{y}){\makebox(0,\value{z})[l]{${#2}$}}%
\end{picture}}%

\def\basicnnear[#1]{\NNEAR{}{}{#100}}%

\newcommand{\nnear}{\@ifnextchar[{\basicnnear}{\basicnnear[67]}}%

\def\basicNnear[#1]#2{\NNEAR{#2}{}{#100}}%

\newcommand{\Nnear}{\@ifnextchar[{\basicNnear}{\basicNnear[67]}}%

\def\basicnneaR[#1]#2{\NNEAR{}{#2}{#100}}%

\newcommand{\nneaR}{\@ifnextchar[{\basicnneaR}{\basicnneaR[67]}}%

\def\basicnnedist[#1]{\NNEDIST{}{}{#100}}%

\newcommand{\nnedist}{\@ifnextchar[{\basicnnedist}{\basicnnedist[67]}}%

\def\basicNnedist[#1]#2{\NNEDIST{#2}{}{#100}}%

\newcommand{\Nnedist}{\@ifnextchar[{\basicNnedist}{\basicNnedist[67]}}%

\def\basicnnedisT[#1]#2{\NNEDIST{}{#2}{#100}}%

\newcommand{\nnedisT}{\@ifnextchar[{\basicnnedisT}{\basicnnedisT[67]}}%

\def\basicnnedotar[#1]{\NNEDOTAR{}{}{#100}}%

\newcommand{\nnedotar}{\@ifnextchar[{\basicnnedotar}{\basicnnedotar[67]}}%

\def\basicNnedotar[#1]#2{\NNEDOTAR{#2}{}{#100}}%

\newcommand{\Nnedotar}{\@ifnextchar[{\basicNnedotar}{\basicNnedotar[67]}}%

\def\basicnnedotaR[#1]#2{\NNEDOTAR{}{#2}{#100}}%

\newcommand{\nnedotaR}{\@ifnextchar[{\basicnnedotaR}{\basicnnedotaR[67]}}%

\def\basicnnemono[#1]{\NNEMONO{}{}{#100}}%

\newcommand{\nnemono}{\@ifnextchar[{\basicnnemono}{\basicnnemono[67]}}%

\def\basicNnemono[#1]#2{\NNEMONO{#2}{}{#100}}%

\newcommand{\Nnemono}{\@ifnextchar[{\basicNnemono}{\basicNnemono[67]}}%

\def\basicnnemonO[#1]#2{\NNEMONO{}{#2}{#100}}%

\newcommand{\nnemonO}{\@ifnextchar[{\basicnnemonO}{\basicnnemonO[67]}}%

\def\basicnneepi[#1]{\NNEEPI{}{}{#100}}%

\newcommand{\nneepi}{\@ifnextchar[{\basicnneepi}{\basicnneepi[67]}}%

\def\basicNneepi[#1]#2{\NNEEPI{#2}{}{#100}}%

\newcommand{\Nneepi}{\@ifnextchar[{\basicNneepi}{\basicNneepi[67]}}%

\def\basicnneepI[#1]#2{\NNEEPI{}{#2}{#100}}%

\newcommand{\nneepI}{\@ifnextchar[{\basicnneepI}{\basicnneepI[67]}}%

\def\basicnnebimo[#1]{\NNEBIMO{}{}{#100}}%

\newcommand{\nnebimo}{\@ifnextchar[{\basicnnebimo}{\basicnnebimo[67]}}%

\def\basicNnebimo[#1]#2{\NNEBIMO{#2}{}{#100}}%

\newcommand{\Nnebimo}{\@ifnextchar[{\basicNnebimo}{\basicNnebimo[67]}}%

\def\basicnnebimO[#1]#2{\NNEBIMO{}{#2}{#100}}%

\newcommand{\nnebimO}{\@ifnextchar[{\basicnnebimO}{\basicnnebimO[67]}}%

\def\basicnneiso[#1]{\NNEAR{\isosign}{}{#100}}%

\newcommand{\nneiso}{\@ifnextchar[{\basicnneiso}{\basicnneiso[67]}}%

\def\basicNneiso[#1]#2{\NNEAR{#2}{\isosign}{#100}}%

\newcommand{\Nneiso}{\@ifnextchar[{\basicNneiso}{\basicNneiso[67]}}%

\def\basicnneisO[#1]#2{\NNEAR{\isosign}{#2}{#100}}%

\newcommand{\nneisO}{\@ifnextchar[{\basicnneisO}{\basicnneisO[67]}}%

\def\basicnneeql[#1]{\NNEEQL{}{}{#100}}%

\newcommand{\nneeql}{\@ifnextchar[{\basicnneeql}{\basicnneeql[67]}}%

\def\basicNneeql[#1]#2{\NNEEQL{#2}{}{#100}}%

\newcommand{\Nneeql}{\@ifnextchar[{\basicNneeql}{\basicNneeql[67]}}%

\def\basicnneeqL[#1]#2{\NNEEQL{}{#2}{#100}}%

\newcommand{\nneeqL}{\@ifnextchar[{\basicnneeqL}{\basicnneeqL[67]}}%

\def\basicnnebiar[#1]{\NNEBIAR{}{}{#100}}%

\newcommand{\nnebiar}{\@ifnextchar[{\basicnnebiar}{\basicnnebiar[67]}}%

\def\basicNnebiar[#1]#2#3{\NNEBIAR{#2}{#3}{#100}}%

\newcommand{\Nnebiar}{\@ifnextchar[{\basicNnebiar}{\basicNnebiar[67]}}%

\def\basicnnebidist[#1]{\NNEBIDIST{}{}{#100}}%

\newcommand{\nnebidist}{\@ifnextchar[{\basicnnebidist}{\basicnnebidist[67]}}%

\def\basicNnebidist[#1]#2#3{\NNEBIDIST{#2}{#3}{#100}}%

\newcommand{\Nnebidist}{\@ifnextchar[{\basicNnebidist}{\basicNnebidist[67]}}%

\def\basicnneadjar[#1]{\NNEADJAR{}{}{#100}}%

\newcommand{\nneadjar}{\@ifnextchar[{\basicnneadjar}{\basicnneadjar[67]}}%

\def\basicNneadjar[#1]#2#3{\NNEADJAR{#2}{#3}{#100}}%

\newcommand{\Nneadjar}{\@ifnextchar[{\basicNneadjar}{\basicNneadjar[67]}}%

\def\basicnneadjdist[#1]{\NNEADJDIST{}{}{#100}}%

\newcommand{\nneadjdist}{\@ifnextchar[{\basicnneadjdist}{\basicnneadjdist[67]}}%

\def\basicNneadjdist[#1]#2#3{\NNEADJDIST{#2}{#3}{#100}}%

\newcommand{\Nneadjdist}{\@ifnextchar[{\basicNneadjdist}{\basicNneadjdist[67]}}%

\newcommand{\SSEAR}[3]{%
\Z=#3%
\divide\Z by 2%
\begin{picture}(0,0)%
\put(-\Z,#3){\line(1,-2){#3}}%
\put(\Z,-#3){\ssehead}%
\truex{200}\truey{800}\truez{600}%
\put(\value{x},\value{x}){\makebox(0,\value{z})[l]{${#1}$}}%
\put(-\value{x},-\value{y}){\makebox(0,\value{z})[r]{${#2}$}}%
\end{picture}}%

\newcommand{\SSEDIST}[3]{%
\Z=#3%
\divide\Z by 2%
\begin{picture}(0,0)%
\put(-\Z,#3){\line(1,-2){#3}}%
\put(\Z,-#3){\ssehead}%
\truex{400}%
\put(0,0){\distsign}%
\truex{200}\truey{800}\truez{600}%
\put(\value{x},\value{x}){\makebox(0,\value{z})[l]{${#1}$}}%
\put(-\value{x},-\value{y}){\makebox(0,\value{z})[r]{${#2}$}}%
\end{picture}}%

\newcommand{\SSEDOTAR}[3]{%
\truex{100}\truey{268}\truez{134}%
\Z=#3%
\divide\Z by 2%
\NUMBEROFDOTS=#3%
\divide\NUMBEROFDOTS by \value{z}%
\advance\NUMBEROFDOTS by 1%
\begin{picture}(0,0)%
\multiput(-\Z,#3)(\value{z},-\value{y}){\NUMBEROFDOTS}%
{\circle*{\value{x}}}%
\put(\Z,-#3){\ssehead}%
\truex{200}\truey{800}\truez{600}%
\put(\value{x},\value{x}){\makebox(0,\value{z})[l]{${#1}$}}%
\put(-\value{x},-\value{y}){\makebox(0,\value{z})[r]{${#2}$}}%
\end{picture}}%

\newcommand{\SSEMONO}[3]{%
\Z=#3%
\divide\Z by 2%
\truetaiL%
\bimolength=#3%
\advance\bimolength by -\truemonotaiL%
\monolength=\bimolength%
\advance\monolength by -\Z%
\secondmonolength=\monolength%
\multiply\secondmonolength by 2%
\begin{picture}(0,0)%
\put(-\monolength,\secondmonolength){\line(1,-2){\bimolength}}%
\put(-\monolength,\secondmonolength){\ssehead}%
\put(\Z,-#3){\ssehead}%
\truex{200}\truey{800}\truez{600}%
\put(\value{x},\value{x}){\makebox(0,\value{z})[l]{${#1}$}}%
\put(-\value{x},-\value{y}){\makebox(0,\value{z})[r]{${#2}$}}%
\end{picture}}%

\newcommand{\SSEEPI}[3]{%
\Z=#3%
\divide\Z by 2%
\trueheaD%
\bimolength=#3%
\advance\bimolength by -\trueepiheaD%
\epilength=\bimolength%
\advance\epilength by -\Z%
\secondepilength=\epilength%
\multiply\secondepilength by 2%
\begin{picture}(0,0)%
\put(-\Z,#3){\line(1,-2){\bimolength}}%
\put(\epilength,-\secondepilength){\ssehead}%
\put(\Z,-#3){\ssehead}%
\truex{200}\truey{800}\truez{600}%
\put(\value{x},\value{x}){\makebox(0,\value{z})[l]{${#1}$}}%
\put(-\value{x},-\value{y}){\makebox(0,\value{z})[r]{${#2}$}}%
\end{picture}}%

\newcommand{\SSEBIMO}[3]{%
\Z=#3%
\divide\Z by 2%
\truetaiL\trueheaD%
\bimolength=#3%
\advance\bimolength by -\truemonotaiL%
\monolength=\bimolength%
\advance\monolength by -\Z%
\advance\bimolength by -\trueepiheaD%
\epilength=\bimolength%
\advance\epilength by -\monolength%
\secondmonolength=\monolength%
\multiply\secondmonolength by 2%
\secondepilength=\epilength%
\multiply\secondepilength by 2%
\begin{picture}(0,0)%
\put(-\monolength,\secondmonolength){\line(1,-2){\bimolength}}%
\put(-\monolength,\secondmonolength){\ssehead}%
\put(\epilength,-\secondepilength){\ssehead}%
\put(\Z,-#3){\ssehead}%
\truex{200}\truey{800}\truez{600}%
\put(\value{x},\value{x}){\makebox(0,\value{z})[l]{${#1}$}}%
\put(-\value{x},-\value{y}){\makebox(0,\value{z})[r]{${#2}$}}%
\end{picture}}%

\newcommand{\SSEEQL}[3]{%
\Z=#3%
\divide\Z by 2%
\begin{picture}(0,0)%
\put(-\Z,#3){\begin{picture}(0,0)%
\truex{44}\truey{89}%
\put(-\value{y},-\value{x}){\line(1,-2){#3}}%
\put(\value{y},\value{x}){\line(1,-2){#3}}%
\end{picture}}%
\truex{200}\truey{800}\truez{600}%
\put(\value{x},\value{x}){\makebox(0,\value{z})[l]{${#1}$}}%
\put(-\value{x},-\value{y}){\makebox(0,\value{z})[r]{${#2}$}}%
\end{picture}}%

\newcommand{\SSEBIAR}[3]{%
\Y=#3%
\divide\Y by 2%
\Z=#3%
\multiply \Z by 2%
\begin{picture}(0,0)%
\put(-\Y,#3){\begin{picture}(0,0)%
\truex{313}\truey{156}%
\put(-\value{x},-\value{y}){\line(1,-2){#3}}%
\put(\value{x},\value{y}){\line(1,-2){#3}}%
\monolength=#3%
\advance\monolength by -\value{x}%
\epilength=#3%
\advance\epilength by \value{x}%
\secondmonolength=\Z%
\advance\secondmonolength by -\value{y}%
\secondepilength=\Z%
\advance\secondepilength by \value{y}%
\put(\monolength,-\secondepilength){\ssehead}%
\put(\epilength,-\secondmonolength){\ssehead}%
\end{picture}}
\truex{400}\truey{1000}\truez{600}%
\put(\value{x},\value{x}){\makebox(0,\value{z})[l]{${#1}$}}%
\put(-\value{x},-\value{y}){\makebox(0,\value{z})[r]{${#2}$}}%
\end{picture}}%

\newcommand{\SSEBIDIST}[3]{%
\Y=#3%
\divide\Y by 2%
\Z=#3%
\multiply \Z by 2%
\begin{picture}(0,0)%
\truex{313}\truey{156}\truez{400}%
\put(-\Y,#3){\begin{picture}(0,0)%
\put(-\value{x},-\value{y}){\line(1,-2){#3}}%
\put(\value{x},\value{y}){\line(1,-2){#3}}%
\monolength=#3%
\advance\monolength by -\value{x}%
\epilength=#3%
\advance\epilength by \value{x}%
\secondmonolength=\Z%
\advance\secondmonolength by -\value{y}%
\secondepilength=\Z%
\advance\secondepilength by \value{y}%
\put(\monolength,-\secondepilength){\ssehead}%
\put(\epilength,-\secondmonolength){\ssehead}%
\end{picture}}
\put(-\value{x},-\value{y}){\distsign}%
\put(\value{x},\value{y}){\distsign}%
\truex{500}\truey{1000}\truez{600}%
\put(\value{x},\value{x}){\makebox(0,\value{z})[l]{${#1}$}}%
\put(-\value{x},-\value{y}){\makebox(0,\value{z})[r]{${#2}$}}%
\end{picture}}%

\newcommand{\SSEADJAR}[3]{%
\Y=#3%
\divide\Y by 2%
\Z=#3%
\multiply \Z by 2%
\begin{picture}(0,0)%
\put(-\Y,#3){\begin{picture}(0,0)%
\truex{313}\truey{156}%
\monolength=#3%
\advance\monolength by -\value{x}%
\epilength=#3%
\advance\epilength by \value{x}%
\secondmonolength=\Z%
\advance\secondmonolength by -\value{y}%
\secondepilength=\Z%
\advance\secondepilength by \value{y}%
\put(-\value{x},-\value{y}){\line(1,-2){#3}}%
\put(\monolength,-\secondepilength){\ssehead}%
\put(\epilength,-\secondmonolength){\line(-1,2){#3}}%
\put(\value{x},\value{y}){\nnwhead}%
\end{picture}}
\truex{400}\truey{1000}\truez{600}%
\put(\value{x},\value{x}){\makebox(0,\value{z})[l]{${#1}$}}%
\put(-\value{x},-\value{y}){\makebox(0,\value{z})[r]{${#2}$}}%
\end{picture}}%

\newcommand{\SSEADJDIST}[3]{%
\Y=#3%
\divide\Y by 2%
\Z=#3%
\multiply \Z by 2%
\begin{picture}(0,0)%
\truex{313}\truey{156}\truez{400}%
\put(-\Y,#3){\begin{picture}(0,0)%
\monolength=#3%
\advance\monolength by -\value{x}%
\epilength=#3%
\advance\epilength by \value{x}%
\secondmonolength=\Z%
\advance\secondmonolength by -\value{y}%
\secondepilength=\Z%
\advance\secondepilength by \value{y}%
\put(-\value{x},-\value{y}){\line(1,-2){#3}}%
\put(\monolength,-\secondepilength){\ssehead}%
\put(\epilength,-\secondmonolength){\line(-1,2){#3}}%
\put(\value{x},\value{y}){\nnwhead}%
\end{picture}}
\put(\value{x},\value{y}){\distsign}%
\put(-\value{x},-\value{y}){\distsign}%
\truex{500}\truey{1000}\truez{600}%
\put(\value{x},\value{x}){\makebox(0,\value{z})[l]{${#1}$}}%
\put(-\value{x},-\value{y}){\makebox(0,\value{z})[r]{${#2}$}}%
\end{picture}}%

\def\basicssear[#1]{\SSEAR{}{}{#100}}%

\newcommand{\ssear}{\@ifnextchar[{\basicssear}{\basicssear[67]}}%

\def\basicSsear[#1]#2{\SSEAR{#2}{}{#100}}%

\newcommand{\Ssear}{\@ifnextchar[{\basicSsear}{\basicSsear[67]}}%

\def\basicsseaR[#1]#2{\SSEAR{}{#2}{#100}}%

\newcommand{\sseaR}{\@ifnextchar[{\basicsseaR}{\basicsseaR[67]}}%

\def\basicssedist[#1]{\SSEDIST{}{}{#100}}%

\newcommand{\ssedist}{\@ifnextchar[{\basicssedist}{\basicssedist[67]}}%

\def\basicSsedist[#1]#2{\SSEDIST{#2}{}{#100}}%

\newcommand{\Ssedist}{\@ifnextchar[{\basicSsedist}{\basicSsedist[67]}}%

\def\basicssedisT[#1]#2{\SSEDIST{}{#2}{#100}}%

\newcommand{\ssedisT}{\@ifnextchar[{\basicssedisT}{\basicssedisT[67]}}%

\def\basicssedotar[#1]{\SSEDOTAR{}{}{#100}}%

\newcommand{\ssedotar}{\@ifnextchar[{\basicssedotar}{\basicssedotar[67]}}%

\def\basicSsedotar[#1]#2{\SSEDOTAR{#2}{}{#100}}%

\newcommand{\Ssedotar}{\@ifnextchar[{\basicSsedotar}{\basicSsedotar[67]}}%

\def\basicssedotaR[#1]#2{\SSEDOTAR{}{#2}{#100}}%

\newcommand{\ssedotaR}{\@ifnextchar[{\basicssedotaR}{\basicssedotaR[67]}}%

\def\basicssemono[#1]{\SSEMONO{}{}{#100}}%

\newcommand{\ssemono}{\@ifnextchar[{\basicssemono}{\basicssemono[67]}}%

\def\basicSsemono[#1]#2{\SSEMONO{#2}{}{#100}}%

\newcommand{\Ssemono}{\@ifnextchar[{\basicSsemono}{\basicSsemono[67]}}%

\def\basicssemonO[#1]#2{\SSEMONO{}{#2}{#100}}%

\newcommand{\ssemonO}{\@ifnextchar[{\basicssemonO}{\basicssemonO[67]}}%

\def\basicsseepi[#1]{\SSEEPI{}{}{#100}}%

\newcommand{\sseepi}{\@ifnextchar[{\basicsseepi}{\basicsseepi[67]}}%

\def\basicSseepi[#1]#2{\SSEEPI{#2}{}{#100}}%

\newcommand{\Sseepi}{\@ifnextchar[{\basicSseepi}{\basicSseepi[67]}}%

\def\basicsseepI[#1]#2{\SSEEPI{}{#2}{#100}}%

\newcommand{\sseepI}{\@ifnextchar[{\basicsseepI}{\basicsseepI[67]}}%

\def\basicssebimo[#1]{\SSEBIMO{}{}{#100}}%

\newcommand{\ssebimo}{\@ifnextchar[{\basicssebimo}{\basicssebimo[67]}}%

\def\basicSsebimo[#1]#2{\SSEBIMO{#2}{}{#100}}%

\newcommand{\Ssebimo}{\@ifnextchar[{\basicSsebimo}{\basicSsebimo[67]}}%

\def\basicssebimO[#1]#2{\SSEBIMO{}{#2}{#100}}%

\newcommand{\ssebimO}{\@ifnextchar[{\basicssebimO}{\basicssebimO[67]}}%

\def\basicsseiso[#1]{\SSEAR{\isosign}{}{#100}}%

\newcommand{\sseiso}{\@ifnextchar[{\basicsseiso}{\basicsseiso[67]}}%

\def\basicSseiso[#1]#2{\SSEAR{#2}{\isosign}{#100}}%

\newcommand{\Sseiso}{\@ifnextchar[{\basicSseiso}{\basicSseiso[67]}}%

\def\basicsseisO[#1]#2{\SSEAR{\isosign}{#2}{#100}}%

\newcommand{\sseisO}{\@ifnextchar[{\basicsseisO}{\basicsseisO[67]}}%

\def\basicsseeql[#1]{\SSEEQL{}{}{#100}}%

\newcommand{\sseeql}{\@ifnextchar[{\basicsseeql}{\basicsseeql[67]}}%

\def\basicSseeql[#1]#2{\SSEEQL{#2}{}{#100}}%

\newcommand{\Sseeql}{\@ifnextchar[{\basicSseeql}{\basicSseeql[67]}}%

\def\basicsseeqL[#1]#2{\SSEEQL{}{#2}{#100}}%

\newcommand{\sseeqL}{\@ifnextchar[{\basicsseeqL}{\basicsseeqL[67]}}%

\def\basicssebiar[#1]{\SSEBIAR{}{}{#100}}%

\newcommand{\ssebiar}{\@ifnextchar[{\basicssebiar}{\basicssebiar[67]}}%

\def\basicSsebiar[#1]#2#3{\SSEBIAR{#2}{#3}{#100}}%

\newcommand{\Ssebiar}{\@ifnextchar[{\basicSsebiar}{\basicSsebiar[67]}}%

\def\basicssebidist[#1]{\SSEBIDIST{}{}{#100}}%

\newcommand{\ssebidist}{\@ifnextchar[{\basicssebidist}{\basicssebidist[67]}}%

\def\basicSsebidist[#1]#2#3{\SSEBIDIST{#2}{#3}{#100}}%

\newcommand{\Ssebidist}{\@ifnextchar[{\basicSsebidist}{\basicSsebidist[67]}}%

\def\basicsseadjar[#1]{\SSEADJAR{}{}{#100}}%

\newcommand{\sseadjar}{\@ifnextchar[{\basicsseadjar}{\basicsseadjar[67]}}%

\def\basicSseadjar[#1]#2#3{\SSEADJAR{#2}{#3}{#100}}%

\newcommand{\Sseadjar}{\@ifnextchar[{\basicSseadjar}{\basicSseadjar[67]}}%

\def\basicsseadjdist[#1]{\SSEADJDIST{}{}{#100}}%

\newcommand{\sseadjdist}{\@ifnextchar[{\basicsseadjdist}{\basicsseadjdist[67]}}%

\def\basicSseadjdist[#1]#2#3{\SSEADJDIST{#2}{#3}{#100}}%

\newcommand{\Sseadjdist}{\@ifnextchar[{\basicSseadjdist}{\basicSseadjdist[67]}}%

\newcommand{\SSWAR}[3]{%
\Z=#3%
\divide\Z by 2%
\begin{picture}(0,0)%
\put(\Z,#3){\line(-1,-2){#3}}%
\put(-\Z,-#3){\sswhead}%
\truex{200}\truey{800}\truez{600}%
\put(-\value{x},\value{x}){\makebox(0,\value{z})[r]{${#1}$}}%
\put(\value{x},-\value{y}){\makebox(0,\value{z})[l]{${#2}$}}%
\end{picture}}%

\newcommand{\SSWDIST}[3]{%
\Z=#3%
\divide\Z by 2%
\begin{picture}(0,0)%
\put(\Z,#3){\line(-1,-2){#3}}%
\put(-\Z,-#3){\sswhead}%
\truex{400}%
\put(0,0){\distsign}%
\truex{200}\truey{800}\truez{600}%
\put(-\value{x},\value{x}){\makebox(0,\value{z})[r]{${#1}$}}%
\put(\value{x},-\value{y}){\makebox(0,\value{z})[l]{${#2}$}}%
\end{picture}}%

\newcommand{\SSWDOTAR}[3]{%
\truex{100}\truey{268}\truez{134}%
\Z=#3%
\divide\Z by 2%
\NUMBEROFDOTS=#3%
\divide\NUMBEROFDOTS by \value{z}%
\advance\NUMBEROFDOTS by 1%
\begin{picture}(0,0)%
\multiput(\Z,#3)(-\value{z},-\value{y}){\NUMBEROFDOTS}%
{\circle*{\value{x}}}%
\put(-\Z,-#3){\sswhead}%
\truex{200}\truey{800}\truez{600}%
\put(-\value{x},\value{x}){\makebox(0,\value{z})[r]{${#1}$}}%
\put(\value{x},-\value{y}){\makebox(0,\value{z})[l]{${#2}$}}%
\end{picture}}%

\newcommand{\SSWMONO}[3]{%
\Z=#3%
\divide\Z by 2%
\truetaiL%
\bimolength=#3%
\advance\bimolength by -\truemonotaiL%
\monolength=\bimolength%
\advance\monolength by -\Z%
\secondmonolength=\monolength%
\multiply\secondmonolength by 2%
\begin{picture}(0,0)%
\put(\monolength,\secondmonolength){\line(-1,-2){\bimolength}}%
\put(\monolength,\secondmonolength){\sswhead}%
\put(-\Z,-#3){\sswhead}%
\truex{200}\truey{800}\truez{600}%
\put(-\value{x},\value{x}){\makebox(0,\value{z})[r]{${#1}$}}%
\put(\value{x},-\value{y}){\makebox(0,\value{z})[l]{${#2}$}}%
\end{picture}}%

\newcommand{\SSWEPI}[3]{%
\Z=#3%
\divide\Z by 2%
\trueheaD%
\bimolength=#3%
\advance\bimolength by -\trueepiheaD%
\epilength=\bimolength%
\advance\epilength by -\Z%
\secondepilength=\epilength%
\multiply\secondepilength by 2%
\begin{picture}(0,0)%
\put(\Z,#3){\line(-1,-2){\bimolength}}%
\put(-\epilength,-\secondepilength){\sswhead}%
\put(-\Z,-#3){\sswhead}%
\truex{200}\truey{800}\truez{600}%
\put(-\value{x},\value{x}){\makebox(0,\value{z})[r]{${#1}$}}%
\put(\value{x},-\value{y}){\makebox(0,\value{z})[l]{${#2}$}}%
\end{picture}}%

\newcommand{\SSWBIMO}[3]{%
\Z=#3%
\divide\Z by 2%
\truetaiL\trueheaD%
\bimolength=#3%
\advance\bimolength by -\truemonotaiL%
\monolength=\bimolength%
\advance\monolength by -\Z%
\advance\bimolength by -\trueepiheaD%
\epilength=\bimolength%
\advance\epilength by -\monolength%
\secondmonolength=\monolength%
\multiply\secondmonolength by 2%
\secondepilength=\epilength%
\multiply\secondepilength by 2%
\begin{picture}(0,0)%
\put(\monolength,\secondmonolength){\line(-1,-2){\bimolength}}%
\put(\monolength,\secondmonolength){\sswhead}%
\put(-\epilength,-\secondepilength){\sswhead}%
\put(-\Z,-#3){\sswhead}%
\truex{200}\truey{800}\truez{600}%
\put(-\value{x},\value{x}){\makebox(0,\value{z})[r]{${#1}$}}%
\put(\value{x},-\value{y}){\makebox(0,\value{z})[l]{${#2}$}}%
\end{picture}}%

\newcommand{\SSWBIAR}[3]{%
\Y=#3%
\divide\Y by 2%
\Z=#3%
\multiply \Z by 2%
\begin{picture}(0,0)%
\put(\Y,#3){\begin{picture}(0,0)%
\truex{313}\truey{156}%
\put(-\value{x},\value{y}){\line(-1,-2){#3}}%
\put(\value{x},-\value{y}){\line(-1,-2){#3}}%
\monolength=#3%
\advance\monolength by -\value{x}%
\epilength=#3%
\advance\epilength by \value{x}%
\secondmonolength=\Z%
\advance\secondmonolength by -\value{y}%
\secondepilength=\Z%
\advance\secondepilength by \value{y}%
\put(-\monolength,-\secondepilength){\sswhead}%
\put(-\epilength,-\secondmonolength){\sswhead}%
\end{picture}}
\truex{300}\truey{1000}\truez{600}%
\put(-\value{x},\value{x}){\makebox(0,\value{z})[r]{${#1}$}}%
\put(\value{x},-\value{y}){\makebox(0,\value{z})[l]{${#2}$}}%
\end{picture}}%

\newcommand{\SSWBIDIST}[3]{%
\Y=#3%
\divide\Y by 2%
\Z=#3%
\multiply \Z by 2%
\begin{picture}(0,0)%
\truex{313}\truey{156}\truez{400}%
\put(\Y,#3){\begin{picture}(0,0)%
\put(-\value{x},\value{y}){\line(-1,-2){#3}}%
\put(\value{x},-\value{y}){\line(-1,-2){#3}}%
\monolength=#3%
\advance\monolength by -\value{x}%
\epilength=#3%
\advance\epilength by \value{x}%
\secondmonolength=\Z%
\advance\secondmonolength by -\value{y}%
\secondepilength=\Z%
\advance\secondepilength by \value{y}%
\put(-\monolength,-\secondepilength){\sswhead}%
\put(-\epilength,-\secondmonolength){\sswhead}%
\end{picture}}
\put(-\value{x},\value{y}){\distsign}%
\put(\value{x},-\value{y}){\distsign}%
\truex{300}\truey{1000}\truez{600}%
\put(-\value{x},\value{x}){\makebox(0,\value{z})[r]{${#1}$}}%
\put(\value{x},-\value{y}){\makebox(0,\value{z})[l]{${#2}$}}%
\end{picture}}%

\newcommand{\SSWADJAR}[3]{%
\Y=#3%
\divide\Y by 2%
\Z=#3%
\multiply \Z by 2%
\begin{picture}(0,0)%
\put(\Y,#3){\begin{picture}(0,0)%
\truex{313}\truey{156}%
\monolength=#3%
\advance\monolength by -\value{x}%
\epilength=#3%
\advance\epilength by \value{x}%
\secondmonolength=\Z%
\advance\secondmonolength by -\value{y}%
\secondepilength=\Z%
\advance\secondepilength by \value{y}%
\put(\value{x},-\value{y}){\line(-1,-2){#3}}%
\put(-\monolength,-\secondepilength){\sswhead}%
\put(-\epilength,-\secondmonolength){\line(1,2){#3}}%
\put(-\value{x},\value{y}){\nnehead}%
\end{picture}}
\truex{300}\truey{1000}\truez{600}%
\put(-\value{x},\value{x}){\makebox(0,\value{z})[r]{${#1}$}}%
\put(\value{x},-\value{y}){\makebox(0,\value{z})[l]{${#2}$}}%
\end{picture}}%

\newcommand{\SSWADJDIST}[3]{%
\Y=#3%
\divide\Y by 2%
\Z=#3%
\multiply \Z by 2%
\begin{picture}(0,0)%
\truex{313}\truey{156}\truez{400}%
\put(\Y,#3){\begin{picture}(0,0)%
\monolength=#3%
\advance\monolength by -\value{x}%
\epilength=#3%
\advance\epilength by \value{x}%
\secondmonolength=\Z%
\advance\secondmonolength by -\value{y}%
\secondepilength=\Z%
\advance\secondepilength by \value{y}%
\put(\value{x},-\value{y}){\line(-1,-2){#3}}%
\put(-\monolength,-\secondepilength){\sswhead}%
\put(-\epilength,-\secondmonolength){\line(1,2){#3}}%
\put(-\value{x},\value{y}){\nnehead}%
\end{picture}}
\put(-\value{x},\value{y}){\distsign}%
\put(\value{x},-\value{y}){\distsign}%
\truex{300}\truey{1000}\truez{600}%
\put(-\value{x},\value{x}){\makebox(0,\value{z})[r]{${#1}$}}%
\put(\value{x},-\value{y}){\makebox(0,\value{z})[l]{${#2}$}}%
\end{picture}}%

\def\basicsswar[#1]{\SSWAR{}{}{#100}}%

\newcommand{\sswar}{\@ifnextchar[{\basicsswar}{\basicsswar[67]}}%

\def\basicSswar[#1]#2{\SSWAR{#2}{}{#100}}%

\newcommand{\Sswar}{\@ifnextchar[{\basicSswar}{\basicSswar[67]}}%

\def\basicsswaR[#1]#2{\SSWAR{}{#2}{#100}}%

\newcommand{\sswaR}{\@ifnextchar[{\basicsswaR}{\basicsswaR[67]}}%

\def\basicsswdist[#1]{\SSWDIST{}{}{#100}}%

\newcommand{\sswdist}{\@ifnextchar[{\basicsswdist}{\basicsswdist[67]}}%

\def\basicSswdist[#1]#2{\SSWDIST{#2}{}{#100}}%

\newcommand{\Sswdist}{\@ifnextchar[{\basicSswdist}{\basicSswdist[67]}}%

\def\basicsswdisT[#1]#2{\SSWDIST{}{#2}{#100}}%

\newcommand{\sswdisT}{\@ifnextchar[{\basicsswdisT}{\basicsswdisT[67]}}%

\def\basicsswdotar[#1]{\SSWDOTAR{}{}{#100}}%

\newcommand{\sswdotar}{\@ifnextchar[{\basicsswdotar}{\basicsswdotar[67]}}%

\def\basicSswdotar[#1]#2{\SSWDOTAR{#2}{}{#100}}%

\newcommand{\Sswdotar}{\@ifnextchar[{\basicSswdotar}{\basicSswdotar[67]}}%

\def\basicsswdotaR[#1]#2{\SSWDOTAR{}{#2}{#100}}%

\newcommand{\sswdotaR}{\@ifnextchar[{\basicsswdotaR}{\basicsswdotaR[67]}}%

\def\basicsswmono[#1]{\SSWMONO{}{}{#100}}%

\newcommand{\sswmono}{\@ifnextchar[{\basicsswmono}{\basicsswmono[67]}}%

\def\basicSswmono[#1]#2{\SSWMONO{#2}{}{#100}}%

\newcommand{\Sswmono}{\@ifnextchar[{\basicSswmono}{\basicSswmono[67]}}%

\def\basicsswmonO[#1]#2{\SSWMONO{}{#2}{#100}}%

\newcommand{\sswmonO}{\@ifnextchar[{\basicsswmonO}{\basicsswmonO[67]}}%

\def\basicsswepi[#1]{\SSWEPI{}{}{#100}}%

\newcommand{\sswepi}{\@ifnextchar[{\basicsswepi}{\basicsswepi[67]}}%

\def\basicSswepi[#1]#2{\SSWEPI{#2}{}{#100}}%

\newcommand{\Sswepi}{\@ifnextchar[{\basicSswepi}{\basicSswepi[67]}}%

\def\basicsswepI[#1]#2{\SSWEPI{}{#2}{#100}}%

\newcommand{\sswepI}{\@ifnextchar[{\basicsswepI}{\basicsswepI[67]}}%

\def\basicsswbimo[#1]{\SSWBIMO{}{}{#100}}%

\newcommand{\sswbimo}{\@ifnextchar[{\basicsswbimo}{\basicsswbimo[67]}}%

\def\basicSswbimo[#1]#2{\SSWBIMO{#2}{}{#100}}%

\newcommand{\Sswbimo}{\@ifnextchar[{\basicSswbimo}{\basicSswbimo[67]}}%

\def\basicsswbimO[#1]#2{\SSWBIMO{}{#2}{#100}}%

\newcommand{\sswbimO}{\@ifnextchar[{\basicsswbimO}{\basicsswbimO[67]}}%

\def\basicsswiso[#1]{\SSWAR{\isosign}{}{#100}}%

\newcommand{\sswiso}{\@ifnextchar[{\basicsswiso}{\basicsswiso[67]}}%

\def\basicSswiso[#1]#2{\SSWAR{#2}{\isosign}{#100}}%

\newcommand{\Sswiso}{\@ifnextchar[{\basicSswiso}{\basicSswiso[67]}}%

\def\basicsswisO[#1]#2{\SSWAR{\isosign}{#2}{#100}}%

\newcommand{\sswisO}{\@ifnextchar[{\basicsswisO}{\basicsswisO[67]}}%

\def\basicsswbiar[#1]{\SSWBIAR{}{}{#100}}%

\newcommand{\sswbiar}{\@ifnextchar[{\basicsswbiar}{\basicsswbiar[67]}}%

\def\basicSswbiar[#1]#2#3{\SSWBIAR{#2}{#3}{#100}}%

\newcommand{\Sswbiar}{\@ifnextchar[{\basicSswbiar}{\basicSswbiar[67]}}%

\def\basicsswbidist[#1]{\SSWBIDIST{}{}{#100}}%

\newcommand{\sswbidist}{\@ifnextchar[{\basicsswbidist}{\basicsswbidist[67]}}%

\def\basicSswbidist[#1]#2#3{\SSWBIDIST{#2}{#3}{#100}}%

\newcommand{\Sswbidist}{\@ifnextchar[{\basicSswbidist}{\basicSswbidist[67]}}%

\def\basicsswadjar[#1]{\SSWADJAR{}{}{#100}}%

\newcommand{\sswadjar}{\@ifnextchar[{\basicsswadjar}{\basicsswadjar[67]}}%

\def\basicSswadjar[#1]#2#3{\SSWADJAR{#2}{#3}{#100}}%

\newcommand{\Sswadjar}{\@ifnextchar[{\basicSswadjar}{\basicSswadjar[67]}}%

\def\basicsswadjdist[#1]{\SSWADJDIST{}{}{#100}}%

\newcommand{\sswadjdist}{\@ifnextchar[{\basicsswadjdist}{\basicsswadjdist[67]}}%

\def\basicSswadjdist[#1]#2#3{\SSWADJDIST{#2}{#3}{#100}}%

\newcommand{\Sswadjdist}{\@ifnextchar[{\basicSswadjdist}{\basicSswadjdist[67]}}%

\newcommand{\NNWAR}[3]{%
\Z=#3%
\divide\Z by 2%
\begin{picture}(0,0)%
\put(\Z,-#3){\line(-1,2){#3}}%
\put(-\Z,#3){\nnwhead}%
\truex{200}\truey{800}\truez{600}%
\put(\value{x},\value{x}){\makebox(0,\value{z})[l]{${#1}$}}%
\put(-\value{x},-\value{y}){\makebox(0,\value{z})[r]{${#2}$}}%
\end{picture}}%

\newcommand{\NNWDIST}[3]{%
\Z=#3%
\divide\Z by 2%
\begin{picture}(0,0)%
\put(\Z,-#3){\line(-1,2){#3}}%
\put(-\Z,#3){\nnwhead}%
\truex{400}%
\put(0,0){\distsign}%
\truex{200}\truey{800}\truez{600}%
\put(\value{x},\value{x}){\makebox(0,\value{z})[l]{${#1}$}}%
\put(-\value{x},-\value{y}){\makebox(0,\value{z})[r]{${#2}$}}%
\end{picture}}%

\newcommand{\NNWDOTAR}[3]{%
\truex{100}\truey{268}\truez{134}%
\Z=#3%
\divide\Z by 2%
\NUMBEROFDOTS=#3%
\divide\NUMBEROFDOTS by \value{z}%
\advance\NUMBEROFDOTS by 1%
\begin{picture}(0,0)%
\multiput(\Z,-#3)(-\value{z},\value{y}){\NUMBEROFDOTS}%
{\circle*{\value{x}}}%
\put(-\Z,#3){\nnwhead}%
\truex{200}\truey{800}\truez{600}%
\put(\value{x},\value{x}){\makebox(0,\value{z})[l]{${#1}$}}%
\put(-\value{x},-\value{y}){\makebox(0,\value{z})[r]{${#2}$}}%
\end{picture}}%

\newcommand{\NNWMONO}[3]{%
\Z=#3%
\divide\Z by 2%
\truetaiL%
\bimolength=#3%
\advance\bimolength by -\truemonotaiL%
\monolength=\bimolength%
\advance\monolength by -\Z%
\secondmonolength=\monolength%
\multiply\secondmonolength by 2%
\begin{picture}(0,0)%
\put(\monolength,-\secondmonolength){\line(-1,2){\bimolength}}%
\put(\monolength,-\secondmonolength){\nnwhead}%
\put(-\Z,#3){\nnwhead}%
\truex{200}\truey{800}\truez{600}%
\put(\value{x},\value{x}){\makebox(0,\value{z})[l]{${#1}$}}%
\put(-\value{x},-\value{y}){\makebox(0,\value{z})[r]{${#2}$}}%
\end{picture}}%

\newcommand{\NNWEPI}[3]{%
\Z=#3%
\divide\Z by 2%
\trueheaD%
\bimolength=#3%
\advance\bimolength by -\trueepiheaD%
\epilength=\bimolength%
\advance\epilength by -\Z%
\secondepilength=\epilength%
\multiply\secondepilength by 2%
\begin{picture}(0,0)%
\put(\Z,-#3){\line(-1,2){\bimolength}}%
\put(-\epilength,\secondepilength){\nnwhead}%
\put(-\Z,#3){\nnwhead}%
\truex{200}\truey{800}\truez{600}%
\put(\value{x},\value{x}){\makebox(0,\value{z})[l]{${#1}$}}%
\put(-\value{x},-\value{y}){\makebox(0,\value{z})[r]{${#2}$}}%
\end{picture}}%

\newcommand{\NNWBIMO}[3]{%
\Z=#3%
\divide\Z by 2%
\truetaiL\trueheaD%
\bimolength=#3%
\advance\bimolength by -\truemonotaiL%
\monolength=\bimolength%
\advance\monolength by -\Z%
\advance\bimolength by -\trueepiheaD%
\epilength=\bimolength%
\advance\epilength by -\monolength%
\secondmonolength=\monolength%
\multiply\secondmonolength by 2%
\secondepilength=\epilength%
\multiply\secondepilength by 2%
\begin{picture}(0,0)%
\put(\monolength,-\secondmonolength){\line(-1,2){\bimolength}}%
\put(\monolength,-\secondmonolength){\nnwhead}%
\put(-\epilength,\secondepilength){\nnwhead}%
\put(-\Z,#3){\nnwhead}%
\truex{200}\truey{800}\truez{600}%
\put(\value{x},\value{x}){\makebox(0,\value{z})[l]{${#1}$}}%
\put(-\value{x},-\value{y}){\makebox(0,\value{z})[r]{${#2}$}}%
\end{picture}}%

\newcommand{\NNWBIAR}[3]{%
\Y=#3%
\divide\Y by 2%
\Z=#3%
\multiply \Z by 2%
\begin{picture}(0,0)%
\put(\Y,-#3){\begin{picture}(0,0)%
\truex{313}\truey{156}%
\put(-\value{x},-\value{y}){\line(-1,2){#3}}%
\put(\value{x},\value{y}){\line(-1,2){#3}}%
\monolength=#3%
\advance\monolength by -\value{x}%
\epilength=#3%
\advance\epilength by \value{x}%
\secondmonolength=\Z%
\advance\secondmonolength by -\value{y}%
\secondepilength=\Z%
\advance\secondepilength by \value{y}%
\put(-\monolength,\secondepilength){\nnwhead}%
\put(-\epilength,\secondmonolength){\nnwhead}%
\end{picture}}
\truex{400}\truey{1000}\truez{600}%
\put(\value{x},\value{x}){\makebox(0,\value{z})[l]{${#1}$}}%
\put(-\value{x},-\value{y}){\makebox(0,\value{z})[r]{${#2}$}}%
\end{picture}}%

\newcommand{\NNWBIDIST}[3]{%
\Y=#3%
\divide\Y by 2%
\Z=#3%
\multiply \Z by 2%
\begin{picture}(0,0)%
\truex{313}\truey{156}\truez{400}%
\put(\Y,-#3){\begin{picture}(0,0)%
\put(-\value{x},-\value{y}){\line(-1,2){#3}}%
\put(\value{x},\value{y}){\line(-1,2){#3}}%
\monolength=#3%
\advance\monolength by -\value{x}%
\epilength=#3%
\advance\epilength by \value{x}%
\secondmonolength=\Z%
\advance\secondmonolength by -\value{y}%
\secondepilength=\Z%
\advance\secondepilength by \value{y}%
\put(-\monolength,\secondepilength){\nnwhead}%
\put(-\epilength,\secondmonolength){\nnwhead}%
\end{picture}}
\put(-\value{x},-\value{y}){\distsign}%
\put(\value{x},\value{y}){\distsign}%
\truex{500}\truey{1000}\truez{600}%
\put(\value{x},\value{x}){\makebox(0,\value{z})[l]{${#1}$}}%
\put(-\value{x},-\value{y}){\makebox(0,\value{z})[r]{${#2}$}}%
\end{picture}}%

\newcommand{\NNWADJAR}[3]{%
\Y=#3%
\divide\Y by 2%
\Z=#3%
\multiply \Z by 2%
\begin{picture}(0,0)%
\put(\Y,-#3){\begin{picture}(0,0)%
\truex{313}\truey{156}%
\monolength=#3%
\advance\monolength by -\value{x}%
\epilength=#3%
\advance\epilength by \value{x}%
\secondmonolength=\Z%
\advance\secondmonolength by -\value{y}%
\secondepilength=\Z%
\advance\secondepilength by \value{y}%
\put(-\value{x},-\value{y}){\line(-1,2){#3}}%
\put(-\epilength,\secondmonolength){\nnwhead}%
\put(-\monolength,\secondepilength){\line(1,-2){#3}}%
\put(\value{x},\value{y}){\ssehead}%
\end{picture}}
\truex{400}\truey{1000}\truez{600}%
\put(\value{x},\value{x}){\makebox(0,\value{z})[l]{${#1}$}}%
\put(-\value{x},-\value{y}){\makebox(0,\value{z})[r]{${#2}$}}%
\end{picture}}%

\newcommand{\NNWADJDIST}[3]{%
\Y=#3%
\divide\Y by 2%
\Z=#3%
\multiply \Z by 2%
\begin{picture}(0,0)%
\truex{313}\truey{156}\truez{400}%
\put(\Y,-#3){\begin{picture}(0,0)%
\monolength=#3%
\advance\monolength by -\value{x}%
\epilength=#3%
\advance\epilength by \value{x}%
\secondmonolength=\Z%
\advance\secondmonolength by -\value{y}%
\secondepilength=\Z%
\advance\secondepilength by \value{y}%
\put(-\value{x},-\value{y}){\line(-1,2){#3}}%
\put(-\epilength,\secondmonolength){\nnwhead}%
\put(-\monolength,\secondepilength){\line(1,-2){#3}}%
\put(\value{x},\value{y}){\ssehead}%
\end{picture}}
\put(\value{x},\value{y}){\distsign}%
\put(-\value{x},-\value{y}){\distsign}%
\truex{500}\truey{1000}\truez{600}%
\put(\value{x},\value{x}){\makebox(0,\value{z})[l]{${#1}$}}%
\put(-\value{x},-\value{y}){\makebox(0,\value{z})[r]{${#2}$}}%
\end{picture}}%

\def\basicnnwar[#1]{\NNWAR{}{}{#100}}%

\newcommand{\nnwar}{\@ifnextchar[{\basicnnwar}{\basicnnwar[67]}}%

\def\basicNnwar[#1]#2{\NNWAR{#2}{}{#100}}%

\newcommand{\Nnwar}{\@ifnextchar[{\basicNnwar}{\basicNnwar[67]}}%

\def\basicnnwaR[#1]#2{\NNWAR{}{#2}{#100}}%

\newcommand{\nnwaR}{\@ifnextchar[{\basicnnwaR}{\basicnnwaR[67]}}%

\def\basicnnwdist[#1]{\NNWDIST{}{}{#100}}%

\newcommand{\nnwdist}{\@ifnextchar[{\basicnnwdist}{\basicnnwdist[67]}}%

\def\basicNnwdist[#1]#2{\NNWDIST{#2}{}{#100}}%

\newcommand{\Nnwdist}{\@ifnextchar[{\basicNnwdist}{\basicNnwdist[67]}}%

\def\basicnnwdisT[#1]#2{\NNWDIST{}{#2}{#100}}%

\newcommand{\nnwdisT}{\@ifnextchar[{\basicnnwdisT}{\basicnnwdisT[67]}}%

\def\basicnnwdotar[#1]{\NNWDOTAR{}{}{#100}}%

\newcommand{\nnwdotar}{\@ifnextchar[{\basicnnwdotar}{\basicnnwdotar[67]}}%

\def\basicNnwdotar[#1]#2{\NNWDOTAR{#2}{}{#100}}%

\newcommand{\Nnwdotar}{\@ifnextchar[{\basicNnwdotar}{\basicNnwdotar[67]}}%

\def\basicnnwdotaR[#1]#2{\NNWDOTAR{}{#2}{#100}}%

\newcommand{\nnwdotaR}{\@ifnextchar[{\basicnnwdotaR}{\basicnnwdotaR[67]}}%

\def\basicnnwmono[#1]{\NNWMONO{}{}{#100}}%

\newcommand{\nnwmono}{\@ifnextchar[{\basicnnwmono}{\basicnnwmono[67]}}%

\def\basicNnwmono[#1]#2{\NNWMONO{#2}{}{#100}}%

\newcommand{\Nnwmono}{\@ifnextchar[{\basicNnwmono}{\basicNnwmono[67]}}%

\def\basicnnwmonO[#1]#2{\NNWMONO{}{#2}{#100}}%

\newcommand{\nnwmonO}{\@ifnextchar[{\basicnnwmonO}{\basicnnwmonO[67]}}%

\def\basicnnwepi[#1]{\NNWEPI{}{}{#100}}%

\newcommand{\nnwepi}{\@ifnextchar[{\basicnnwepi}{\basicnnwepi[67]}}%

\def\basicNnwepi[#1]#2{\NNWEPI{#2}{}{#100}}%

\newcommand{\Nnwepi}{\@ifnextchar[{\basicNnwepi}{\basicNnwepi[67]}}%

\def\basicnnwepI[#1]#2{\NNWEPI{}{#2}{#100}}%

\newcommand{\nnwepI}{\@ifnextchar[{\basicnnwepI}{\basicnnwepI[67]}}%

\def\basicnnwbimo[#1]{\NNWBIMO{}{}{#100}}%

\newcommand{\nnwbimo}{\@ifnextchar[{\basicnnwbimo}{\basicnnwbimo[67]}}%

\def\basicNnwbimo[#1]#2{\NNWBIMO{#2}{}{#100}}%

\newcommand{\Nnwbimo}{\@ifnextchar[{\basicNnwbimo}{\basicNnwbimo[67]}}%

\def\basicnnwbimO[#1]#2{\NNWBIMO{}{#2}{#100}}%

\newcommand{\nnwbimO}{\@ifnextchar[{\basicnnwbimO}{\basicnnwbimO[67]}}%

\def\basicnnwiso[#1]{\NNWAR{\isosign}{}{#100}}%

\newcommand{\nnwiso}{\@ifnextchar[{\basicnnwiso}{\basicnnwiso[67]}}%

\def\basicNnwiso[#1]#2{\NNWAR{#2}{\isosign}{#100}}%

\newcommand{\Nnwiso}{\@ifnextchar[{\basicNnwiso}{\basicNnwiso[67]}}%

\def\basicnnwisO[#1]#2{\NNWAR{\isosign}{#2}{#100}}%

\newcommand{\nnwisO}{\@ifnextchar[{\basicnnwisO}{\basicnnwisO[67]}}%

\def\basicnnwbiar[#1]{\NNWBIAR{}{}{#100}}%

\newcommand{\nnwbiar}{\@ifnextchar[{\basicnnwbiar}{\basicnnwbiar[67]}}%

\def\basicNnwbiar[#1]#2#3{\NNWBIAR{#2}{#3}{#100}}%

\newcommand{\Nnwbiar}{\@ifnextchar[{\basicNnwbiar}{\basicNnwbiar[67]}}%

\def\basicnnwbidist[#1]{\NNWBIDIST{}{}{#100}}%

\newcommand{\nnwbidist}{\@ifnextchar[{\basicnnwbidist}{\basicnnwbidist[67]}}%

\def\basicNnwbidist[#1]#2#3{\NNWBIDIST{#2}{#3}{#100}}%

\newcommand{\Nnwbidist}{\@ifnextchar[{\basicNnwbidist}{\basicNnwbidist[67]}}%

\def\basicnnwadjar[#1]{\NNWADJAR{}{}{#100}}%

\newcommand{\nnwadjar}{\@ifnextchar[{\basicnnwadjar}{\basicnnwadjar[67]}}%

\def\basicNnwadjar[#1]#2#3{\NNWADJAR{#2}{#3}{#100}}%

\newcommand{\Nnwadjar}{\@ifnextchar[{\basicNnwadjar}{\basicNnwadjar[67]}}%

\def\basicnnwadjdist[#1]{\NNWADJDIST{}{}{#100}}%

\newcommand{\nnwadjdist}{\@ifnextchar[{\basicnnwadjdist}{\basicnnwadjdist[67]}}%

\def\basicNnwadjdist[#1]#2#3{\NNWADJDIST{#2}{#3}{#100}}%

\newcommand{\Nnwadjdist}{\@ifnextchar[{\basicNnwadjdist}{\basicNnwadjdist[67]}}%

\newcommand{\EENEAR}[3]{%
\Y=#3%
\divide \Y by 2%
\Z=\Y%
\divide \Z by 3%
\begin{picture}(0,0)%
\put(-\Y,-\Z){\line(3,1){#3}}%
\put(\Y,\Z){\eenehead}%
\truex{200}\truey{800}\truez{600}%
\put(-\value{x},\value{x}){\makebox(0,\value{z})[r]{${#1}$}}%
\put(\value{x},-\value{y}){\makebox(0,\value{z})[l]{${#2}$}}%
\end{picture}}%

\def\basiceenear[#1]{\EENEAR{}{}{#100}}%

\newcommand{\eenear}{\@ifnextchar[{\basiceenear}{\basiceenear[211]}}%

\def\basicEenear[#1]#2{\EENEAR{#2}{}{#100}}%

\newcommand{\Eenear}{\@ifnextchar[{\basicEenear}{\basicEenear[211]}}%

\def\basiceeneaR[#1]#2{\EENEAR{}{#2}{#100}}%

\newcommand{\eeneaR}{\@ifnextchar[{\basiceeneaR}{\basiceeneaR[211]}}%

\newcommand{\EESEAR}[3]{%
\Y=#3%
\divide \Y by 2%
\Z=\Y%
\divide \Z by 3%
\begin{picture}(0,0)%
\put(-\Y,\Z){\line(3,-1){#3}}%
\put(\Y,-\Z){\eesehead}%
\truex{200}\truey{800}\truez{600}%
\put(\value{x},\value{x}){\makebox(0,\value{z})[l]{${#1}$}}%
\put(-\value{x},-\value{y}){\makebox(0,\value{z})[r]{${#2}$}}%
\end{picture}}%

\def\basiceesear[#1]{\EESEAR{}{}{#100}}%

\newcommand{\eesear}{\@ifnextchar[{\basiceesear}{\basiceesear[211]}}%

\def\basicEesear[#1]#2{\EESEAR{#2}{}{#100}}%

\newcommand{\Eesear}{\@ifnextchar[{\basicEesear}{\basicEesear[211]}}%

\def\basiceeseaR[#1]#2{\EESEAR{}{#2}{#100}}%

\newcommand{\eeseaR}{\@ifnextchar[{\basiceeseaR}{\basiceeseaR[211]}}%

\newcommand{\WWNWAR}[3]{%
\Y=#3%
\divide \Y by 2%
\Z=\Y%
\divide \Z by 3%
\begin{picture}(0,0)%
\put(\Y,-\Z){\line(-3,1){#3}}%
\put(-\Y,\Z){\wwnwhead}%
\truex{200}\truey{800}\truez{600}%
\put(\value{x},\value{x}){\makebox(0,\value{z})[l]{${#1}$}}%
\put(-\value{x},-\value{y}){\makebox(0,\value{z})[r]{${#2}$}}%
\end{picture}}%

\def\basicwwnwar[#1]{\WWNWAR{}{}{#100}}%

\newcommand{\wwnwar}{\@ifnextchar[{\basicwwnwar}{\basicwwnwar[211]}}%

\def\basicWwnwar[#1]#2{\WWNWAR{#2}{}{#100}}%

\newcommand{\Wwnwar}{\@ifnextchar[{\basicWwnwar}{\basicWwnwar[211]}}%

\def\basicwwnwaR[#1]#2{\WWNWAR{}{#2}{#100}}%

\newcommand{\wwnwaR}{\@ifnextchar[{\basicwwnwaR}{\basicwwnwaR[211]}}%

\newcommand{\WWSWAR}[3]{%
\Y=#3%
\divide \Y by 2%
\Z=\Y%
\divide \Z by 3%
\begin{picture}(0,0)%
\put(\Y,\Z){\line(-3,-1){#3}}%
\put(-\Y,-\Z){\wwswhead}%
\truex{200}\truey{800}\truez{600}%
\put(-\value{x},\value{x}){\makebox(0,\value{z})[r]{${#1}$}}%
\put(\value{x},-\value{y}){\makebox(0,\value{z})[l]{${#2}$}}%
\end{picture}}%

\def\basicwwswar[#1]{\WWSWAR{}{}{#100}}%

\newcommand{\wwswar}{\@ifnextchar[{\basicwwswar}{\basicwwswar[211]}}%

\def\basicWwswar[#1]#2{\WWSWAR{#2}{}{#100}}%

\newcommand{\Wwswar}{\@ifnextchar[{\basicWwswar}{\basicWwswar[211]}}%

\def\basicwwswaR[#1]#2{\WWSWAR{}{#2}{#100}}%

\newcommand{\wwswaR}{\@ifnextchar[{\basicwwswaR}{\basicwwswaR[211]}}%

\newcommand{\NNNEAR}[3]{%
\Y=#3%
\divide \Y by 2%
\Z=\Y%
\multiply \Z by 3%
\begin{picture}(0,0)%
\put(-\Y,-\Z){\line(1,3){#3}}%
\put(\Y,\Z){\nnnehead}%
\truex{100}\truez{600}%
\put(-\value{x},\value{x}){\makebox(0,\value{z})[r]{${#1}$}}%
\put(\value{x},-\value{z}){\makebox(0,\value{z})[l]{${#2}$}}%
\end{picture}}%

\def\basicnnnear[#1]{\NNNEAR{}{}{#100}}%

\newcommand{\nnnear}{\@ifnextchar[{\basicnnnear}{\basicnnnear[71]}}%

\def\basicNnnear[#1]#2{\NNNEAR{#2}{}{#100}}%

\newcommand{\Nnnear}{\@ifnextchar[{\basicNnnear}{\basicNnnear[71]}}%

\def\basicnnneaR[#1]#2{\NNNEAR{}{#2}{#100}}%

\newcommand{\nnneaR}{\@ifnextchar[{\basicnnneaR}{\basicnnneaR[71]}}%

\newcommand{\SSSWAR}[3]{%
\Y=#3%
\divide \Y by 2%
\Z=\Y%
\multiply \Z by 3%
\begin{picture}(0,0)%
\put(\Y,\Z){\line(-1,-3){#3}}%
\put(-\Y,-\Z){\ssswhead}%
\truex{100}\truez{600}%
\put(-\value{x},\value{x}){\makebox(0,\value{z})[r]{${#1}$}}%
\put(\value{x},-\value{z}){\makebox(0,\value{z})[l]{${#2}$}}%
\end{picture}}%

\def\basicssswar[#1]{\SSSWAR{}{}{#100}}%

\newcommand{\ssswar}{\@ifnextchar[{\basicssswar}{\basicssswar[71]}}%

\def\basicSsswar[#1]#2{\SSSWAR{#2}{}{#100}}%

\newcommand{\Ssswar}{\@ifnextchar[{\basicSsswar}{\basicSsswar[71]}}%

\def\basicssswaR[#1]#2{\SSSWAR{}{#2}{#100}}%

\newcommand{\ssswaR}{\@ifnextchar[{\basicssswaR}{\basicssswaR[71]}}%

\newcommand{\SSSEAR}[3]{%
\Y=#3%
\divide \Y by 2%
\Z=\Y%
\multiply \Z by 3%
\begin{picture}(0,0)%
\put(-\Y,\Z){\line(1,-3){#3}}%
\put(\Y,-\Z){\sssehead}%
\truex{200}\truez{600}%
\put(\value{x},\value{x}){\makebox(0,\value{z})[l]{${#1}$}}%
\put(-\value{x},-\value{z}){\makebox(0,\value{z})[r]{${#2}$}}%
\end{picture}}%

\def\basicsssear[#1]{\SSSEAR{}{}{#100}}%

\newcommand{\sssear}{\@ifnextchar[{\basicsssear}{\basicsssear[71]}}%

\def\basicSssear[#1]#2{\SSSEAR{#2}{}{#100}}%

\newcommand{\Sssear}{\@ifnextchar[{\basicSssear}{\basicSssear[71]}}%

\def\basicssseaR[#1]#2{\SSSEAR{}{#2}{#100}}%

\newcommand{\ssseaR}{\@ifnextchar[{\basicssseaR}{\basicssseaR[71]}}%

\newcommand{\NNNWAR}[3]{%
\Y=#3%
\divide \Y by 2%
\Z=\Y%
\multiply \Z by 3%
\begin{picture}(0,0)%
\put(\Y,-\Z){\line(-1,3){#3}}%
\put(-\Y,\Z){\nnnwhead}%
\truex{200}\truez{600}%
\put(\value{x},\value{x}){\makebox(0,\value{z})[l]{${#1}$}}%
\put(-\value{x},-\value{z}){\makebox(0,\value{z})[r]{${#2}$}}%
\end{picture}}%

\def\basicnnnwar[#1]{\NNNWAR{}{}{#100}}%

\newcommand{\nnnwar}{\@ifnextchar[{\basicnnnwar}{\basicnnnwar[71]}}%

\def\basicNnnwar[#1]#2{\NNNWAR{#2}{}{#100}}%

\newcommand{\Nnnwar}{\@ifnextchar[{\basicNnnwar}{\basicNnnwar[71]}}%

\def\basicnnnwaR[#1]#2{\NNNWAR{}{#2}{#100}}%

\newcommand{\nnnwaR}{\@ifnextchar[{\basicnnnwaR}{\basicnnnwaR[71]}}%

\newcommand{\NEENEAR}[3]{%
\Y=#3%
\divide \Y by 2%
\Z=#3%
\divide \Z by 3%
\begin{picture}(0,0)%
\put(-\Y,-\Z){\line(3,2){#3}}%
\put(\Y,\Z){\neenehead}%
\truex{200}\truey{800}\truez{600}%
\put(-\value{x},\value{x}){\makebox(0,\value{z})[r]{${#1}$}}%
\put(\value{x},-\value{y}){\makebox(0,\value{z})[l]{${#2}$}}%
\end{picture}}%

\def\basicneenear[#1]{\NEENEAR{}{}{#100}}%

\newcommand{\neenear}{\@ifnextchar[{\basicneenear}{\basicneenear[215]}}%

\def\basicNeenear[#1]#2{\NEENEAR{#2}{}{#100}}%

\newcommand{\Neenear}{\@ifnextchar[{\basicNeenear}{\basicNeenear[215]}}%

\def\basicneeneaR[#1]#2{\NEENEAR{}{#2}{#100}}%

\newcommand{\neeneaR}{\@ifnextchar[{\basicneeneaR}{\basicneeneaR[215]}}%

\newcommand{\SEESEAR}[3]{%
\Y=#3%
\divide \Y by 2%
\Z=#3%
\divide \Z by 3%
\begin{picture}(0,0)%
\put(-\Y,\Z){\line(3,-2){#3}}%
\put(\Y,-\Z){\seesehead}%
\truex{200}\truey{800}\truez{600}%
\put(\value{x},\value{x}){\makebox(0,\value{z})[l]{${#1}$}}%
\put(-\value{x},-\value{y}){\makebox(0,\value{z})[r]{${#2}$}}%
\end{picture}}%

\def\basicseesear[#1]{\SEESEAR{}{}{#100}}%

\newcommand{\seesear}{\@ifnextchar[{\basicseesear}{\basicseesear[215]}}%

\def\basicSeesear[#1]#2{\SEESEAR{#2}{}{#100}}%

\newcommand{\Seesear}{\@ifnextchar[{\basicSeesear}{\basicSeesear[215]}}%

\def\basicseeseaR[#1]#2{\SEESEAR{}{#2}{#100}}%

\newcommand{\seeseaR}{\@ifnextchar[{\basicseeseaR}{\basicseeseaR[215]}}%

\newcommand{\NWWNWAR}[3]{%
\Y=#3%
\divide \Y by 2%
\Z=#3%
\divide \Z by 3%
\begin{picture}(0,0)%
\put(\Y,-\Z){\line(-3,2){#3}}%
\put(-\Y,\Z){\nwwnwhead}%
\truex{200}\truey{800}\truez{600}%
\put(\value{x},\value{x}){\makebox(0,\value{z})[l]{${#1}$}}%
\put(-\value{x},-\value{y}){\makebox(0,\value{z})[r]{${#2}$}}%
\end{picture}}%

\def\basicnwwnwar[#1]{\NWWNWAR{}{}{#100}}%

\newcommand{\nwwnwar}{\@ifnextchar[{\basicnwwnwar}{\basicnwwnwar[215]}}%

\def\basicNwwnwar[#1]#2{\NWWNWAR{#2}{}{#100}}%

\newcommand{\Nwwnwar}{\@ifnextchar[{\basicNwwnwar}{\basicNwwnwar[215]}}%

\def\basicnwwnwaR[#1]#2{\NWWNWAR{}{#2}{#100}}%

\newcommand{\nwwnwaR}{\@ifnextchar[{\basicnwwnwaR}{\basicnwwnwaR[215]}}%

\newcommand{\SWWSWAR}[3]{%
\Y=#3%
\divide \Y by 2%
\Z=#3%
\divide \Z by 3%
\begin{picture}(0,0)%
\put(\Y,\Z){\line(-3,-2){#3}}%
\put(-\Y,-\Z){\swwswhead}%
\truex{200}\truey{800}\truez{600}%
\put(-\value{x},\value{x}){\makebox(0,\value{z})[r]{${#1}$}}%
\put(\value{x},-\value{y}){\makebox(0,\value{z})[l]{${#2}$}}%
\end{picture}}%

\def\basicswwswar[#1]{\SWWSWAR{}{}{#100}}%

\newcommand{\swwswar}{\@ifnextchar[{\basicswwswar}{\basicswwswar[215]}}%

\def\basicSwwswar[#1]#2{\SWWSWAR{#2}{}{#100}}%

\newcommand{\Swwswar}{\@ifnextchar[{\basicSwwswar}{\basicSwwswar[215]}}%

\def\basicswwswaR[#1]#2{\SWWSWAR{}{#2}{#100}}%

\newcommand{\swwswaR}{\@ifnextchar[{\basicswwswaR}{\basicswwswaR[215]}}%

\newcommand{\NENNEAR}[3]{%
\Y=#3%
\divide \Y by 2%
\Z=#3%
\multiply \Z by 3%
\divide \Z by 4%
\begin{picture}(0,0)%
\put(-\Y,-\Z){\line(2,3){#3}}%
\put(\Y,\Z){\nennehead}%
\truex{100}\truez{600}%
\put(-\value{x},\value{x}){\makebox(0,\value{z})[r]{${#1}$}}%
\put(\value{x},-\value{z}){\makebox(0,\value{z})[l]{${#2}$}}%
\end{picture}}%

\def\basicnennear[#1]{\NENNEAR{}{}{#100}}%

\newcommand{\nennear}{\@ifnextchar[{\basicnennear}{\basicnennear[143]}}%

\def\basicNennear[#1]#2{\NENNEAR{#2}{}{#100}}%

\newcommand{\Nennear}{\@ifnextchar[{\basicNennear}{\basicNennear[143]}}%

\def\basicnenneaR[#1]#2{\NENNEAR{}{#2}{#100}}%

\newcommand{\nenneaR}{\@ifnextchar[{\basicnenneaR}{\basicnenneaR[143]}}%

\newcommand{\SWSSWAR}[3]{%
\Y=#3%
\divide \Y by 2%
\Z=#3%
\multiply \Z by 3%
\divide \Z by 4%
\begin{picture}(0,0)%
\put(\Y,\Z){\line(-2,-3){#3}}%
\put(-\Y,-\Z){\swsswhead}%
\truex{100}\truez{600}%
\put(-\value{x},\value{x}){\makebox(0,\value{z})[r]{${#1}$}}%
\put(\value{x},-\value{z}){\makebox(0,\value{z})[l]{${#2}$}}%
\end{picture}}%

\def\basicswsswar[#1]{\SWSSWAR{}{}{#100}}%

\newcommand{\swsswar}{\@ifnextchar[{\basicswsswar}{\basicswsswar[143]}}%

\def\basicSwsswar[#1]#2{\SWSSWAR{#2}{}{#100}}%

\newcommand{\Swsswar}{\@ifnextchar[{\basicSwsswar}{\basicSwsswar[143]}}%

\def\basicswsswaR[#1]#2{\SWSSWAR{}{#2}{#100}}%

\newcommand{\swsswaR}{\@ifnextchar[{\basicswsswaR}{\basicswsswaR[143]}}%

\newcommand{\SESSEAR}[3]{%
\Y=#3%
\divide \Y by 2%
\Z=#3%
\multiply \Z by 3%
\divide \Z by 4%
\begin{picture}(0,0)%
\put(-\Y,\Z){\line(2,-3){#3}}%
\put(\Y,-\Z){\sessehead}%
\truex{200}\truez{600}%
\put(\value{x},\value{x}){\makebox(0,\value{z})[l]{${#1}$}}%
\put(-\value{x},-\value{z}){\makebox(0,\value{z})[r]{${#2}$}}%
\end{picture}}%

\def\basicsessear[#1]{\SESSEAR{}{}{#100}}%

\newcommand{\sessear}{\@ifnextchar[{\basicsessear}{\basicsessear[143]}}%

\def\basicSessear[#1]#2{\SESSEAR{#2}{}{#100}}%

\newcommand{\Sessear}{\@ifnextchar[{\basicSessear}{\basicSessear[143]}}%

\def\basicsesseaR[#1]#2{\SESSEAR{}{#2}{#100}}%

\newcommand{\sesseaR}{\@ifnextchar[{\basicsesseaR}{\basicsesseaR[143]}}%

\newcommand{\NWNNWAR}[3]{%
\Y=#3%
\divide \Y by 2%
\Z=#3%
\multiply \Z by 3%
\divide \Z by 4%
\begin{picture}(0,0)%
\put(\Y,-\Z){\line(-2,3){#3}}%
\put(-\Y,\Z){\nwnnwhead}%
\truex{200}\truez{600}%
\put(\value{x},\value{x}){\makebox(0,\value{z})[l]{${#1}$}}%
\put(-\value{x},-\value{z}){\makebox(0,\value{z})[r]{${#2}$}}%
\end{picture}}%

\def\basicnwnnwar[#1]{\NWNNWAR{}{}{#100}}%

\newcommand{\nwnnwar}{\@ifnextchar[{\basicnwnnwar}{\basicnwnnwar[143]}}%

\def\basicNwnnwar[#1]#2{\NWNNWAR{#2}{}{#100}}%

\newcommand{\Nwnnwar}{\@ifnextchar[{\basicNwnnwar}{\basicNwnnwar[143]}}%

\def\basicnwnnwaR[#1]#2{\NWNNWAR{}{#2}{#100}}%

\newcommand{\nwnnwaR}{\@ifnextchar[{\basicnwnnwaR}{\basicnwnnwaR[143]}}%

\newcommand{\Necurve}[2]%
{\begin{picture}(0,0)%
\truex{1300}\truey{2000}\truez{200}%
\put(0,\value{x}){\oval(#200,\value{y})[t]}%
\put(0,\value{x}){\makebox(0,0){\begin{picture}(#200,0)%
\put(#200,0){\line(0,-1){\value{z}}}%
\put(#200,-\value{z}){\shead}%
\put(0,0){\line(0,-1){\value{z}}}\end{picture}}}%
\truex{2500}%
\put(0,\value{x}){\makebox(0,0)[b]{${#1}$}}%
\end{picture}}%

\def\basicnecurvar[#1]{\Necurve{}{#1}}

\newcommand{\necurvar}{\@ifnextchar[{\basicnecurvar}{\basicnecurvar[160]}}%

\def\basicNecurvar[#1]#2{\Necurve{#2}{#1}}%

\newcommand{\Necurvar}{\@ifnextchar[{\basicNecurvar}{\basicNecurvar[160]}}%

\newcommand{\Nwcurve}[2]%
{\begin{picture}(0,0)%
\truex{1300}\truey{2000}\truez{200}%
\put(0,\value{x}){\oval(#200,\value{y})[t]}%
\put(0,\value{x}){\makebox(0,0){\begin{picture}(#200,0)%
\put(#200,0){\line(0,-1){\value{z}}}%
\put(0,0){\line(0,-1){\value{z}}}%
\put(0,-\value{z}){\shead}%
\end{picture}}}%
\truex{2500}%
\put(0,\value{x}){\makebox(0,0)[b]{${#1}$}}%
\end{picture}}%

\def\basicnwcurvar[#1]{\Nwcurve{}{#1}}

\newcommand{\nwcurvar}{\@ifnextchar[{\basicnwcurvar}{\basicnwcurvar[160]}}%

\def\basicNwcurvar[#1]#2{\Nwcurve{#2}{#1}}%

\newcommand{\Nwcurvar}{\@ifnextchar[{\basicNwcurvar}{\basicNwcurvar[160]}}%

\newcommand{\Securve}[2]%
{\begin{picture}(0,0)%
\truex{1300}\truey{2000}\truez{200}%
\put(0,-\value{x}){\oval(#200,\value{y})[b]}%
\put(0,-\value{x}){\makebox(0,0){\begin{picture}(#200,0)%
\put(#200,0){\line(0,1){\value{z}}}%
\put(0,0){\line(0,1){\value{z}}}%
\put(#200,\value{z}){\nhead}%
\end{picture}}}%
\truex{2500}%
\put(0,-\value{x}){\makebox(0,0)[t]{${#1}$}}%
\end{picture}}%

\def\basicsecurvar[#1]{\Securve{}{#1}}

\newcommand{\securvar}{\@ifnextchar[{\basicsecurvar}{\basicsecurvar[160]}}%

\def\basicSecurvar[#1]#2{\Securve{#2}{#1}}%

\newcommand{\Securvar}{\@ifnextchar[{\basicSecurvar}{\basicSecurvar[160]}}%

\newcommand{\Swcurve}[2]%
{\begin{picture}(0,0)%
\truex{1300}\truey{2000}\truez{200}%
\put(0,-\value{x}){\oval(#200,\value{y})[b]}%
\put(0,-\value{x}){\makebox(0,0){\begin{picture}(#200,0)%
\put(#200,0){\line(0,1){\value{z}}}%
\put(0,0){\line(0,1){\value{z}}}%
\put(0,\value{z}){\nhead}%
\end{picture}}}%
\truex{2500}%
\put(0,-\value{x}){\makebox(0,0)[t]{${#1}$}}%
\end{picture}}%

\def\basicswcurvar[#1]{\Swcurve{}{#1}}

\newcommand{\swcurvar}{\@ifnextchar[{\basicswcurvar}{\basicswcurvar[160]}}%

\def\basicSwcurvar[#1]#2{\Swcurve{#2}{#1}}%

\newcommand{\Swcurvar}{\@ifnextchar[{\basicSwcurvar}{\basicSwcurvar[160]}}%

\newcommand{\Escurve}[2]%
{\begin{picture}(0,0)%
\truex{1400}\truey{2000}\truez{200}%
\put(\value{x},0){\oval(\value{y},#200)[r]}%
\put(\value{x},0){\makebox(0,0){\begin{picture}(0,#200)%
\put(0,0){\line(-1,0){\value{z}}}%
\put(0,#200){\line(-1,0){\value{z}}}%
\put(-\value{z},0){\whead}%
\end{picture}}}%
\truex{2500}%
\put(\value{x},0){\makebox(0,0)[l]{${#1}$}}%
\end{picture}}%

\def\basicescurvar[#1]{\Escurve{}{#1}}

\newcommand{\escurvar}{\@ifnextchar[{\basicescurvar}{\basicescurvar[160]}}%

\def\basicEscurvar[#1]#2{\Escurve{#2}{#1}}%

\newcommand{\Escurvar}{\@ifnextchar[{\basicEscurvar}{\basicEscurvar[160]}}%

\newcommand{\Encurve}[2]%
{\begin{picture}(0,0)%
\truex{1400}\truey{2000}\truez{200}%
\put(\value{x},0){\oval(\value{y},#200)[r]}%
\put(\value{x},0){\makebox(0,0){\begin{picture}(0,#200)%
\put(0,0){\line(-1,0){\value{z}}}%
\put(0,#200){\line(-1,0){\value{z}}}%
\put(-\value{z},#200){\whead}%
\end{picture}}}%
\truex{2500}%
\put(\value{x},0){\makebox(0,0)[l]{${#1}$}}%
\end{picture}}%

\def\basicencurvar[#1]{\Encurve{}{#1}}

\newcommand{\encurvar}{\@ifnextchar[{\basicencurvar}{\basicencurvar[160]}}%

\def\basicEncurvar[#1]#2{\Encurve{#2}{#1}}%

\newcommand{\Encurvar}{\@ifnextchar[{\basicEncurvar}{\basicEncurvar[160]}}%

\newcommand{\Wscurve}[2]%
{\begin{picture}(0,0)%
\truex{1300}\truey{2000}\truez{200}%
\put(-\value{x},0){\oval(\value{y},#200)[l]}%
\put(-\value{x},0){\makebox(0,0){\begin{picture}(0,#200)%
\put(0,0){\line(1,0){\value{z}}}%
\put(0,#200){\line(1,0){\value{z}}}%
\put(\value{z},0){\ehead}%
\end{picture}}}%
\truex{2400}%
\put(-\value{x},0){\makebox(0,0)[r]{${#1}$}}%
\end{picture}}%

\def\basicwscurvar[#1]{\Wscurve{}{#1}}

\newcommand{\wscurvar}{\@ifnextchar[{\basicwscurvar}{\basicwscurvar[160]}}%

\def\basicWscurvar[#1]#2{\Wscurve{#2}{#1}}%

\newcommand{\Wscurvar}{\@ifnextchar[{\basicWscurvar}{\basicWscurvar[160]}}%

\newcommand{\Wncurve}[2]%
{\begin{picture}(0,0)%
\truex{1300}\truey{2000}\truez{200}%
\put(-\value{x},0){\oval(\value{y},#200)[l]}%
\put(-\value{x},0){\makebox(0,0){\begin{picture}(0,#200)%
\put(0,0){\line(1,0){\value{z}}}%
\put(\value{z},#200){\ehead}%
\put(0,#200){\line(1,0){\value{z}}}%
\end{picture}}}%
\truex{2400}%
\put(-\value{x},0){\makebox(0,0)[r]{${#1}$}}%
\end{picture}}%

\def\basicwncurvar[#1]{\Wncurve{}{#1}}

\newcommand{\wncurvar}{\@ifnextchar[{\basicwncurvar}{\basicwncurvar[160]}}%

\def\basicWncurvar[#1]#2{\Wncurve{#2}{#1}}%

\newcommand{\Wncurvar}{\@ifnextchar[{\basicWncurvar}{\basicWncurvar[160]}}%

\newcommand{\crosslength}[2]{%
\settowidth{\firstitem}{#1}%
\settowidth{\seconditem}{#2}%
\ifdim\firstitem < \seconditem%
\itemlength=\seconditem%
\else%
\itemlength=\firstitem%
\fi%
\divide\itemlength by 2%
\hspace{\itemlength}}%

\catcode`\@=12

\newcommand{\PreOrdC}{\mathsf{PreOrd}(\calC)}

\newcommand{\id}{\mathsf{id}}
\newcommand{\ONE}{\mathbf{1}}
\newcommand{\Cat}{\mathsf{Cat}}
\newcommand{\Isom}{\mathsf{Iso}}
\newcommand{\Aut}{\mathsf{Aut}}
\newcommand{\PreOrd}{\mathsf{PreOrd}}
\newcommand{\Ob}{\mathsf{Ob}}
\newcommand{\tr}{{\triangleright}}
\newcommand{\ArB}{\mathsf{Ar}(\calB)}
\newcommand{\Arrows}{\mathsf{Ar}}

\newcommand{\proof}{\medbreak\noindent%
{\it Proof}\hspace{1em plus .5em}}

\newcommand{\qed}{\mbox{}\hfill\mbox{}%
\mbox{$\square$}\medbreak}

\newcommand{\calC}{{\mathcal C}}
\newcommand{\calT}{{\mathcal T}}
\newcommand{\calF}{{\mathcal F}}

\newcommand{\calA}{{\mathcal A}}
\newcommand{\calZ}{{\mathcal Z}}
\newcommand{\calS}{{\mathcal S}}
\newcommand{\calQ}{{\mathcal Q}}
\newcommand{\calX}{{\mathcal X}}
\newcommand{\calG}{{\mathcal G}}
\newcommand{\calY}{{\mathcal Y}}
\newcommand{\calB}{{\mathcal B}}
\newcommand{\calK}{{\mathcal K}}

\newtheorem{theorem}{Theorem}[section]
\newtheorem{definition}[theorem]{Definition}
\newtheorem{proposition}[theorem]{Proposition}
\newtheorem{lemma}[theorem]{Lemma}
\newtheorem{remark}[theorem]{Remark}
\newtheorem{corollary}[theorem]{Corollary}

\title{Groupoids and skeletal categories\\form a pretorsion theory in $\Cat$}
\author{Francis Borceux, 
Federico Campanini\footnote{Partially supported by Ministero dell’Istruzione, dell'Università e della Ricerca (Progetto di ricerca di rilevante interesse nazionale “Categories, Algebras: Ring-Theoretical and Homological Approaches (CARTHA)).}
\\Marino Gran\footnote{This work was supported by the \emph{Fonds de la Recherche Scientifique - FNRS} under Grant  CDR
no. J.0080.23}, 
Walter Tholen\footnote{The fourth author acknowledges partial financial assistance by Natural Sciences and Engineering Council of Canada under the Discovery Grants Program, no. 501260.}}
\date{Universit\'e Catholique de Louvain\\York University in Toronto}
\begin{document}

\maketitle

\begin{abstract}
We describe a pretorsion theory  in the category $\Cat$ of small categories: the torsion objects are the groupoids, while the torsion-free objects are the skeletal categories,  {\em i.e.}, those categories in which every isomorphism is an automorphism. We infer these results from two unexpected properties of coequalizers in $\Cat$ that identify pairs of objects: they are faithful and reflect isomorphisms.
\end{abstract}

\section{Introduction}
A {\em pretorsion theory}  in a category $\calC$ consists of two classes of objects, respectively called the {\em torsion} and {\em torsion-free} objects, together with axioms which mirror those for a {\em torsion theory} in an abelian category. The objects which are both {\em torsion} and {\em torsion-free} are called {\em trivial}.  When $\calC$ has a zero object which is taken to be the only trivial object, we recover the classical notion. Hence, in direct generalization of its abelian origins, the current setting is applicable in an arbitrary category, even in the absence of a zero object. 

The original example of a torsion theory is, of course, that in the category of abelian groups, with the usual notions of torsion group and torsion-free group. Tens of papers have been devoted to torsion theories in various non-abelian contexts, including \cite{BG, CDT, RT, GJ, GKV, GC, PVdL, D, HKvR}. For a first easy example of a pretorsion theory, consider the category  $\PreOrd$ of preordered sets. Choosing the equivalence relations as torsion objects and the partial orders as torsion-free objects, one obtains a pretorsion theory on $\PreOrd$ (see \cite{fafi}) with the discrete objects as trivial ones. This fact has been generalized to the category $\PreOrdC$ of preordered objects in a Barr-exact category $\calC$ (see \cite{FFG2, bocagr, bocagr2}). Other examples of pretorsion theories have been studied in \cite{FFG, bocagr3, FHZ, xa, GM}.

A preordered set may be seen as a small category with at most one arrow between any two objects. On the other hand, given a small category $\calC$, one gets a preorder on its set of objects by declaring $A\leq B$ when there is at least one morphism from $A$ to $B$. In \cite{xa}, a small category is called  ``torsion'' or ``torsion free'' when it is the case for the corresponding preordered set of its objects; it is proved that this yields a pretorsion theory on $\Cat$.

In the present paper, we adopt a totally different approach, based on the structure of the category of arrows, not just on the existence of arrows. In our pretorsion theory on $\Cat$, a small category $\calC$ is
\begin{description}
\item[torsion] when $\calC$ is a groupoid, {\em i.e.}, when every morphism in $\calC$ is an isomorphism;
\item[torsion-free] when $\calC$ is a skeletal category, {\em i.e.}, when every isomorphism in $\calC$ is an automorphism {(see \cite{macl})}.
\end{description}
Of course, when the category $\calC$ is just a preordered set, we recapture the situation studied in \cite{fafi}. This was also the case for the pretorsion theory studied in \cite{xa}.

The results in the present paper rely heavily on a careful study of some particular coequalizers in $\Cat$: the coequalizers of pairs of functors defined on a discrete category. Such coequalizers have properties that are atypical for a quotient functor, the most striking ones being their faithfulness and the reflection of isomorphisms. Establishing these properties is quite cumbersome and reduces to so-called  ``word problems'' on formal chains of arrows in the  quotient graph used to construct the coequalizer. The work of  John Isbell on {\em Dominions} (1968, see \cite{isbell}) and the {\em Diplomarbeit} of Reinhard B\"orger (1977, see \cite{borger}) provided us with the necessary tools for handling these problems. 

We conclude the paper with showing the existence of the so-called {\it $\calZ$-kernels} and {\it $\calZ$-cokernels} in $\Cat$, that is,  kernels and cokernels defined relatively to the ideal of trivial morphisms: those morphisms factoring through a trivial object.  Thanks to the results in \cite{bocagr3} it then follows that there is a \emph{stable category} naturally associated with this pretorsion theory, that satisfies an interesting universal property.

\section{Pretorsion theories}

Let us recall that an ideal $\calZ$ in a category $\calC$ is a class of arrows such that
for every arrow $f\in\calZ$, one has $fu\in\calZ$ and 
$vf\in\calZ$, for all arrows $u$, $v$ composable with $f$ \cite{Ehresmann}. When $\calC$ has a zero object, the zero morphisms constitute an ideal.

Given an ideal $\calZ$ in a category $\calC$, an arrow $k$ is the $\calZ$-kernel of an arrow $f$ when $fk\in\calZ$ and, if $fm\in\calZ$ for some arrow $m$, then $m$ factors uniquely through $k$. The uniqueness condition forces $k$ to be a monomorphism. When $\calZ$ is the ideal of zero morphisms, we recapture the usual notion of kernel. There is of course a dual notion of $\calZ$-cokernel. A pair of composable morphisms
$$
K\Ar k A \Ar q Q
$$
is a short $\calZ$-exact sequence when $k$ is the $\calZ$-kernel of $q$ and $q$ is the $\calZ$-cokernel of $k$.

 The following definition was introduced in \cite{fafi} and then thoroughly investigated in \cite{FFG}:
\begin{definition}\label{pretor}
A pretorsion theory in a category $\calC$ consists of a pair $(\calT,\calF)$ of classes of objects, both of them closed under isomorphisms, whose elements are the {\em torsion} and the {\em torsion-free} objects of the pretorsion theory, respectively. The objects in $\calT\cap\calF$ are called {\em trivial}, and the ideal $\calZ$ of {\em trivial morphisms} is that of those arrows factoring through a trivial object.\\
These data must satisfy the following two axioms:
\begin{description}
\item[PT1] every arrow $f\colon A\ar B$ with $A\in\calT$ and $B\in\calF$ is trivial;
\item[PT2] for every object $A\in\calC$, there exists a short $\calZ$-exact sequence
$$
K\Ar k A \Ar q Q
$$
with $K\in\calT$ and $Q\in\calF$.
\end{description}
\end{definition}

\section{Some coequalizers in $\Cat$ revisited} \label{Section_3}

 Let us first exhibit a result borrowed from J.\ Isbell (see \cite{isbell}) and  R.~B\"orger (see \cite{borger}). 

\begin{definition}[Def. 4.3 in \cite{borger}]\label{freem}
Let $\calB$ be a small category. Let $M$ be the free monoid on the set of arrows of $\calB$. An element of $M$ -- that is, an arbitrary finite sequence of arrows in $\calB$ -- is {\em reduced} when
\begin{itemize}
\item the sequence does not contain any pair of consecutive arrows which turn out to be  composable in $\calB$;
\item the sequence does not contain any identity morphism of $\calB$.
\end{itemize}
\end{definition}

Let us clarify that, when speaking of a {\em a pair of consecutive composable  arrows} in a sequence, we always refer to a sequence $(f_1,\ldots,f_n)$  and an index $i<n$ such that the codomain of $f_i$ is equal to the domain of $f_{i+1}$. To avoid any confusion with the (standard) direction in which we write arrows,
for clarity we shall sometimes write the composite $f_{i+1}f_i$ in the reverse order, as $f_i\tr f_{i+1}$, so that the shorter sequence in which the pair $(f_i,f_{i+1})$ is replaced by the composite of the two arrows, may then be written as
$$(f_1, \cdots,f_{i-1},f_i\tr f_{i+1}, f_{i+2}, \cdots,f_n).
$$

The unit element of $M$ is of course the empty sequence, which we denote by $()$. By a congruence on the monoid $M$ is meant an equivalence relation on $M$ such that $x\sim y$ implies $xz\sim yz$ and $zx\sim zy$ for all $z$.

\begin{proposition}[Statement 1.1 in \cite{isbell}, Satz 4.4 in \cite{borger}]\label{reduced}~\\
Under the conditions of Definition~\ref{freem}, consider the smallest congruence $S$ on $M$ such that 
\begin{itemize}
\item $\bigl((u,v),u\tr v\bigr)\in S$ for every pair $(u,v)$ of composable morphisms in $\calB$;
\item $\bigl(\id_B,()\bigr)\in S$ for every object $B\in\calB$.
\end{itemize}
Then every element of $M$ is $S$-equivalent to a unique reduced element.
\end{proposition}

\proof
In 1968, after referring to the work of Mersch (see \cite{mersch}), Isbell states as item 1.1 of his paper \cite{isbell} a result which is essentially our Proposition~\ref{reduced}; however, the proof is only sketched: it makes explicit the representation in terms of reduced chains, leaving off any details of a proof. Isbell simply says that these are analogous to those for free groups (see \cite{artin} and \cite{vdw}). 

In {\em Satz 4.4} of his 1977 {\em Diplomarbeit} \cite{borger}, B\"orger presents a long and carefully written proof, with all technical details. Since this work has never been published, we took the liberty of making a scan of it available, together with a detailed sketch of his proof: see the appendix to this paper.
\qed

Let us now switch to coequalizers in $\Cat$. Their construction has been described by many authors, with varying levels of detail: see for example Section~~5.1 in \cite{BOR}, Section 4 in \cite{BBP}, Section~11  in \cite{tho}, and so on. All of them use the notion of {\em congruence} on a category $\calG$, that is: an equivalence relation on each hom-set $\calG(X,Y)$ such that $f\sim g$ implies $fh\sim gh$ and $kf\sim kg$ for all arrows $h$, $k$ composable with $f$. 

Given a pair $(F,G)$ of functors between small categories
$$
\calA \Biar FG \calB \Ar Q \calQ
$$
their coequalizer $Q$ can thus be obtained in the following way.
\begin{itemize}
\item As far as objects are concerned, $\Ob(\calQ)$ is the quotient of $\Ob(\calB)$ by the equivalence relation generated by $F(A)\sim G(A)$, for every object $A\in\calA$.
\item First one constructs a graph $\calG_0$ on this quotient set of objects, by putting every arrow $f\colon A \ar B$ in $\calB$ as an arrow from $[A]$ to $[B]$ in $\calG_0$.
\item Next one considers the category $\calG$ having  the same objects as $\calG_0$; its arrows are the non-empty finite chains of ``composable'' morphisms of $\calG_0$ (in the sense  that the $\calG_0$-codomain of a $\calG_0$-morphism is the $\calG_0$-domain of the next one), together with an empty chain from each object of $\calG_0$ to itself;  the composition is just concatenation.
\item The coequalizer $\calQ$ is the quotient of $\calG$ by the congruence generated by:
\begin{description}
\item[Q1] when two consecutive morphisms in a $\calG_0$-chain are composable in $\calB$, the chain is equivalent to the one obtained when replacing the corresponding pair by its composite in $\calB$;
\item[Q2] when an identity morphism appears in a $\calG_0$-chain, the chain is equivalent to the one obtained when dropping that identity; 
\item[Q3] when a morphism of the form $F(f)$ appears in a $\calG_0$-chain, the chain is equivalent to the one obtained when replacing $F(f)$ by $G(f)$.
\end{description}
\end{itemize}
The following result is then an immediate consequence of the Isbell-B\"orger result presented above.

\begin{proposition}\label{red}
With the notation above, let us assume that the category $\calA$ is discrete. Given an arrow in the coequalizer $\calQ$, there then exists a unique reduced arrow in $\calG$ representing it.
\end{proposition}

\proof
A non-identity arrow  in $\calQ$ is represented by a triple $(A,\xi,B)$ where, with the notation of Definition \ref{freem}, $\xi\in M$ is a sequence of consecutive arrows in $\calG_0$, the first one having $A$ as domain in $\calB$, and the last one having $B$ as codomain in $\calB$. An identity arrow in $\calQ$ is represented by a triple $\bigl(B,(),B\bigr)$, where $()$ is the empty sequence and $B$ is an object of $\calB$.

Since $\calA$ is discrete, condition Q2 above indicates at once that condition Q3 can be omitted, because the only possible values for $f$ in Q3 are identities. Thus the congruence defining $\calQ$ from $\calG$  is constructed using only conditions Q1 and Q2, just as the congruence on $M$ in Proposition \ref{reduced}. Therefore, the result follows at once from Proposition~\ref{reduced}.
\qed

Let us now infer various interesting consequences from this last proposition.

\begin{proposition}\label{faith}
The coequalizer of two functors defined on a discrete category is faithful.
\end{proposition}

\proof
With the notation above,  let $f,g\colon A\biar B$ be two non-identity morphisms in $\calB$  such that $Q(f)=Q(g)$. Then $(f)$ and $(g)$ are reduced morphisms in $\calG$ representing the same arrow of $\calQ$, thus they are equal by Proposition \ref{red}.

Next, if $f=\id_B$ and $g\colon B\ar B$, with $Q(f)=Q(g)$, $g$ must be an identity, otherwise we would again have two reduced morphisms $()$ and $(g)$ of $\calG$ representing the same morphism of $\calQ$.
\qed

\begin{lemma}\label{isos}
Consider the coequalizer $\calQ$ of two functors 
$F,G\colon\calA  \biar \calB$   defined on a discrete category $\calA$.
With the notation above, a morphism in $\calQ$ is an isomorphism if, and only if, its reduced form in $\calG$ is empty or composed of isomorphisms in $\calB$.
\end{lemma}

\proof
By \ref{red}, consider the reduced form of a non-identity  isomorphism 
$$
( f_1, \ldots,f_n),~~f_i\colon A_i\ar B_i,~~Q(B_i)=Q(A_{i+1}),~~B_i\not=A_{i+1},
$$
in $\calQ$ and the reduced form of its inverse
$$
( g_1, \ldots,g_m),~~g_i\colon C_i\ar D_i,~~Q(D_i)=Q(C_{i+1}),~~D_i\not=C_{i+1};
$$
thus, in particular, $Q(A_1)=Q(D_m)$ and $Q(C_1)=Q(B_n)$. The sequence
$$
( f_1, \ldots, f_n, g_1, \ldots g_m )
$$
must therefore be equivalent to the empty sequence. If $B_n\not=C_1$, this last sequence is reduced and equivalent to the empty sequence, which is impossible by the uniqueness condition in Proposition \ref{red}. Thus $B_n=C_1$ and we can shorten the sequence to
$$
( f_1, \ldots, f_{n-1},  f_n\tr g_1, g_2, \ldots g_m ).
$$
But the domain of $ f_n\tr g_1$ is $A_n\not=B_{n-1}$, and the codomain of $f_n\tr g_1$ is $D_1\not=C_2$. So, if $f_n\tr g_1$ is not an identity, the shortened sequence is reduced, which contradicts again the fact that it is equivalent to the empty sequence. Thus $f_n\tr g_1$ is an identity in $\calA$ and, in particular, $A_n=D_1$.
Next, looking at the other composite
$$
( g_1,\cdots,g_m,f_1,\ldots,f_n )
$$
which must also be an identity morphism  in $\calQ$, we conclude that $g_1$ is a monomorphism in $\calQ$ and $f_n$ is an epimorphism in $\calQ$. But by  Proposition \ref{faith} the functor $Q$ is faithful, thus in particular it reflects monomorphisms and epimorphisms. Therefore $g_1$ is a monomorphism in $\calB$ and $f_n$ is an epimorphism in $\calB$. And since $f_n\tr g_1$ is an identity in $\calB$, $g_1$ is both a retraction and a monomorphism while $f_n$ is both a section and an epimorphism; they are thus inverse isomorphisms.

So, one can further shorten the situation and obtain the sequence
$$
( f_1, \ldots, f_{n-1}, g_2, \ldots g_m ).
$$
One repeats inductively the same process as above, up to the moment when we have used all the components of one of the two original reduced sequences. In this way we end up,  let us say,  with
$$
( f_1,\ldots,f_k)
$$
(the case of any remaining $g_i$s is, of course, analogous). But if not empty, this sequence is reduced, since the sequence $( f_1,\ldots,f_n)$ was reduced, and it must therefore be  equivalent to the empty sequence. This is impossible, again by the uniqueness condition in Proposition \ref{red}. Thus, we in fact ended up with an empty sequence.

This concludes the proof of one of the two stated implications; the other one is obvious.
\qed

\begin{proposition}\label{refl}
The coequalizer of two functors defined on a discrete category reflects isomorphisms.
\end{proposition}

\proof
With the same notation as above, let $f$ be a morphism in $\calB$ such that $Q( f)$ is an isomorphism. If $f$ is an identity, there is nothing to prove.  Otherwise, the morphism $( f)$ of $\calG$ is in reduced form and, by assumption, an isomorphism in $\calQ$. By  Lemma \ref{isos}, $f$ is an isomorphism in $\calB$.
\qed

\section{The pretorsion theory in $\Cat$}

The following result extends, to small categories, the pretorsion theory studied in \cite{fafi} in the case of preordered sets, {\em i.e.}, of those categories having at most one arrow between any two objects.
We shall denote by $\mathsf{Grpd}, \mathsf{SkCat}$ and $\mathsf{SkGrpd}$ the (full) subcategories of $\mathsf{Cat}$ whose objects are groupoids, skeletal categories and skeletal groupoids, respectively.

\begin{theorem}\label{graut}
The pair $(\mathsf{Grpd}, \mathsf{SkCat})$ is a pretorsion theory in $\mathsf{Cat}$. 
\end{theorem}

\proof
The trivial objects are thus those categories in which all arrows are automorphisms, i.e. the skeletal groupoids. We shall write $\Isom(\calC)$ for the groupoid of isomorphisms of a small category $\calC$ and  $\Aut(\calC)$ for the groupoid of its automorphisms. $\Aut(\calC)$ is thus a trivial category, in the sense of Definition~\ref{pretor}. So, a functor $F\colon \calA\ar\calB$ is trivial when it factors through $\Aut(\calB)$. In particular every functor $F\colon \calG\ar\calC$ from a groupoid $\calG$ to a torsion-free category $\calC$ is trivial. This takes care of axiom PT1.

Consider now an arbitrary small category $\calC$; it contains the groupoid $\Isom(\calC)$. We form the following coequalizer in $\Cat$:
$$
\coprod_{\sigma\mbox{~iso}}\ONE \Biar{d_0}{d_1}\calC\Ar Q \calQ
$$
where the left-hand category is a copower of the terminal category $\ONE$, indexed by the set of all isomorphisms of $\calC$. This is of course a discrete category. The functors $d_0$ and $d_1$ are those which, on the component indexed by an isomorphism $\sigma$, map the unique object of $\ONE$ respectively to the domain and codomain of $\sigma$. On one hand, all the isomorphisms of $\calC$ are thus mapped by $Q$ to automorphisms in $\calQ$. On the other hand, by Lemma \ref{isos}, an isomorphism in $\calQ$ is an identity or a composite of images of isomorphisms in $\calC$. It is therefore a composite of automorphisms in $\calQ$ and, thus, an automorphism in $\calQ$. Consequently, $\calQ$ is skeletal.

It remains to to be proved that we have  obtained a short $\calZ$-exact sequence $(i,Q)$:
\begin{diagram}[70]
\calX ||
\sdotar | \Sear G ||
\Isom(\calC) | \Ear i | \calC | \Eepi Q | \calQ ||
| | |\seaR H | \sdotar ||
| | | | \calY ||
\end{diagram}
Consider first a functor $G\colon \calX\ar\calC$ such that $QG$ is trivial. This means that, for every arrow $x\in\calX$, $G(x)$ is mapped by $Q$ to an automorphism and, thus, an isomorphism. But by Proposition \ref{refl}, every $G(x)$ is then an isomorphism in $\calC$ and, thus, $G$ factors through $\Isom(\calC)$. That factorization is unique since $i$ is an inclusion functor. On the other hand, if $H\colon \calC \ar \calY$ is such that $Hi$ is trivial, then every isomorphism of $\calC$ is mapped by $H$ to an automorphism, thus its domain and its codomain are identified by $H$. This proves that $Hd_0=Hd_1$, and we get the expected unique factorization through the coequalizer $Q$ of $(d_0,d_1)$.
\qed

To the best of our knowledge, the following consequence of Theorem \ref{graut} has not yet been stated in the literature.

\begin{corollary}\label{skeletal}
The full subcategory $\mathsf{SkCat}$ of skeletal categories is reflective in the category $\mathsf{Cat}$ of small categories.
\end{corollary}

\proof
In a category $\calC$ provided with a pretorsion  theory, the full subcategory of torsion objects is  coreflective and that of torsion-free objects is reflective (see \cite{FFG}).
\qed

\begin{remark}\label{skeleton}
Comparison with the notion of skeleton.
\end{remark}

 Applying a strong (potentially class-based) version of the axiom of choice one easily sees that every category $\calC$ is equivalent to a skeletal full subcategory $\calS$ of $\calC$ (see \cite{pareig}), called its {\em skeleton} (see \cite{macl}).  It is important to note that, in the notation of the proof of Theorem \ref{graut}, such a skeleton $\calS$ (uniquely determined only up to equivalence) is by no means the $\calZ$-cokernel $\calQ$ and, thus, the skeletal reflection of $\calC$, as given (without the use of any choice principle) by Corollary \ref{skeletal}. For example, consider the category $\calC$ with two objects, the identities on these, and a unique isomorphism $(f,f^{-1})$ between the two objects. Both candidates for a skeleton of $\calC$ are isomorphic to the terminal category $\ONE$. However, in the $\calZ$-cokernel $\calQ$ as in Theorem \ref{graut}   with its single object and its identity morphism, there are also all the powers $f^n$, $(f^{-1})^n$. It is thus the monoid $({\Bbb Z},+)$, viewed as a one-object category.
\qed

\begin{remark}\label{comp} 
An ``internalization'' of Theorem \ref{graut}.
 \end{remark}

The construction of the pretorsion theory on preordered sets (see \cite{fafi}) may easily be carried over to the case of preordered objects in a Barr-exact category  (see \cite{FFG2}), because it refers only to finite limits and coequalizers of kernel pairs  (see Barr's metatheorem  \cite{barr}). The case of the pretorsion  theory in Theorem \ref{graut} is strikingly different, because the arguments that we have developed -- and in particular B\"orger's result in the Appendix -- are highly set theoretical and cannot be carried out as such in a quite arbitrary category. 
In most classical algebraic categories (such as the categories of groups, Lie algebras, rings, modules, etc.) the internal categories are always groupoids. This is actually the case in any Mal'tsev variety \cite{Smith}, thus internalizing Theorem \ref{graut} to this context is equivalent to proving that the category $\mathsf{SkGrpd}({\mathcal C})$ of internal skeletal groupoids -- {\em  i.e.} those whose domain and codomain morphisms are equal -- is epireflective in the category $\mathsf{Grpd}({\mathcal C})$ of internal groupoids. Now, when $\mathcal C$ is a Mal'tsev variety, the category $\mathsf{Grpd}({\mathcal C})$ is again a Mal'tsev variety, since it is a subvariety of the variety of reflexive graphs in $\mathcal C$  (see Corollary $2.4$ in \cite{GR}). It is then easy to see that the subcategory $\mathsf{SkGrpd}({\mathcal C})$ of skeletal groupoids in $\mathcal C$ is a subvariety of $\mathsf{Grpd}({\mathcal C})$. Indeed, the subcategory $\mathsf{SkGrpd}({\mathcal C})$ is determined by the additional identity  expressing the fact that the (unary) operations induced by the domain and codomain morphisms have to be equal. By the Birkhoff theorem $\mathsf{SkGrpd}({\mathcal C})$ is then a subvariety and, in particular, it is epireflective in $\mathsf{Grpd}({\mathcal C})$.

\section{The existence of $\calZ$-kernels and $\calZ$-cokernels} 

Let us now prove the existence of all $\calZ$-kernels and $\calZ$-cokernels in $\Cat$, with respect to  the ideal $\calZ$ of trivial morphisms determined by the subcategory $ \mathsf{SkGrpd}$ of skeletal groupoids, as in Theorem~\ref{graut}.

\begin{proposition}\label{zker}
Let $F\colon \calA \to \calB$ be a functor in $\Cat$. Its $\calZ$-kernel is given by the functor $k$ in the following pullback
\begin{diagram}[80]
\calK | \ear | \Aut(\calB) ||
\Smono k | | \smono ||
\calA| \eaR F | \calB ||
\end{diagram}
\end{proposition}

\proof
Given $G\colon \calX\ar \calA$, the functor $FG$ is trivial precisely when it factors through $\Aut(\calB)$, thus through the pullback.
\qed

The case of $\calZ$-cokernels is more involved. Inspired by considerations in \cite{GJ}, we prove first:

\begin{lemma}\label{zcokerid}
The $\calZ$-cokernel of the identity functor on a small category $\calA$ exists and can be constructed in the following way:
\begin{itemize}
\item consider first the category of fractions inverting all the arrows of $\calA$
$$
p\colon \calA \ar \calG=\calA\bigl[\Arrows(\calA)^{-1}\bigr];
$$
\item consider next the canonical short $\calZ$-exact sequence of $\calG$:
$$
\calK\Ar k \calG \Ar q \calQ.
$$
\end{itemize}
Then $k$ is isomorphic to the identity on $\calG$ while the composite 
$$
\calA \Ar p \calG \Ar q \calQ
$$
is the $\calZ$-cokernel of $F$.
\end{lemma}

\proof
The construction of the $\calZ$-exact sequence in Theorem~\ref{graut} shows at once that the $\calZ$-kernel part of the sequence is the identity on $\calG$. But by construction of $\calQ$, all the (iso)morphisms of $\calG$ are mapped to automorphisms, so that $qp$ must be trivial. 

Next, given a functor $G\colon\calA\ar\calX$ such that $G=G\circ\id_{\calA}$ is trivial, $G$ inverts all the arrows of $\calA$ and thus factors uniquely as a functor $H$ through the groupoid $\calG$ of fractions. 
\begin{diagram}[80]
\calA| \eeql | \calA | \Ear p | \calG | \Ear q | \calQ ||
    | | |     | \eseaR G |      \Sear H | \saR L ||
| | | | |  | \calX ||
\end{diagram}
But $H=H\circ\id_{\calG}$ is trivial because $G$ and thus $H$ map all morphisms of $\calA$, and thus also their inverses in $\calG$, to automorphisms in $\calX$. Therefore $H$ factors uniquely as a functor $L$ through the $\calZ$-cokernel $q$ of $\id_{\calG}$. The uniqueness is obvious since both $p$ and $q$ are epimorphisms. 
\qed

\begin{proposition}\label{zcoker}
The $\calZ$-cokernel of an arbitrary  functor $F\colon \calA \to \calB$ in $\Cat$ exists and is given by the pushout of the $\calZ$-cokernel of $\id_{\calA}$ along $F$.
\end{proposition}

\proof
This can easily be verified directly. With the notation above consider the pushout
\begin{diagram}[80]
\calA | \eeql | \calA | \Ear F | \calB ||
| | \Sepi {qp} | | \sepI U ||
| | \calQ | \eaR V | \calC ||
\end{diagram}
where $U$ is an epimorphism, since so is $qp$; this takes already care of the uniqueness condition.
If a functor $W\colon \calB \ar \calY$ is such that $WF=WF\id_{\calA}$ is trivial, then $WF$ factors uniquely  through the $\calZ$-cokernel $qp$ of $\id_{\calA}$ and thus further through the pushout. \qed

\begin{remark}\label{multipointed}
\emph{The argument used in the proof of Proposition~\ref{zcoker} is  a special instance of a known result. Indeed, from Proposition~\ref{zker} and Lemma~\ref{zcokerid}, it follows that $\Cat$ is a multipointed category (in the sense of Grandis and Janelidze \cite[Section~1.3]{GJ})) with respect to the full subcategory $\mathsf{SkGrpd}$ of skeletal groupoids. Accordingly, Proposition \ref{zcoker}  is then a consequence of their Proposition~1.3. It is worth noting that the notion of torsion theory in multipointed categories, as defined in \cite{GJ}, is a particular instance of the notion of pretorsion theory \cite{fafi,FFG}. In particular, in the ``multipointed context", the subcategory of trivial objects is required to be both epireflective and monocoreflective in the ambient category \cite[Section~1.5]{GJ}. This is the case for the pretorsion theory $(\mathsf{Grpd}, \mathsf{SkCat})$ studied in this paper, but not in general. Indeed, some examples of pretorsion theories that are not torsion theories in a multipointed category can be found in \cite[Remark~7.3]{FFG}.}
\end{remark}

\begin{corollary}
The constructions of the $\cal Z$-kernels and $\cal Z$-cokernels of the identity morphisms given in Proposition~\ref{zker} and Lemma~\ref{zcokerid} define the monocoreflection and the epireflection, respectively, of $\mathsf{Cat}$ to $\mathsf{SkGrpd}$.
\end{corollary}

\proof
This follows from Remark~\ref{multipointed} and \cite[Section~1.5]{GJ}.
\qed

\begin{remark}\label{stable}
\emph{It is possible to construct the ``universal coproduct-preserving stable category" associated with the pretorsion  theory described in Theorem~\ref{graut}. Indeed, this fact has been observed in Section~6 of \cite{bocagr3}.}
\end{remark}

\section{Appendix: Sketch of the proof of  \ref{reduced}}

This appendix presents a sketch of the proof of Proposition \ref{reduced}, which appears as Statement 1.1 in \cite{isbell} and whose full proof can be found under Satz 4.4 in \cite{borger}.

\begin{proposition}\label{proof}
Let $\calB$ be a small category and $M$ the free monoid on the set $\ArB$ of arrows of $\calB$. An element of $M$, that is a finite sequence of elements of $\ArB$, is {\em reduced} when it does not contain any identity morphism of $\calB$, nor any consecutive pair of arrows which are composable in $\calB$.\\
Consider the smallest congruence on $M$ identifying a sequence containing an identity with the sequence obtained when dropping that identity, and a sequence containing a pair of consecutive composable morphisms in $\calB$ with the sequence where this pair is replaced by the corresponding composite.\\
Then every element of $M$ is equivalent to a unique reduced element. 
\end{proposition}

\proof
As a matter of convention, let us use Latin letters for the arrows of $\calB$ and Greek letters for arbitrary elements of $M$. The existence of a reduced element obtained from $\alpha\in M$, in the way indicated above, is obvious. The point is to prove the uniqueness.

An element of $M$ is a finite sequence $\alpha$ of elements of $\ArB$; we shall often consider its length  $L(\alpha)\in\mathbb N$ . We shall also write $R\subseteq M$ for the set of reduced elements of $M$.

The first step of the proof is to observe, by induction on the length of $\beta$, that the following formulae define inductively a mapping
$$
\Phi\colon \ArB \times R \ar R,
$$
where $\Phi(a,\beta)$ is such that
\begin{enumerate}
\item $\Phi(a,\beta)=\Phi(ab,\nu)$  when $\beta=b\nu$, with  $\nu\in R$, while  $b\in \ArB$ is composable with $a$ in $\calB$;
\item $\Phi(a,\beta)=\beta$ when $a$ is an identity arrow in $\calB$ and  Case~1 does not apply;
\item $\Phi(a,\beta)=a\beta$ otherwise.
\end{enumerate}

The second step of the proof is, again by induction on the length of $\nu$, to define a binary operation
$$
\star \colon M\times R \ar R
$$
such that, when writing $()$ for the empty sequence (the unit of $M$), one has
$$
()\star\nu=\nu,~~a\beta\star\nu=\Phi(a,\beta\star\nu).
$$

One considers then the equivalence relation $S$ on $M$ defined by
$$
S=\{(\alpha,\beta)\mid \forall \nu\in R~~\alpha\star\nu=\beta\star\nu\}
$$
and  proves, by induction on the length of $\alpha$, that
$$
\alpha\star(\beta\star\nu)=\alpha\beta\star\nu,~~\alpha,\beta\in M,~\nu\in R. 
$$
This equality, together with the fact that $\gamma\star\nu\in R$ holds for every $\gamma \in R$, proves that $S$ is a congruence on $M$.
The equality implies also, via a new induction on the length of $\nu$ as appearing in the definition of $S$, that
$$
\bigl((a,b),ab\bigr)\in S,~~\bigl(\id_B,()\bigr)\in S,~~\mbox{for all}~a,b\mbox{ composable in }\calB, ~B\in\calB.
$$

By yet another induction on the length of $\nu$, one observes that $\nu\star ()=\nu$ holds for every $\nu\in R$. Then, given $\alpha\in M$ and $\nu,\eta\in R$, by definition of $S$, we obtain the implication
$$
(\alpha,\nu)\in S ~\mbox{and}~(\alpha,\eta)\in S ~\Longrightarrow ~\nu=\nu\star ()=\eta\star () = \eta.
$$
This proves the uniqueness condition in the statement.
\qed

\vspace{10mm}

\noindent
francis.borceux@uclouvain.be
\hfill tholen@york.ca\\
federico.campanini@uclouvain.be
\hfill Department of Mathematics\\
marino.gran@uclouvain.be
\hfill  and Statistics\\
Institut de Recherche 
\hfill York University\\
\hspace*{5mm}en Math\'ematique et Physique
\hfill Toronto ON, M3J 1P3\\
Universit\'e Catholique de Louvain
\hfill Canada\\
1348 Louvain-la-Neuve\\
Belgium

\end{document}